\pdfoutput=1
\documentclass[10pt]{article}
\usepackage{cancel}
\usepackage{graphicx}
\usepackage[colorlinks=true, pdfstartview=FitV, linkcolor=blue, citecolor=blue,urlcolor=blue]{hyperref}

\usepackage{framed,color}

\usepackage{amsbsy}
\usepackage{amssymb,amsmath}
 \usepackage{amsfonts}
\usepackage{graphicx}

\usepackage{cancel}
\usepackage{mathtools}
\usepackage{amsmath}
\usepackage{amsthm}
\usepackage{amsfonts}
\usepackage{amssymb}
\usepackage{mathrsfs}
\usepackage{graphicx}
\usepackage[T1]{fontenc}
\usepackage[OT2,T1]{fontenc}
\usepackage[latin1]{inputenc}
 \usepackage{lmodern}
\usepackage[all,cmtip]{xy}

\usepackage{bbm}
\usepackage[vcentermath]{youngtab}
\usepackage{ stmaryrd }
\usepackage{calc}

\usepackage{setspace}
\usepackage{epigraph}
\usepackage{mathdots}
\usepackage{mathabx}
\usepackage{braket}
\usepackage{longtable}
\usepackage{float}

\def\eqref#1{(\ref{#1})}
\newtheorem{theorem}{Theorem}[section]
\newtheorem{example}{Example}[section]
\newtheorem{exercise}{Exercise}[section]

\newtheorem{lemma}{Lemma}[section]
\newtheorem{remark}{Remark}[section]

\newtheorem{proposition}{Proposition}[section]
\newtheorem{corollary}{Corollary}[section]
\newtheorem{definition}{Definition}[section]

\def\bre{\begin{remark}}
\def\ere{\end{remark}}
\def\bth{\begin{theorem}}
\def\eth{\end{theorem}}
\def\bcr{\begin{corollary}}
\def\ecr{\end{corollary}}
\def\bex{\begin{example}\small}
\def\eex{\end{example}}
\def\bexr{\begin{exercise}\small}
\def\eexr{\end{exercise}}
\def\ble{\begin{lemma}}
\def\ele{\end{lemma}}

\def\bde{\begin{definition}}
\def\ede{\end{definition}}
\def\bpr{\begin{proposition}}
\def\epr{\end{proposition}}
\def\be{\begin{equation}}
\def\ee{\end{equation}}
\def\bea{\begin{eqnarray}}
\def\eea{\end{eqnarray}}
\def\beas{\begin{eqnarray*}}
\def\eeas{\end{eqnarray*}}

\def\vY{\vec{Y}}
\def\vPsi{\vec{\Psi}}
\def\vpsi{\vec{\psi}}

\def\lp{\lambda^\prime}

\newlanguage\fakelanguage
\newcommand\cyr{\fontencoding{OT2}\fontfamily{wncyr}\selectfont
   \language\fakelanguage}
\DeclareTextFontCommand{\textcyr}{\cyr}

\numberwithin{equation}{section}

\numberwithin{theorem}{section}

\numberwithin{proposition}{section}

\numberwithin{definition}{section}

\numberwithin{remark}{section}

\numberwithin{lemma}{section}

\numberwithin{corollary}{section}

\date{}
\begin{document}
\baselineskip=14pt

\vspace{0.2cm}
\begin{center}
\begin{Large}
\fontsize{17pt}{27pt}
\selectfont

\textbf{Isomonodromic Laplace Transform with 
Coalescing    Eigenvalues and 
 Confluence of Fuchsian Singularities}
\end{Large}
\\
\bigskip
\begin{large} {Davide  Guzzetti}\end{large}
\\{ SISSA,  Via Bonomea, 265,  34136 Trieste -- Italy.  E-MAIL: guzzetti@sissa.it}
\\
\bigskip
\end{center}
\begin{small}
{\bf Abstract:} 
We consider a Pfaffian system expressing isomonodromy of an irregular system of Okubo type, depending on complex deformation parameters $u=(u_1,...,u_n)$, which are eigenvalues of the leading matrix at the irregular singularity. At the same time, we consider a Pfaffian system of non-normalized Schlesinger type  expressing isomonodromy of a Fuchsian system, whose poles are the deformation parameters $u_1,...,u_n$. The parameters vary in a polydisc  containing   a  {\it coalescence locus for the eigenvalues} of the leading matrix of the irregular system,  corresponding  to  {\it confluence of the Fuchsian singularities}. We construct isomonodromic  {\it selected and singular vector solutions} of the Fuchsian Pfaffian system  together with their {\it isomonodromic connection coefficients}, so extending a result of \cite{BJL4} and \cite{guz2016} to the isomonodromic case, including confluence of singularities. Then, we introduce an isomonodromic  Laplace transform of the  selected and singular vector solutions, allowing to obtain  isomonodromic fundamental solutions for the irregular system,  and their  Stokes matrices expressed   in terms of connection coefficients. These facts,   in addition to extending  \cite{BJL4,guz2016} to the isomonodromic case (with coalescences/confluences),  allow to prove by means of Laplace transform the main result of \cite{CDG}, namely  the   analytic  theory of {\it non-generic  isomonodromic deformations} of the  irregular system  with coalescing eigenvalues. 
\vspace{0.7cm}

\end{small}
\vskip 15pt

\tableofcontents
\section{Introduction}

In this  paper, I answer  a  question asked when I presented the results of \cite{CDG} and the related paper \cite{Guzz-SIGMA}. Paper  \cite{CDG} deals with the extension of the theory of 
isomonodromic deformations of  the differential  system \eqref{24nov2018-1} below,  
 in presence of a coalescence phenomenon involving the eigenvalues of the leading matrix
  $\Lambda$. 
  These eigenvalues are the deformation parameters. The question is if we 
  can obtain some results of \cite{CDG}   in terms of the Laplace transform relating  
  system
   \eqref{24nov2018-1} to a Fuchsian one, such as system \eqref{03} below.
    The latter has 
   simple poles at the eigenvalues of $\Lambda$, so that the  coalescence of the eigenvalues will 
   correspond to the  confluence of the Fuchsian singularities. So the question is if combining 
  integrable deformations of Fuchsian systems, confluence of singularities and Laplace 
   transform, we can  obtain the results of \cite{CDG}.  
The positive   answer  is   {\bf Theorem  \ref{20agosto2020-5}} of this paper.   In order to achieve it, we   extend to the case depending on deformation parameters, including their coalescence,  one main result of \cite{BJL4} and \cite{guz2016}  concerning the existence of selected and singular vector solutions of a Pfaffian Fuchsian system associated with \eqref{03}  (see the system \eqref{19agosto2020-5} below), and their connection coefficients, which will be isomonodromic. This will be obtained  in  {\bf Theorem  \ref{30agosto2020-5}} and  {\bf Proposition  \ref{13ottobre2020-1}}.

\vskip 0.2 cm 
In \cite{CDG}, the isomonodromy deformation theory of an $n$ dimensional differential system with    Fuchsian singularity at $z=0$ and singularity of the second kind at  $z=\infty$ of Poincar\'e rank 1
\begin{equation}
\label{24nov2018-1}
\frac{dY}{dz}=\left(\Lambda(u)+\frac{A(u)}{z}\right)Y, \quad\quad \Lambda(u)=\hbox{\rm diag}(u_1,...,u_n),
\end{equation}
has been  considered\footnote{With the notation $\widehat{A}_1(u)$ for  $A(u)$.}. The deformation parameters  $u=(u_1,...,u_n)$ vary in a  polydisc where  the matrix $A(u)$ is {\it holomorphic}. One of the main results of \cite{CDG} is the extension of the theory of isomonodromic deformations of  \eqref{24nov2018-1} to the {\it non-generic case} when $\Lambda$  has coalescing 
eigenvalues but remains diagonalizable. This means that the polydisc contains a locus of {\it coalescence points} such that $u_i=u_j$ for some $1\leq i\neq j\leq n$.  In this 
case, $z=\infty$ is sometimes called {\it resonant irregular singularity}.   On a sufficiently small domain in the polydisc, the well know theory of isomonodromy deformations applies and allows to  define constant   monodromy data. 
Theorem 1.1 and corollary 1.1 of \cite{CDG} say that the these data are well defined and constant  on the whole polydisc, including the coalescence locus,  if the  entries of  $A(u)$ satisfy the {\it vanishing conditions} 
\be
\label{2aprile2021-5}
(A(u))_{ij}\to 0 \hbox{ when $u$ tends to a coalescence point such that    $u_i-u_j\to 0$ at this point}.
\ee
 More precisely, if conditions  \eqref{2aprile2021-5} are satisfied, the following results (reviewed in Theorem \ref{18agosto2020-8}  of Section \ref{23agosto2020-12} below)  hold.  
\begin{itemize}

\item[(I)]
 Fundamental matrix solutions in Levelt form at $z=0$ and solutions with prescribed ``canonical'' asymptotic behaviour  in Stokes sectors at $z=\infty$ are holomorphic of $u$ in the polydisc. Also the coefficients of the formal solution determining the asymptotics  at $\infty$ are holomorphic. 
 
 \item[(II)]     {\it Essential monodromy data}, such as Stokes matrices, the central connection matrix, the formal monodromy exponent at infinity and the Levelt exponents at $z=0$ are well defined and constant on the whole polydisc, including coalescence points. 
 
 The Stokes matrices (labelled by $\nu\in \mathbb{Z}$) satisfy the vanishing conditions $$(\mathbb{S}_\nu)_{ij}=(\mathbb{S}_\nu)_{ji}=0,  \hbox{ $i\neq j$,  if there is  a coalescence point  such that $u_i=u_j$.}
 $$
 
 \item[(III)] The  constant essential monodromy data  can be computed  from  the system ``frozen'' at a fixed coalescence point. In particular, if the constant diagonal entries of $A$  are partly non-resonant (see Corollary \ref{26nov2018-14}), then there is no ambiguity in this computation, being the formal solution unique. 
\end{itemize}

The results above have been established in \cite{CDG} by {\it direct} analysis of system  \eqref{24nov2018-1}, of its Stokes phenomenon and its isomonodromic deformations.

\bre{\rm   If $A(u)$ is holomorphic on the polydisc and \eqref{24nov2018-1} is an isomonodromic family  on the polydisc {\it minus} the coalescence locus (in the sense of integrability of an associated Pfaffian system \eqref{28nov2018-1} introduced later), then \eqref{2aprile2021-5}  are automatically satisfied and Theorem 1.1 of \cite{CDG} holds.  This is not mentioned  in \cite{CDG}.   I thank the referee for this observation. More details are in Remark \ref{6aprile2021-2}.
}\ere

For future use, we denote by $\lp_1,\dots,\lp_n$ the diagonal entries of $A(u)$, and 
$$
B:=\hbox{\rm diag}(A(u))=\hbox{\rm diag}(\lp_1,\dots,\lp_n).
$$
We will see that these $\lp_k$  are constant in the isomonodromic case.

\vskip 0.3 cm 
From another perspective,   if $u$ is {\it fixed} and $u_i\neq u_j$ for   $i\neq j$, namely  for a system \eqref{24nov2018-1} {\it not depending on parameters} with {\it pairwise distinct eigenvalues} of $\Lambda$, it is well known that columns of fundamental matrix solutions with prescribed asymptotics in Stokes sectors at $z=\infty$   can be obtained by Laplace-type integrals of certain selected column-vector solutions of an $n$-dimensional Fuchsian system of the type 
\be
\label{23agosto2020-24}
\frac{d\Psi}{ d\lambda}=\sum_{k=1}^n \frac{B_k }{ \lambda-u_k}\Psi,~~~~~B_k:=-E_k(A+I).
\ee
Here, $E_k$ is the elementary matrix whose entries are zero, except for $(E_k)_{kk}=1$. 
These facts are studied  in the seminal paper \cite{BJL4} in the generic case of non-integer  diagonal entries $\lp_k$ of $A$. The 
   results of \cite{BJL4} have been extended in \cite{guz2016} to the general case, when the entries $\lp_k$  take any complex value. 
 
\vskip 0.2 cm 
The purpose of the present paper is to introduce an {\it isomonodromic  Laplace transform}  relating \eqref{24nov2018-1} to an isomonodromic Fuchsian system
\be
\label{03}
\frac{d\Psi}{ d\lambda}=\sum_{k=1}^n \frac{B_k(u) }{ \lambda-u_k}\Psi,~~~~~B_k:=-E_k(A(u)+I).
\ee
when $u_1,...,u_n$ vary in a polydisc containing a locus of  coalescence points.  More precisely, the Laplace transform will relate  solutions of the integrable  Pfaffian systems \eqref{28nov2018-1} and   \eqref{19agosto2020-5} introduced later, associated with  \eqref{24nov2018-1} and \eqref{03} respectively. 
 The two main goals will   be: 
\begin{itemize}
\item {\bf Theorem  \ref{30agosto2020-5}}, which   characterises {\it selected vector solutions} and {\it singular vector solutions} of \eqref{03} and \eqref{19agosto2020-5}, so extending the results of \cite{BJL4} and \cite{guz2016} to the case depending on isomonodromic deformation parameters, including coalescing  Fuchsian singularities $u_1,...,u_n$. 

\item {\bf Theorem \ref{20agosto2020-5}}, in which  the Laplace transform of the vector solutions of Theorem  \ref{30agosto2020-5} allows  to obtain the main results (I), (II) and (III) of \cite{CDG} in presence of  coalescing eigenvalues $u_1,...,u_n$ of $\Lambda(u)$.

\end{itemize}
In details.

\vskip 0.2 cm

$\bullet$ In Proposition \ref{11agosto2020-5} we  establish the equivalence between the ``strong'' isomonodromic deformations (non-normalized Schlesinger deformations) of \eqref{03} and the ''strong'' isomonodromic deformations of \eqref{24nov2018-1}. 

$\bullet$ Then, we  study isomonodromy   deformations  of \eqref{03}  when $u$ varies in a polydisc containing a  coalescence locus. 
  {\bf Theorem \ref{30agosto2020-5}}, provides   selected  and singular vector solutions, which are the isomonodromic analogue of solutions  introduced in \cite{BJL4,guz2016}, respectively  denoted by $\vec{\Psi}_k(\lambda,u~ |\nu)$ and $\vec{\Psi}_k^{(sing)}(\lambda,u~|\nu)$, $k=1,...,n$, the latter being singular at $\lambda=u_k$. As will be explained later,   $\nu\in\mathbb{Z}$  labels the directions of branch cuts in   the punctured $\lambda$-plane  at the poles $u_1,...,u_n$. These solutions allow to introduce  {\it connection coefficients} $c_{jk}^{(\nu)}$, defined by 
$$ 
\vec{\Psi}_k(\lambda,u~ |\nu)=\vec{\Psi}^{(sing)}_j(\lambda,u~ |\nu) c_{jk}^{(\nu)} +\hbox{holomorphic part at $\lambda=u_j$},
\quad\forall~ j\neq k.
$$
The above is the deformation parameters dependent analogue of the definition of connection coefficients in \cite{guz2016}. 

$\bullet$ In  {\bf Proposition  \ref{13ottobre2020-1}}, we  prove that the  $c_{jk}^{(\nu)}$ are {\bf isomonodromic connection coefficients}, namely independent of  $u$.  When there is a coalescence $u_j=u_k$   in the polydisc, they satisfy $$
 c_{jk}^{(\nu)}=0,\quad j\neq k.
 $$

$\bullet$ In {\bf Theorem \ref{20agosto2020-5}}, the Laplace transform of the vectors $\vec{\Psi}_k(\lambda,u~ |\nu)$ or $\vec{\Psi}_k^{(sing)}(\lambda,u~ |\nu)$ yields the columns of the isomonodromic fundamental matrix solutions $Y_\nu(z,u)$ of \eqref{24nov2018-1}, labelled by $\nu\in\mathbb{Z}$,  uniquely determined by a prescribed asymptotic behaviour in certain $u$-independent  sectors $\widehat{\mathcal{S}}_\nu$, of central opening angle greater than $\pi$.  The analytic properties for the matrices $Y_\nu(z,u)$ will be proved, so re-obtaining the result (I) above. 
 In order to describe the Stokes phenomenon, only three solutions $Y_\nu(z,u)$, $Y_{\nu+\mu}(z,u)$ and $Y_{\nu+2\mu}(z,u)$  suffice. 
 The labelling will be  explained later. The  Stokes matrices $\mathbb{S}_{\nu+k\mu}$, $k=0,1$,  
 defined by a relation $Y_{\nu+(k+1)\mu}=Y_{\nu+k\mu} \mathbb{S}_{\nu+k\mu}$, will be expressed in terms of the coefficients $c_{jk}^{(\nu)}$ in formula \eqref{16ottobre2020-5}.  This  extends to the isomonodromic case, including coalescences, an analogous  expression appearing in \cite{BJL4,guz2016} and implies the results in  (II)  above.

$\bullet$ In Section \ref{16agosto2020-1}, we   re-obtain the result (III),  that system \eqref{24nov2018-1},  "frozen" by fixing $u$  equal to the most coalescence point  $u^c$ in the polydisc (see Section \ref{23agosto2020-12} for $u^c$), admits a unique formal solution if and only if the (constant) diagonal entries $\lambda_j^\prime$ of $A$  satisfy 
$\lambda_i^\prime-\lambda_j^\prime\not\in \mathbb{Z}\backslash\{0\}$ for every $i\neq j$ such that    $u^c_i=u^c_j$.   In this case we prove that  the  selected vector solutions of the Fuchsian system  \eqref{03} at $u=u^c$,  needed to perform the Laplace transforms, are uniquely determined. On the other hand, if  some $\lambda_i^\prime-\lambda_j^\prime\in \mathbb{Z}\backslash\{0\}$ corresponding to $u_i^c=u_j^c$, then   there is a family of solutions of the Fuchsian system  \eqref{03} at a coalescence point, depending on a finite number of parameters: this facts is responsible, through the Laplace transform, of the existence of a family of formal solutions at the coalescence point. 

\vskip 0.2 cm

In \cite{Dub1,Dub2}, B. Dubrovin related system \eqref{24nov2018-1} to an isomonodromic system of type \eqref{03}, in the specific case when they respectively yield the  flat sections of the deformed  connection of a semisimple Dubrovin-Frobenius manifold and the  flat sections
 of the intersection form (extended Gauss-Manin system). In \cite{Dub1,Dub2}, the solutions  of \eqref{24nov2018-1}  
 are expressed by Laplace transform of the isomonodromic system \eqref{03}, but the eigenvalues $u_1,...,u_n$ are assumed to be 
  pairwise distinct  in a sufficiently small domain (analogous to the polydisc $\mathbb{D}(u^0)$ to be introduced later). 
 Moreover, $A$ is   skew-symmetric, so its diagonal elements are zero ($A$ is denoted by  $V$ and $\Lambda$  by $U$  in \cite{Dub1,Dub2}).  By a Coxeter-type identity, the entries of the monodromy matrices for the selected solutions of  \eqref{03} (which 
 are part of the {\it monodromy of the Dubrovin-Frobenius manifold})  are expressed in terms of the entries of the Stokes matrices. See also \cite{Monica,Dub4}. 

In  proposition 2.5.1 of \cite{GGI}, the authors prove  (I)   when system \eqref{24nov2018-1} is associated with a  Dubrovin-Frobenius manifold with semisimple coalescence points, and $A$ is skew-symmetric (in \cite{GGI} the irregular singularity  is at $z=0$). 
 Their proof  contains the core idea that the analytic properties of a  solution $Y(z,u)$  in (I)  
 are obtainable, by a Laplace transform, from  the analytic properties of a fundamental 
 matrix solution $\Psi(\lambda,u)$ of the Fuchsian Pfaffian system associated
  with \eqref{03} (see their   lemma 2.5.3). The latter is a particular case of the  
 Pfaffian systems studied  in \cite{YT}. On the other hand,   the  analysis of selected 
 and singular vector solutions of the Fuchsian Pfaffian system, required in our paper to  
 cover all possible cases (all possible $A$), is not necessary in  \cite{GGI}, due to the skew-symmetry of $A$, and the specific form of their Pfaffian system (see their equation (2.5.2); their discussion  is equivalent  our case $\lp_j=-1$ for all $j=1,...,n$).  Moreover, points (II) and  (III) are not discussed in \cite{GGI}   by means of the  Laplace transform. 
 
 In the present paper, by an  isomonodromic Laplace transform, we prove (I), (II) and  (III), and at the same time we generalise the results of \cite{BJL4,guz2016} to the isomonodromic case with coalescences, with  no assumptions on the eigenvalues and the diagonal entries of $A$.  This  analytic construction, to the best of our knowledge, cannot be found in the literature. 

The approach of the present paper may also  be used  to extend  the results of \cite{Dub1,Dub2} described above, relating the deformed flat connection and the intersection form,  namely Stokes matrices and monodromy group of the Dubrovin-Frobenius manifold,    in case of semisimple coalescent Frobenius structures studied in \cite{CDG1,Eretico2,C-et-all,Cotti1}. 

For further comments and reference on the use of the Laplace transform, the confluence of singularities and related topics, see the introduction of \cite{guz2016} and \cite{BridgeTole,Loday,MazzoIrr,REMY, Scha0, Scha1,Scha2, Scha, Klimes1, Klimes2, Klimes3, HLR} . 

\vskip 0.2 cm 
 Stickily related to ours  are the important  results of \cite{sabbah}. In \cite{CDG} (and in the present paper by Laplace transform), we have answered the question if the integrable deformation  \eqref{28nov2018-1} of  system \eqref{24nov2018-1}  extends from a polydisc (or a small open set) not containing coalescence points to a wider domain intersecting (a stratum of) the coalescence locus, and we have characterized the  monodromy data. The converse question is answered in \cite{sabbah}, namley if  an integrable deformation \eqref{28nov2018-1}   of   $(\Lambda(u^c) +A(u^c)/z)dz$ exists and is unique, having formal normal form 
  $ d(z\Lambda(u))+B/z ~dz$, where $B$ is the diagonal of $A(u^c)$.  More broadly, the  question of    \cite{sabbah} is the existence and uniqueness of   integrable deformations of meromorphic connections on $\mathbb{P}^1$  with irregular singularity, when a prefixed restriction is given at a 
  single point  $t_o$ in the space of  deformation parameters $T$, allowed to be a degenerate point,  namely a   coalescence point in our case (in \cite{sabbah}, deformation parameters are called $t\in T$). One asks if a connection $\omega(z,t_0)$  given at $t_o\in T$ can be deformed to $\omega(z,t)$, and if this deformation is unique.\footnote{ The notation $\omega$ and $G$ is  not taken from \cite{sabbah}.}
 Concerning uniqueness,  for a fixed normal form $\omega_0(z,t)$, the problem is to classify isomorphism pairs $(\omega,G)$ consisting of an integrable connection $\omega(z,t)$ (with poles in $T\times \{z=0\}$,  being $z=0$ used in  \cite{sabbah}, while $z=\infty$ is used in our works) and  a formal gauge transformation $G(z,t)$  (formal in $z$ but holomorphic in $t$), transforming $\omega(z,t)$ to $\omega_0(z,t)$.  
 In a general context, a uniqueness theorem is proved in \cite{Tey18b}: two pairs are isomorphic 
 (meaning that the composition of a gauge of one pair with the inverse gauge of the other pair is 
 convergent w.r.t. $z$) if and only if their restriction to any specific value $t_o$ are isomorphic. 
 Thus, the $t$-extension of a pair in a neighbourhood of $t_0$ is unique up to isomorphism. 
 The proof in \cite{Tey18b} makes use of the results of  Kedlaya \cite{Ked10,Ked11} and Mochizuki 
 \cite{Moc09,Moc11a,Moc11b,Moc14}, which allows to blow up $T\times\{0\}$,  and of the higher dimensional 
 asymptotic analysis in poli-sectors   for the formal gauge transformations, that is Majima's asympotic analysis  \cite{Mag84} for Pfaffian systems with irregular singularitues. In  
 \cite{sabbah}, the uniqueness result is proved for a restricted class of integrable connections,
  in which our \eqref{28nov2018-1} is contained (with irregular singularity at $z=0$ instead of
   $\infty$). So, given a block-diagonal  normal form $\omega_0(z,t)$ and a pair  consisting 
   of $\omega(t_o,z)$ and a formal  gauge $G(t_o,z)$, it is proved that the pair can be 
   deformed (existence) in a unique way (uniqueness) to $\omega(z,t), G(z,t)$, such that 
   $G[\omega]=\omega_0$. The strategy is to use a sequence of Kedlaya-Mochizuki blow-ups 
   to obtain a good normal form   
   (see also \cite{sab93,sab00}).  Then, Majima results on asymptic analysis can be used and adapted. 
 In our specific case, theorem 4.9  of \cite{sabbah}  means the existence and uniqueness of  the integrable deformation \eqref{28nov2018-1} of $(\Lambda(u^c) +A(u^c)/z)dz$, formally equivalent to  $ d(z\Lambda(u))+B/z dz$. These facts generalize   results  of   Malgrange \cite{Malg3,Malg86}  for irregular singularities to the case of coalescence points.

 Theorem 4.9, obtained   in \cite{sabbah} in geometric terms,  has been successively  proved  in \cite{Cotti20} by  analytic methods. 
  In \cite{Cotti20}, the  integrable deformation is obtained from prefixed monodromy data  at a coalescence point, using  the analytic $L^p$ theory a Riemnann-Hilber boundary value problems. 
  Both  authors of \cite{sabbah} and \cite{Cotti20} apply their results to semisimpe  Dubrovin-Frobenius manifolds. In particular,  \cite{Cotti20} proves that any formal semisimple Frobenius manifold is the completion of a pointed germ of an analytic Dubrovin-Frobenius manifold. The result is extended to $F$-manifolds in the recent work \cite{Cotti21}.
  
  A geometric formulation of the Laplace transform we have used here,  together with a synthetic proof of part of Theorem 1.1 of \cite{CDG},  is the object of the  recent work \cite{sabbah2021}.

\subsection*{Acknowledgements}   I remember with gratitude  professor B. Dubrovin, who gave us the initial motivation to investigate the problem of non-generic isomonodromic deformations, contributing in our discussions with  insight, experience and enthusiasm. 
I would like to thank professor H. Iritani for kindly  drawing my attention to the proof of proposition 2.5.1 of \cite{GGI}, and professor  C. Sabbah for his clear explanation of \cite{sabbah,Tey18b}  given at a talk I attended,  which allowed me   to appreciate the importance of  the geometric approach.   I also thank doctor G. Cotti for several stimulating discussions in the background of this paper.  Finally, I  would like to thank very much the anonymous   referees for valuable comments, suggestions and corrections that have  improved the paper. 
 The author is a member of  the European Union's H2020 research and innovation programme under the Marie Sk\l{l}odowska-Curie grant No. 778010 {\it IPaDEGAN}.

\vskip 0.3 cm 
The corresponding author states that there is no conflict of interest.

\section{Review of Background Material}

This section contains known material to motivate and understand our paper. For $X$ a topological space, we denote by  $\mathcal{R}(X)$ its universal covering. For $\alpha<\beta\in \mathbb{R}$, a sector  is written as 
$$ 
S(\alpha, \beta):=\{z\in \mathcal{R}(\mathbb{C}\backslash\{0\}) \hbox{ such that } \alpha<\arg z <\beta\} .
$$

\subsection{Background 1: Isomonodromy Deformations of  \eqref{24nov2018-1} with coalescing eigenvalues.}
\label{23agosto2020-12}

We review some results of \cite{CDG,Guzz-SIGMA} (see also \cite{Eretico1,Alcala,guzz-Bed}). 
 Consider a differential system \eqref{24nov2018-1}  with an $n\times n$
 with matrix coefficient   $A(u)$ holomorphic in  a polydisc 
 \be
 \label{22novembre2020-5}
 \mathbb{D}(u^c):=\{u\in\mathbb{C}^n\hbox{ such that } \max_{1\leq j\leq n} |u_j-u_j^c|\leq \epsilon_0\},\quad \epsilon_0>0,
 \ee
centered at a {\it coalescence point} $u^c=(u_1^c,...,u_n^c)$, so called because 
 $$u_i^c=u_j^c\quad  \hbox{ for some } i\neq j.
$$ 
The eigenvalues of $\Lambda(u)$ coalesce at $u^c$ and  also  along the following  {\it coalescence locus} $$
 \Delta:=\mathbb{D}(u^c)\cap \Bigl(\bigcup_{i\neq j}\{u_i-u_j=0\}\Bigr),
 $$
 We assume that $\mathbb{D}(u^c)$ is sufficiently small so that $u^c$ is {\it the most coalescent point}. Namely, if $u_j^c\neq u_k^c$ for some $j\neq k$,  then $u_j\neq u_k$ for all  $u\in \mathbb{D}(u^c)$.  A more precise characterisation of the radius $\epsilon_0$ of the polydisc will be given in Section \ref{30agosto2020-1}. 
 For  $u^0\in \mathbb{D}(u^c)\backslash\Delta$, let  $$\mathbb{D}(u^0)\subset ( \mathbb{D}(u^c)\backslash\Delta)$$ be a (smaller) polydisc centered at $u^0$, not containing coalescence points. 
 
\subsubsection{Deformations in $\mathbb{D}(u^0)$}
\label{18agosto2020-6}

If $\mathbb{D}(u^0)$ is sufficiently small, 
 the isomonodromic theory of Jimbo, Miwa and Ueno \cite{JMU} assures that the essential  monodromy data of \eqref{24nov2018-1} (see Definition \ref{24agosto2020-1} below) are constant over $\mathbb{D}(u^0)$ and can be computed   fixing   $u=u^0$.

In order to give fundamental solutions with  ``canonical''  form  at $z=\infty$,   in $\mathcal{R}(\mathbb{C}\backslash\{0\})$  we introduce the Stokes rays of  $\Lambda(u^0)$, defined by  
 $$
\Re ((u_j^0-u_k^0)z)=0, \quad\Im((u_j^0-u_k^0)z)
<0,
\quad 
1\leq j\neq k \leq n.
$$   
 Let   
\be
\label{29agosto2020-2}
\arg z =\tau^{(0)}
\ee 
 be a direction  which does not coincide with any of the Stokes rays of $\Lambda(u^0)$, called {\it admissible at $u^0$}.  
Each  sector of amplitude $\pi$, whose boundaries are not Stokes rays of $\Lambda(u^0)$, contains a certain number   $\mu^{(0)}\geq 1$  of  Stokes rays of $\Lambda(u^0)$,  with  angular directions 
$$
\arg z ~=~ \tau_{0},\tau_{1},~...~,\tau_{\mu^{(0)}-1},\quad\hbox{ with } \quad  \tau_{0},<\tau_{1}<...<\tau_{\mu^{(0)}-1}
$$
that we decide to label from $0$ to $\mu^{(0)}-1$. They are {\it basic rays}, since they  generate all the Stokes rays in $\mathcal{R}(\mathbb{C}\backslash\{0\})$ associated with $\Lambda(u^0)$ by the formula $$ 
\arg z=\tau_{\nu}:=\tau_{\nu_0} +k\pi,\quad\quad \nu_0\in\{ 0,..., \mu^{(0)}-1\}, \quad \quad \nu=\nu_0+k\mu^{(0)},\quad k\in\mathbb{Z}. 
$$
The choice to label a specific Stokes ray with $0$, as $\tau_0$ above, is arbitrary, and  it induces the labelling  $\nu\in\mathbb{Z}$ for all other rays.  Suppose the labelling has been chosen. Then, for some $\nu\in\mathbb{Z}$, we have 
\be
\label{22novembre2020-1} 
\tau_\nu<\tau^{(0)}<\tau_{\nu+1}.
\ee
Equivalently, given  $\tau^{(0)}$, one can choose a $\nu$ and decide to call $\tau_\nu$ and $\tau_{\nu+1}$ the Stokes   rays  satisfying \eqref{22novembre2020-1}. This induces the labelling of all other rays (notice that $\mu^{(0)}$ is {\it not} a choice!).

\begin{figure}
\centerline{\includegraphics[width=0.6\textwidth]{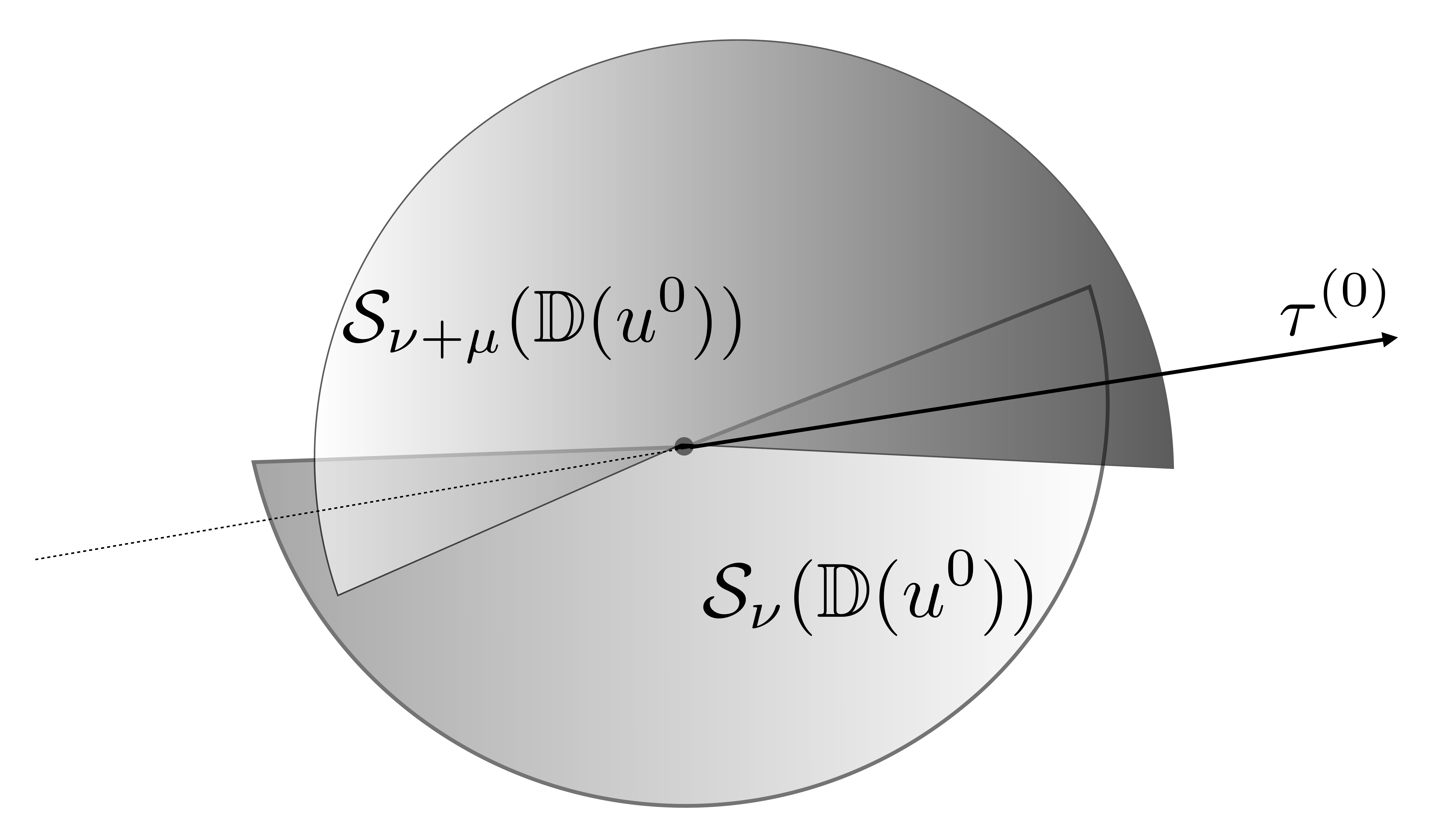}}
\caption{Successive sectors $\mathcal{S}_\nu(\mathbb{D}(u^0))$ and $\mathcal{S}_{\nu+\mu}(\mathbb{D}(u^0))$. Their intersection (in the right part of the figure) does not contain Stokes rays. It contains the admissible direction $\arg z=\tau^{(0)}$.}
\label{2settembre2020-12}
\end{figure}

Similarly, we consider the Stokes rays $
\Re ((u_j-u_k)z)=0$, $\Im((u_j-u_k)z)
<0$ of  $\Lambda(u)$. If $\mathbb{D}(u^0)$ is sufficiently small, when $u$ varies the Stokes rays of $\Lambda(u)$ 
rotate without crossing $\arg z=\tau^{(0)}$ mod $\pi$.  
For $k\in\mathbb{Z}$, we take the sector   $S\bigl(\tau^{(0)}+(k-1)\pi,\tau^{(0)}+k\pi\bigr)$  and  extend it in angular amplitude up to the nearest Stokes rays of $\Lambda(u)$ outside. The resulting (open) sector will be denoted by $\mathcal{S}_{\nu+k\mu^{(0)}}(u)$, and  we define  
$$
\mathcal{S}_{\nu+k\mu^{(0)}}( \mathbb{D}(u^0)):=
\bigcap_{u\in \mathbb{D}(u^0)}\mathcal{S}_{\nu+k\mu^{(0)}}(u).$$  
 The reason for the labelling is that $S\bigl(\tau^{(0)}+(k-1)\pi,\tau^{(0)}+k\pi\bigr)\subset S(\tau_{\nu+k\mu^{(0)}}-\pi,~\tau_{\nu+k\mu^{(0)}+1})$ and consequently 
$$
{\mathcal{S}}_{\nu+k\mu^{(0)}}(\mathbb{D}(u^0))\subset S(\tau_{\nu+k\mu^{(0)}}-\pi,~\tau_{\nu+k\mu^{(0)}+1})\equiv 
S(\tau_{[\nu+k\mu^{(0)}]-\mu^{(0)}},~\tau_{[\nu+k\mu^{(0)}]+1}).
$$
By construction, $\mathcal{S}_\nu( \mathbb{D}(u^0))$  has central angular opening greater than $\pi$. See figure \ref{2settembre2020-12}.

\begin{proposition}[Sibuya \cite{Sh2}, \cite{Sh4}, \cite{Sh11}; see also  \cite{JMU},  \cite{CDG}, \cite{Guzz-SIGMA}]
\label{25nov2018-6}
 Let $\mathbb{D}(u^0)$, not containing coalescence points,  be sufficiently small so that the  Stokes rays  of $\Lambda(u)$ do not cross\footnote{ As $u$ varies,  $\arg z=\tau^{(0)}+h_0\pi$, for a $h_0\in\mathbb{Z}$,  is not crossed by a Stokes ray    if and only if    $\arg z=\tau^{(0)}+h\pi$ is not crossed  $\forall ~h\in\mathbb{Z}$.}
  the admissible rays $\arg z=\tau^{(0)}+h\pi$, $h\in\mathbb{Z}$, as $u$ varies in $\mathbb{D}(u^0)$. System  \eqref{24nov2018-1} has a {\rm unique} formal solution 
 \be
\label{15agosto2020-9}
 Y_F(z,u)=F(z,u) z^{B(u)} \exp\{z\Lambda(u)\},\quad\quad B(u):=\hbox{\rm diag}(A_{11}(u),...,A_{nn}(u)),
\ee
   where 
  \be
  \label{15agosto2020-8}
  F(z,u)=I+\sum_{k=1}^\infty F_k(u) z^{-k}
  \ee
   is a formal series, with  holomorphic matrix coefficients  $F_k(u)$.For every $\nu\in\mathbb{Z}$,  there exist {\rm unique} fundamental matrix solutions
\begin{equation}
\label{25nov2018-4}
Y_\nu(z,u)=\widehat{Y}_\nu(z,u)z^{B(u)} \exp\{z\Lambda(u)\}
\end{equation}
 of \eqref{24nov2018-1}, holomorphic on $\mathcal{R}\bigl(\mathbb{C}\backslash\{0\}\times \mathbb{D}(u^0)\bigr)\equiv \mathcal{R}(\mathbb{C}\backslash\{0\})\times \mathbb{D}(u^0)$, such that uniformly in $u\in\mathbb{D}(u^0)$  the following asymptotic behaviour  holds
\begin{equation}
\label{26nov2018-4}
\widehat{Y}_\nu(z,u)\sim F(z,u)\quad  \hbox{ for  $z\to \infty$ in $\mathcal{S}_\nu( \mathbb{D}(u^0))$}.
\end{equation} 
\end{proposition}

The coefficients $F_k$ are computed recursively \cite{Wasow,CDG}
\begin{align}
\label{26nov2018-3}
&(F_1)_{i j}=\frac{A_{ij}}{u_j-u_i},~i\neq j, \quad\quad \quad (F_1)_{ii}=-\sum_{j\neq i}A_{ij}(F_1)_{ji},
\\
 \label{5aprile2018-1}
&(F_k)_{i j}=\frac{1}{u_j-u_i}\left\{
\Bigl(A_{i i}-A_{jj}+k-1
\Bigr)(F_{k-1})_{i j}
+\sum_{p\neq i}A_{i p}(F_{k-1})_{p j}\right\},
\quad
i\neq j;
\\
\label{5aprile2018-2}
&k(F_{k})_{i i}=-\sum_{j\neq i }A_{i j}(F_{k})_{ji}.
\end{align}
Holomorphic {\bf Stokes matrices} $\mathbb{S}_\nu(u)$, $\nu\in\mathbb{Z}$,  are the connection matrices  defined by 
\be
\label{23agosto2020-20}
Y_{\nu+\mu^{(0)}}(z,u)=Y_\nu(z,u) \mathbb{S}_\nu(u),\quad \quad z\in\mathcal{S}_\nu(\mathbb{D}(u^0))\cap \mathcal{S}_{\nu+\mu^{(0)}}(\mathbb{D}(u^0)).
\ee
Notice that $\mathcal{S}_\nu(\mathbb{D}(u^0))\cap \mathcal{S}_{\nu+\mu^{(0)}}(\mathbb{D}(u^0))$ does not contain Stokes rays of $\Lambda(u)$, for every $u\in\mathbb{D}(u^0)$. 

\vskip 0.2 cm 

At every fixed $u\in \mathbb{D}(u^0)$, system \eqref{24nov2018-1}
 admits a  fundamental matrix solution in {\it Levelt form} 
\begin{equation}
\label{25nov2018-7}
Y^{(0)}(z,u)=G^{(0)}(u)\Bigl(I+\sum_{j=1}^\infty \Psi_j(u)z^j\Bigr)z^Dz^L,
\end{equation}
where the series is convergent absolutely  in every ball $|z|<N$, for every $N>0$.  Here,  $D$ is diagonal with integer entries (called valuations), $L$ has eigenvalues  with  real part lying in $[0,1)$, and 
$D+\lim_{z\to 0}z^D L z^{-D}$ 
is  a  Jordan form of $A$. 
A {\bf central connection matrix} $C_\nu(u)$ is defined by 
\be
\label{23agosto2020-21}
Y_\nu(z,u)=Y^{(0)}(z,u)C_\nu(u).
\ee

A pair of  Stokes matrices $\mathbb{S}_\nu$, $\mathbb{S}_{\nu+\mu^{(0)}}$, together with $B$, $C_\nu$ and $L$  are sufficient to calculate all the other $\mathbb{S}_{\nu^\prime}$ and $C_{\nu^\prime}$, for all $\nu^\prime\in\mathbb{Z}$ (see \cite{BJL1,CDG}). 
The monodromy  matrices  at $z=0$  are    $$M:=e^{2\pi i L} \quad \hbox{ and } \quad  e^{2\pi i B} (\mathbb{S}_\nu \mathbb{S}_{\nu+\mu^{(0)}})^{-1}=C_\nu^{-1}MC_\nu
$$
  for  $Y^{(0)}$  and $Y_\nu$ respectively. Hence,  it makes sense to define {\it strong} isomonodromy deformations, as follows.

\begin{definition}
\label{24agosto2020-1}
Fixed a $\nu\in\mathbb{Z}$, we call {\bf essential monodromy data} the matrices 
  $$\mathbb{S}_\nu, \quad \mathbb{S}_{\nu+\mu^{(0)}}, \quad B, \quad C_\nu, \quad L, \quad D.
$$ The deformation $u$ is {\bf strongly isomonodromic} on $\mathbb{D}(u^0)$, if the essential monodromy data are constant  on $\mathbb{D}(u^0)$. 
\end{definition}
We introduced the terminology {\it strong} in \cite{Guzz-SIGMA}, to mean that  {\it all the  essential monodromy data} are constant, contrary to the case of {\it weak}  isomonodromic deformations, which {\it only preserve monodromy matrices} of a certain fundamental matrix solution. 
For a deformation to be {\bf weakly isomonodromic} it is  necessary and sufficient  that \eqref{24nov2018-1} is the $z$-component of a certain  Pfaffian system 
 $dY=\omega(z,u) Y$, Frobenius integrable (i.e. $d\omega=\omega\wedge \omega$). 
If $\omega$ is of very specific form, the defomation becomes strongly isomonodromic, according to the following
\bth
\label{9agosto2020-1}
System \eqref{24nov2018-1} is strongly isomonodromic in $\mathbb{D}(u^0)$ if and only  $Y_\nu(z,u)$, for every $\nu$,  and $Y^{(0)}(z,u)$, satisfy the Frobenius integrable Pfaffian system 
 \begin{equation}
 \label{28nov2018-1}
 dY=\omega(z,u) Y,\quad \quad \omega(z,u)=\left(\Lambda(u)+\frac{A(u)}{z}\right)dz+\sum_{k=1}^n
  \omega_k(z,u) du_k,
 \end{equation}
with the matrix  coefficients  (here $F_1$ is in \eqref{26nov2018-3})
\be
\label{9agosto2020-2}
\omega_k(z,u)= zE_k+ \omega_k(u),\quad\quad \omega_k(u)= [F_1(u),E_k].
\ee
Equivalently, \eqref{24nov2018-1} is strongly isomonodromic if and only if \footnote{Conditions \eqref{9agosto2020-2} and \eqref{9agosto2020-2-bis} imply Frobenius integrability of  \eqref{28nov2018-1}, so that the deformation is strongly isomonodromic. Conversely, given \eqref{28nov2018-1} with $ \omega_k(z,u)$ holomorphic in $\mathbb{C}\times \mathbb{D}(u^0)$, with $z=\infty$ at most a pole, then the integrability $d\omega(z,u)=\omega(z,u)\wedge \omega(z,u)$, which is necessary condition for isomonodromicity,  implies that $\omega_k(z,u)=zE_k+\omega_k(0,u)$ and \eqref{9agosto2020-2-bis}. Computations give that $\omega_k(0,u)=[F_1(u),E_k]+\mathcal{D}_k(u)$, where $\mathcal{D}_k(u)$ is an arbitrary diagonal holomorphic matrix.  Imposing that  $Y^{(0)}(z,u)$ and all the    $Y_\nu(z,u)$ satisfy \eqref{28nov2018-1}, then $\mathcal{D}_k(u)=0$ and $\omega_k(0,u)=[F_1(u),E_k]$.}  
  $A$ satisfies
 \be
 \label{9agosto2020-2-bis}
 dA=\sum_{j=1}^n \Bigl[\omega_k(u) ,A\Bigr]du_k .
 \ee
  If the deformation is strongly isomonodromic, then $Y^{(0)}(z,u)$ in \eqref{25nov2018-7}  is holomorphic on $\mathcal{R}(\mathbb{C}\backslash\{0\})\times \mathbb{D}(u^0)$, with holomorphic matrix coefficients   $\Psi_j(u)$,  and   
 the series is convergent  uniformly w.r.t. $u\in \mathbb{D}(u^0)$.  Moreover,  $G^{(0)}(u)$  is a holomorphic   fundamental solution of the integrable Pfaffian system 
 \begin{equation}
 \label{26nov2018-6}
 dG=\Bigl(\sum_{j=1}^n\omega_k(u)du_k\Bigr)G,
 \end{equation}
and  $A(u)$  is holomorphically similar to the Jordan form $J=G^{(0)}(u)^{-1}A(u)G^{(0)}(u)$.
\eth
The above theorem is analogous to the characterisation of isomonodromic deformations in \cite{JMU}, but includes also possible  resonances in $A$ (see \cite{CDG} and Appendix B of \cite{Guzz-SIGMA}).  
Notice that $\omega(z,u)$ in \eqref{28nov2018-1}-\eqref{9agosto2020-2} has components 
 \be
 \label{22agosto2020-6}
\omega_k(u)= \left(
 \frac{A_{ij}(\delta_{ik}-\delta_{jk})}{u_i-u_j}
 \right)_{i,j=1}^n
 =\left(
\begin{array}{ccccc}
0 & 0&\frac{ -A_{1k}}{u_1-u_k} &  0& 0
 \\
0 &0 &\vdots & 0& 0
 \\
 \frac{ A_{k1}}{u_k-u_1} & \cdots&0 & \cdots&\frac{ A_{kn}}{u_k-u_n}
 \\
0 &0 &\vdots &0  & 0
 \\
 0 &0 &\frac{ -A_{nk}}{u_n-u_k} &  0& 0
\end{array}
\right)
\ee
Notice that $B=\hbox{\rm diag}(A(u))=\hbox{\rm diag}(\lp_1,...,\lp_n)$ is constant because   \eqref{9agosto2020-2-bis} and \eqref{22agosto2020-6} imply
$$
\frac{\partial A_{ii}}{\partial u_j}=0, \quad \forall i,j=1,...,n.
$$ 

\subsubsection{Deformations in $\mathbb{D}(u^c)$ with coalescences}
\label{20agosto2020-15}
When the polydisc contains  a coalescence locus $\Delta$, the analysis presents  problematic issues.
\begin{itemize}
\item 
 A  fundamental matrix solution $Y(z,u)$ holomorphic on  $\mathcal{R}\bigl((\mathbb{C}\backslash\{0\})\times (\mathbb{D}(u^c))\backslash \Delta)\bigr)$, may be singular at $\Delta$, namely the limit for $u\to u^*\in \Delta$ along any direction may diverge, and  $\Delta$  is in general a {\it branching} locus \cite{Miwa}.  

\item The monodromy data associated with a fundamental matrix solution $\mathring{Y}(z)$ of
 \begin{equation}
\label{24nov2018-2}
\frac{dY}{dz}=\left(\Lambda(u^c)+\frac{A(u^c)}{z}\right)Y,
\end{equation} 
differ from those of any fundamental solution  $Y(z,u)$ of \eqref{24nov2018-1} at   $u\not \in\Delta$  (\cite{BJL2}, \cite{BJL3}, \cite{CDG}).  
\end{itemize}

In $\mathcal{R}(\mathbb{C}\backslash\{0\})$, we introduce the Stokes rays of $\Lambda(u^c)$
$$
\Re ((u_i^c-u_k^c)z)=0, \quad \Im((u_i^c-u_k^c)z)
<0,
\quad 
u_i\neq u_k,
$$
and an {\it admissible direction at $u^c$} 
\be
\label{20agosto2020-4}
\arg z=\tau,
\ee
 such that none of the  Stokes rays  at $u=u^c$ take this direction. Notice that $\tau$ is associated with $u^c$,  differently from $\tau^{(0)}$ of Section \ref{18agosto2020-6}. 
We choose $\mu$ basic Stokes rays of $\Lambda(u^c)$. These are all and the only Stokes rays lying in a  sector of amplitude $\pi$, whose boundaries are not Stokes rays of $\Lambda(u^c)$. Notice that $\mu$ is different from $\mu^{(0)}$ used in Section \ref{18agosto2020-6}.
 We   label their directions $\arg(z)$ as follows:
 $$\tau_0<\tau_1<...<\tau_{\mu-1}. $$
    The directions of all the other Stokes rays of $\Lambda(u^c)$  in $\mathcal{R}(\mathbb{C}\backslash\{0\})$ are consequently  labelled by an integer $\nu\in \mathbb{Z}$ 
\be
\label{15novembre2020-1}
\arg z=\tau_\nu:=\tau_{\nu_0}+k\pi,\quad\hbox{ with $\nu_0\in\{0,...,\mu-1\}$ and $\nu:=\nu_0+k\mu$}.
\ee
They satisfy $\tau_{\nu}<\tau_{\nu+1}$.

Analogously, at any  other $u\in \mathbb{D}(u^c)$, we define Stokes rays $\Re ((u_i-u_j)z)=0$, $\Im((u_i-u_j)z)<0$ of $\Lambda(u)$. 
They  behave differently from the case of $\mathbb{D}(u^0)$. Indeed, if $u$ varies in $\mathbb{D}(u^c)$,  some Stokes 
rays  cross   the admissible directions $\arg z= \tau$ mod $\pi$, as follows. Let $i,j,k$ be  such 
that  $u_i^c=u_j^c\neq u_k^c$. Then, as $u$ moves away from $u^c$, a Stokes ray of $\Lambda(u^c)$ 
 characterized by $\Re((u_i^c-u_k^c)z)=0$ generates three rays. Two of them are  $\Re((u_i-u_k)z)=0$ and $\Re((u_j-u_k)z)=0$. If $\mathbb{D}(u^c)$ is sufficiently small (as in \eqref{30agosto2020-6} below), they do not 
 cross $\arg z= \tau$ mod $\pi$ as $u$ varies in $\mathbb{D}(u^c)$. The third ray is    
 $\Re((u_i-u_j)z)=0$.   When $u$ varies in $\mathbb{D}(u^c)$  making a complete loop  $(u_i-u_j)\mapsto (u_i-u_j)e^{2\pi i}$   around the locus 
 $\{u\in \mathbb{D}(u^c)~|~u_i-u_j=0\}\subset \Delta$,  the third  ray crosses 
  $\arg z= \tau$ mod $2\pi$ and  $\arg z= \tau-\pi$ mod $2\pi$.  
This  identifies a  {\it crossing locus} 
$X(\tau)\subset \mathbb{D}(u^c)$  of points $u$ such that there exists a Stokes ray of $\Lambda(u)$ (so infinitely many in  $\mathcal{R}(\mathbb{C}\backslash\{0\})$)  with  direction  $\tau$ mod $\pi$.

\begin{proposition}[\cite{CDG}]
Each connected component of 
$\mathbb{D}(u^c)\backslash (\Delta\cup X(\tau))$  is  simply connected and homeomorphic to a ball, so it is a  topological cell.  
\end{proposition}

\noindent
Thus, the choice of $\tau$ induces a {\bf cell decomposition} of $\mathbb{D}(u^c)$. Each cell is called  {\bf $\boldsymbol{\tau}$-cell}. 
If $u$ varies in the interior of  a $\tau$-cell,  no Stokes rays cross the admissible directions $\arg z= \tau+h\pi$, $h\in\mathbb{Z}$, but if $u$ varies in the whole $\mathbb{D}(u^c)$, then  $X(\tau)$ is crossed, and thus  Proposition \ref{25nov2018-6} does not hold.

  To overcome this difficulty, we first take a point  $u^0$ in a $\tau$-cell, and consider a polydisc  $\mathbb{D}(u^0)$ {\it  contained  in the $\tau$-cell},    satisfying  the assumptions of sub-section \ref{18agosto2020-6}. Accordingly, we can define as before   the sectors  ${\mathcal{S}}_{\nu+k\mu}(u)$  of angular amplitude greater than $\pi$, and $$
{\mathcal{S}}_{\nu+k\mu}(\mathbb{D}(u^0))=
\bigcap_{u\in \mathbb{D}(u^0)}\mathcal{S}_{\nu+k\mu}(u)\subset \{\tau_{\nu+k\mu}-\pi<\arg z <\tau_{\nu+k\mu+1}\}.
$$
{Notice that here we are using $\tau$ and $\mu$ in place of $\tau^{(0)}$ and $\mu^{(0)}$. 
With the above sectors, monodromy data in  \eqref{23agosto2020-20}-\eqref{23agosto2020-21} can be defined  in $\mathbb{D}(u^0)$.  
   
    Since $A(u)$ is holomorphic in $\mathbb{D}(u^0)$, then $\omega_k(z,u)$    is holomorphic   on $\mathbb{D}(u^c)\backslash\Delta$. Thus, the fundamental matrix solutions $Y_\nu(z,u)$, $Y^{(0)}(z,u)$ of sub-section \ref{18agosto2020-6} extend analytically on $\mathcal{R}\bigl((\mathbb{C}\backslash\{0\})\times (\mathbb{D}(u^c))\backslash \Delta)\bigr)\neq \mathcal{R}(\mathbb{C}_z\backslash\{0\})\times(\mathbb{D}(u^c))\backslash \Delta)$, and  $\Delta$ may be a branching locus for them. 
 
 \bpr[\cite{CDG}]
 \label{6aprile2021-3}
 \label{6aprile2021-1}  $\omega(z,u)$  in \eqref{9agosto2020-2} and  \eqref{22agosto2020-6} is holomorphic on the whole  $\mathbb{D}(u^c)$ if and only if $$
A_{ij}(u)=\mathcal{O}(u_i-u_j)\to 0 \quad\hbox{ whenever $(u_i-u_j)\to 0$ for $u$ approaching $\Delta$}.
$$
$A(u)$ is holomorphically similar on $\mathbb{D}(u^c)$  to a Jordan form $J$  if and only if the above  vanishing conditions hold.  Similarity is realized by a fundamental matrix solution of \eqref{26nov2018-6}, which exists holomorphic on the whole $\mathbb{D}(u^c)$. 
\epr
 
 The extension of the theory of isomonodromy deformations on the whole  $\mathbb{D}(u^c)$  is given in \cite{CDG} by the following theorem, which is a detailed exposition of the points (I) and  (II)  of the Introduction, while point (III) is expressed by Corollary \ref{26nov2018-14} below. 
 
 \bth[\cite{CDG}]
\label{18agosto2020-8}
Let $A(u)$ be holomorphic on  $\mathbb{D}(u^c)$.  Assume that system \eqref{24nov2018-1} is strongly isomonodromic on $\mathbb{D}(u^0)$ contained in a $\tau$-cell of  $\mathbb{D}(u^c)$, so that Theorem \ref{9agosto2020-1} holds.  Moreover, assume that $A$ satisfies the vanishing conditions 
\be 
\label{26nov2018-7}
A_{ij}(u)=\mathcal{O}(u_i-u_j)\to 0 \quad\hbox{ whenever $(u_i-u_j)\to 0$ for $u$ approaching $\Delta$}.
\ee
Then, the following statements hold. 

\hfill \break
{\bf Part I.} 

\begin{itemize} 

\item[(I,1)] $Y^{(0)}(z,u)$ and the $Y_\nu(z,u)$, $\nu\in\mathbb{Z}$ admit analytic continuation as holomorphic functions on $\mathcal{R}(\mathbb{C}\backslash \{0\})\times \mathbb{D}(u^c)$. 
The coalescence locus $\Delta$ is neither  a singularity locus nor a  branching locus.

\item[(I,2)]  The coefficients $F_k(u)$ of $ Y_F(z,u)$, given in \eqref{26nov2018-3}-\eqref{5aprile2018-1}-\eqref{5aprile2018-2}, are  holomorphic of $u\in \mathbb{D}(u^c)$. 

\item[(I,3)]  The fundamental matrix solutions $Y_\nu(z,u)$  have asymptotics $Y_\nu(z,u) \sim Y_F(z,u)$ uniformly in $u\in\mathbb{D}(u^c)$, for $z\to \infty$ in a wide sector $\widehat{\mathcal{S}}_{\nu}$ containing  $  \mathcal{S}_{\nu}(\mathbb{D}(u^0))$, to be defined later in \eqref{16ottobre2020-1}.

\end{itemize}

\hfill \break
{\bf Part II.}  

\begin{itemize}
\item[(II,1)]
 the essential monodromy data $\mathbb{S}_\nu$, $\mathbb{S}_{\nu+\mu}$, $B=\hbox{\rm diag}(A(u^c))$, $C_\nu$, $L$, $D$, initially defined on $\mathbb{D}(u^0)$ by relations \eqref{23agosto2020-20}-\eqref{23agosto2020-21},  are well defined and constant on the whole $\mathbb{D}(u^c)$. They satisfy
$$ 
\mathbb{S}_\nu=\mathring{\mathbb{S}}_\nu,\quad \mathbb{S}_{\nu+\mu}=\mathring{\mathbb{S}}_{\nu+\mu}, \quad L=\mathring{L}, \quad C_\nu=\mathring{C}_\nu, \quad D=\mathring{D},
$$
where 
\item[(II,2)]
$\mathring{\mathbb{S}}_\nu$, $\mathring{\mathbb{S}}_{\nu+\mu}$  are the Stokes matrices of fundamental solutions $\mathring{Y}_\nu(z)$, $\mathring{Y}_{\nu+\mu}(z)$, $\mathring{Y}_{\nu+2\mu}(z)$ of  \eqref{24nov2018-2}  having asymptotic behaviour $\mathring{Y}_F(z)=Y_F(z,u^c)$, for $z\to \infty$ respectively on sectors 
$\tau_\nu-\pi <\arg z< \tau_{\nu+1}$, $\tau_\nu<\arg z < \tau_{\nu+\mu+1}$ and $\tau_{\nu+\mu}<\arg z<\tau_{\nu+2\mu+1}$;
\item[(II,3)]
$\mathring{L}$, $\mathring{D}$ are  the exponents of a fundamental solution 
$\mathring{Y}(z)=\mathring{G}\left(I+\sum_{j=1}^\infty \mathring{\Psi}_jz^j\right)z^{\mathring{D}}z^{\mathring{L}}$  of  \eqref{24nov2018-2}   in Levelt form;

\item[(II,4)] $\mathring{C}_\nu$ connects $\mathring{Y}_\nu(z)=\mathring{Y}(z)\mathring{C}_\nu$. 

\item[(II,5)] The Stokes matrices satisfy the vanishing conditions
$$ 
(\mathbb{S}_\nu)_{ij}=(\mathbb{S}_\nu)_{ji}=0,\quad (\mathbb{S}_{\nu+\mu})_{ij}=(\mathbb{S}_{\nu+\mu})_{ji}=0    \quad \forall~1\leq i\neq j\leq n \hbox{ such that } u_i^c=u_j^c.
$$
\end{itemize}

\eth

\begin{corollary}[\cite{CDG}]
\label{26nov2018-14}
If $A_{ii}-A_{jj}\not\in \mathbb{Z}\backslash\{0\}$ for every $i\neq j$ such that $u_i^c=u_j^c$, then   the formal solution $\mathring{Y}_F(z)$ of \eqref{24nov2018-2} is unique and coincides with $Y_F(z,u^c)$.

\end{corollary}

The assumption of Corollary \ref{26nov2018-14} will be called {\bf partial non-resonance}. If it holds, 
  {\it (II,1)}  says that  in order to obtain the essential monodromy data  of  \eqref{24nov2018-1} it suffices to compute  $ \mathring{\mathbb{S}}_\nu$, $\mathring{\mathbb{S}}_{\nu+\mu}$, $\mathring{L}$,  $\mathring{C}_\nu$ and $\mathring{D}$  for \eqref{24nov2018-2}, which is simper than  \eqref{24nov2018-1}, because $A_{ij}(u^c)=0$ for $i,j$ such that  $u_i^c=u_j^c$ . This  allows in some cases the explicit  computation of  monodromy data. An important example with algebro-geometric implications can be found in \cite{CDG1}.

 \bre
 \label{6aprile2021-2}
 {\rm
The following statement, not mentioned in \cite{CDG}, holds.
 
 {\it If \eqref{24nov2018-1}  is an isomonodromic family on the polydisc minus the coalescence locus,   in the sense that 
$dY=\omega Y$  in \eqref{28nov2018-1}-\eqref{9agosto2020-2} 
is Frobenius integrable on $\mathbb{D}(u^c)\backslash \Delta$, and if  $A(u)$ is holomorphic on  $\mathbb{D}(u^c)$, then the vanishing conditions \eqref{26nov2018-7} hold automatically and  \eqref{24nov2018-1}  is isomonodromic on $\mathbb{D}(u^c)$ in the strong sense, namely Theorem \ref{18agosto2020-8} holds.}

\vskip 0.1 cm 
I thank the referee for suggesting to write  the above statement. 
The sketch of the proof is as follows:   integrability  $d\omega=\omega\wedge \omega$ on $\mathbb{D}(u^c)\backslash \Delta$ 
  implies \eqref{9agosto2020-2-bis}, namely 
\be
 \label{2aprile2021-1}
 \frac{\partial A}{\partial u_j}= [\omega_j(u),A],\quad j=1,...,n;\quad u\in \mathbb{D}(u^c)\backslash \Delta.
\ee
 We want to prove that $A_{ij}(u)\to 0$ for $u_i-u_j\to 0$, for $i\neq j$. From \eqref{2aprile2021-1} and \eqref{22agosto2020-6}  we explicitly obtain  
 $$  \frac{\partial A_{i\ell}}{\partial u_j} = \frac{ (u_i-u_\ell)A_{ij}A_{j\ell} }{(u_i-u_j)(u_\ell -u_j)}, \quad  \hbox{ for $ j\neq i, \ell $ and $i\neq \ell$}.
 $$
  The left hand-side is holomorphic everywhere on $\mathbb{D}(u^c)$ by assumption on $A$, and so must  be the right hand-side. This implies that holomorphically  $ A_{ij} =O(u_i-u_j)$ for $u_i-u_j\to 0$. 
  Then, Theorem \ref{18agosto2020-8} holds and we conclude. $\Box$

 }
 \ere

The  difficulty  in  proving Theorem \ref{18agosto2020-8}  is  the analysis of the Stokes phenomenon at $z=\infty$. 
 On the other hand, coalescences does not affect the analysis at the Fuchsian singularity $z=0$,  so it is  not an issue for the proof of the statements concerning  $Y^{(0)}(z,u)$,    $L$ , $D$ and $C_\nu$ (as far as the contribution of  $Y^{(0)}$ is concerned). See  Proposition 17.1 of \cite{CDG}, and   the proof of  Theorem 4.9 in \cite{Guzz-SIGMA}. For this reason, in the present paper we will not deal with  $Y^{(0)}(z,u)$,    $L$ , $D$, $C_\nu$ and  {\it (II,3)-(II,4)} above.

In Theorem \ref{20agosto2020-5} we  introduce  {\it an isomonodromic  Laplace transform in order to prove  the statements of  Theorem \ref{18agosto2020-8} above,  concerning the Stokes phenomenon}, namely {\it (I,1), (I,2), (I,3)}
  and {\it (II,1), (II,2),  (II,5)}.

  \subsection{Background 2: Laplace Transform, Connection Coefficients and Stokes Matrices}
\label{19agosto2020-3}

 In this section, we {\it fix} $u\in\mathbb{D}(u^c)\backslash\Delta$.  Accordingly, system \eqref{24nov2018-1} is to be considered as a system {\it not depending on deformation parameters}, with leading matrix $\Lambda$ having {\it pairwise distinct eigenvalues}, and system \eqref{03} is equivalent to  \eqref{23agosto2020-24}, which does not  depend on parameters. For simplicity of notations, let us  fix for example 
$$u=u^{0},\quad \hbox{ as in Section \ref{18agosto2020-6}}.$$

 \vskip 0.2 cm 
 
 Solutions  $Y_\nu(z)$ of  \eqref{24nov2018-1}  with canonical asymptotics $Y_F(z)$ ($u=u^{0}$ fixed is not indicated) can be expressed in terms of convergent Laplace-type integrals \cite{Bi,Ince}, where the integrands are solutions of the Fuchsian system\footnote{The notation $A_0$ and $A_1$ is used in \cite{guz2016} for $\Lambda$ and $A$. In \cite{BJL4} the notation for $\Lambda$ is the same, while $A$ is denoted by $A_1$. The notation $\lambda_1,...,\lambda_n$ is used in \cite{BJL4,Guzz-SIGMA} for $u_1,...,u_n$. There is a misprint in the first page of \cite{guz2016} where it is said that $A_1\in GL(n,\mathbb{C})$; the correct statement is  $A_1\in Mat(n,\mathbb{C})$.}
\be
\label{02}
(\Lambda-\lambda)\frac{d \Psi}{ d\lambda}=  (A+I)\Psi,~~~~~I:=\hbox{ identity matrix}
\ee
 Indeed, let $\vPsi(\lambda)$ be a vector valued function and define 
$$
\vY(z)=\int_\gamma e^{\lambda z}\vPsi(\lambda)d\lambda, 
$$
where $\gamma$ is a suitable path. Then, substituting into \eqref{24nov2018-1}, we have 
$$
(z \Lambda+A)\int_\gamma e^{\lambda z}\vPsi(\lambda)d\lambda=z\frac{d}{ dz}\int_\gamma e^{\lambda z}\vPsi(\lambda)d\lambda
=
z\int_\gamma \lambda e^{\lambda z}\vPsi(\lambda)d\lambda
.$$
This implies that
$$
A \int_\gamma e^{\lambda z}\vPsi(\lambda)d\lambda= \int_\gamma \frac{d(e^{\lambda z}) }{ d\lambda} ~(\lambda-\Lambda)\vPsi(\lambda)d\lambda=
$$
\be
\label{31agosto2020-2}
=e^{\lambda z}(\lambda-\Lambda)\vPsi(\lambda)\Bigl|_\gamma-\int_\gamma
e^{\lambda z} \left[
(\lambda-\Lambda)\frac{d\vPsi(\lambda)}{ d\lambda}+\vPsi(\lambda)
\right]d\lambda.
\ee
If $\gamma$ is such that $e^{\lambda z}(\lambda-\Lambda)\vPsi(\lambda)\Bigl|_\gamma=0$, and if the function $\vPsi(\lambda)$ solves (\ref{02}), then $\vY(z)$ solves \eqref{24nov2018-1}.

\vskip 0.3 cm 
 Multiplying to the left by $(\Lambda-\lambda)^{-1}$,  system (\ref{02}) becomes  \eqref{23agosto2020-24},
\be
\label{29agosto2020-1}
\frac{d\Psi}{ d\lambda}=\sum_{k=1}^n \frac{B_k }{ \lambda-u_k^0}\Psi,~~~~~B_k:=-E_k(A+I).
\ee 
In order to define matrix solutions of  of \eqref{29agosto2020-1} as single valued functions, we consider the $\lambda$-plane with branch-cuts. Let  $\eta^{(0)}\in\mathbb{R}$  satisfy   
 \be
\label{30agosto2020-7}
\eta^{(0)}\neq \arg(u_j^0-u_k^0) \hbox{ mod } \pi,~~~\hbox{ for all }1\leq j,k\leq n.
\ee
We fix  parallel  cuts  $L_k(\eta^{(0)})$,  namely half-lines oriented  from $u_k^0$ to $\infty$ in direction $\arg(\lambda-u_k^0)=\eta^{(0)}$, $1\leq k \leq n$. See figure \ref{cutss}.   
Conditions  \eqref{30agosto2020-7} mean that a cut $L_k$ does not contain another pole $u_j^0$, $j\neq k$. For this reason $\eta^{(0)}$ is called  {\it admissible direction} at $u^0$.  Then,  we choose  a branch of the logarithms $ \ln(\lambda- u_k^0)=\ln|\lambda- u_k^0|+i\arg(\lambda-u_k^0)$   by 
 \be
\label{30agosto2020-7-bis}
 \eta^{(0)}-2\pi <\arg(\lambda-u_k^0)<\eta^{(0)},\quad \quad k=1,...,n.
\ee
 Following \cite{BJL4}, the $\lambda$-plane with these cuts and choices of the logarithms is denoted by $\mathcal{P}_{\eta^{(0)}}$.  
 Matrix solutions of \eqref{29agosto2020-1} are well defined as single-valued functions   of $\lambda\in \mathcal{P}_{\eta^{(0)}}$. 

\bre 
\label{3aprile2021-1}
{\rm  $\mathcal{P}_{\eta^{(0)}}$ can be identified with  one of the countably many components of    $$
\mathcal{R}^\prime:=\mathcal{R}( \mathbb{C}\backslash\{u_1^0,...,u_n^0\})-(\hbox{lift of all half-lines $L_k$}).$$ Each component is obtained by a deck transformation starting from  one. Fix  one component, for example  $\mathcal{P}_{\eta^{(0)}}$,   and define $2n$ letters 
$$ 
l_k:=\hbox{ cross a lift of $L_k$ from the right}, \quad l_k^{-1}:=\hbox{ cross a lift of $L_k$ from the left},\quad k=1,...,n,
$$
where ``right'' or ``left'' refers to the orientation of $L_k$.  The other components are reached by crossing the cuts, so that  there is a one to one correspondence between finite sequences $\{l_{j_1}^{\pm 1},...,l_{j_m}^{\pm 1}\}$ not containing successively a $l_k^{\pm 1}$ and $l_k^{\mp 1}$, and  components of $\mathcal{R}^\prime$ (here  $j_1,...,j_m\in \{1,...,n\}$ and $m\in\mathbb{N}$).  The relations \eqref{30agosto2020-7-bis}  alone do not identify a component of $\mathcal{R}^\prime$ (as incorrectly written in \cite{guz2016}, page 387). For example, the word $l_1l_2l_1^{-1}l_2^{-1}$ leads to a new component of $\mathcal{R}^\prime$ where   $\arg(\lambda-u_1^0),...,\arg(\lambda-u_n^0)$ take the same values of the starting component.\footnote{As well known, the analytic continuation,  starting from the plane $\mathcal{P}_{\eta^{(0)}}$,  of a fundamental matrix solution  of \eqref{29agosto2020-1} defines a function  $\Psi$ on   $\mathcal{R}(\mathbb{C}\backslash\{u_1^0,...,u_n^0\})$. For example, if $\lambda^\prime $ is the lift of  $\lambda\in \mathcal{P}_{\eta^{(0)}}$ to the component of $\mathcal{R}^\prime$ identified by  the word $l_1l_2l_1^{-1}l_2^{-1}$, then   $\Psi(\lambda^\prime)=\Psi(\lambda)M_2^{-1}M_1^{-1}M_2M_1$, where  $M_j$ is the  monodromy matrix associated with $l_j$.} I thank the referee for  this remark.  
}
 \ere

 Stokes matrices for  \eqref{24nov2018-1}, for fixed and pairwise distinct $u_1^0,...,u_n^0$, can been  expressed in terms of connection coefficients  of selected  solutions of \eqref{29agosto2020-1}.  The explicit relations have been obtained 
 in  \cite{BJL4} for the generic case when all $\lp_1,...,\lp_n \not \in \mathbb{Z}$; and  in \cite{guz2016} for  the general   case with  no restrictions on $\lp_1,...,\lp_n$ and  $A$.

\subsubsection*{Selected Vector Solutions}

The Laplace transform involves  three types of vector solutions of \eqref{29agosto2020-1}, denoted in  \cite{guz2016} respectively by $\vec{\Psi}_k(\lambda)$,  $\vec{\Psi}_k^{*}(\lambda)$ and  $\vec{\Psi}_k^{(sing)}(\lambda)$ , for $k=1,...,n$ (in \cite{BJL4} the notation used is $Y_k$ and $Y_k^*$, while  $Y_k^{(sing)}$  does not appear, since it reduces to $Y_k$  in the generic case $\lp_k\not\in\mathbb{Z}$). We will not describe here the $\vec{\Psi}_k^{*}(\lambda)$, which play mostly a technical role. Let 
 \begin{align*}
 &\mathbb{N}=\{0,1,2,...\} \hbox{ integers},\quad \mathbb{Z}_{-}=
\{-1,-2,-3,...\}\hbox{ negative integers}, 
\\
& \vec{e}_k = \hbox{ standard $k$-th unit column vector in $\mathbb{C}^n$}.
\end{align*}
 
It is proved in \cite{guz2016} that for every $k\in\{1,...,n\}$  there are at least $n-1$  independent vector solutions  holomorphic at $\lambda=u_k^0$. The remaining independent solution is singular at $\lambda=u_k^0$, except for some exceptional  cases possibly occurring when $\lp_k\leq -2$ is integer. In such cases, there exist $n$  holomorphic solutions at $\lambda=u_k^0$ (such cases never occur if none of the eigenvalues of $A$ is a negative integer).  The selected vector solutions  $\vec{\Psi}_k$ are obtained  as follows. 
 
 \begin{itemize}
  \item If $\lp_k\leq -2$ is integer and we are in an exceptional case when there are no singular solutions at $u_k^0$, namely 
  $$ \vPsi_k^{(sing)}(\lambda)\equiv 0, 
  $$
   then $\vec{\Psi}_k$ is the unique analytic solution with the following normalization:
  $$ 
  \vPsi_k(\lambda) =\left(\frac{(-1)^{\lp_k}}{ (-\lp_k-1)!}\vec{e}_k+\sum_{l\geq 1} \vec{b}_l^{~(k)}(\lambda-u_k^0)^l\right)(\lambda-u_k^0)^{-\lp_k-1}.
  $$
  
  \item In all other cases, there is a solution $\vec{\Psi}_k^{(sing)}$ with  singular behaviour at $\lambda=u_k^0$. This is determined up to a multiplicative factor and the addition of an arbitrary linear combination of the remaining $n-1$ regular at $\lambda=u_k^0$  solutions, denoted below with $\hbox{reg}(\lambda-u_k^0)$. In \cite{guz2016},  it has the following structure   
  \be 
  \label{20agosto2020-16}
  \vPsi_k^{(sing)}(\lambda)=\left\{
\begin{array}{cc}
\vpsi_k(\lambda)(\lambda-u_k^0)^{-\lp_k-1}+\hbox{reg}(\lambda-u_k^0) ,& \lp_k\not\in\mathbb{Z},\\
\\
\vpsi_k(\lambda)\ln(\lambda-u_k^0)+\hbox{reg}(\lambda-u_k^0), & \lp_k\in \mathbb{Z}_{-},\\
\\
\dfrac{ P_k(\lambda)
}{ (\lambda-u_k^0)^{\lp_k+1}}
+\vpsi_k(\lambda)\ln(\lambda-u_k^0)+\hbox{reg}(\lambda-u_k^0),& \lp_k\in\mathbb{N} .
\end{array}
\right.
\ee
Here  $\vpsi_k(\lambda)$ is analytic at $u_k^0$ and $P_k(\lambda)=\sum_{l=0}^{\lp_k} \vec{b}_l^{~(k)} (\lambda-u_k^0)^l$ is a polynomial of degree $\lp_k$. 
We choose the following  {\it normalization} at $\lambda=u_k^0$
 $$  
 \left\{
\begin{array}{cc}
\vpsi_k(\lambda)=\Gamma(\lp_k+1)\vec{e}_k+\sum_{l\geq 1} \vec{b}_l^{~(k)}(\lambda-u_k^0)^l,& \lp_k\not\in\mathbb{Z},\\
\\
\vpsi_k(\lambda)= \left(\dfrac{(-1)^{\lp_k}}{ (-\lp_k-1)!}\vec{e}_k+\sum_{l\geq 1} \vec{b}_l^{~(k)}(\lambda-u_k^0)^l\right)(\lambda-u_k^0)^{-\lp_k-1}& \lp_k\in \mathbb{Z}_{-},\\
\\
P_k(\lambda)=\lp_k!~\vec{e}_k +O(\lambda-u_k^0) & \lp_k\in\mathbb{N} ,
\end{array}
\right.
$$
The coefficients $ \vec{b}_l^{~(k)}\in\mathbb{C}^n$ are uniquely determined by the normalization.  Then the {\it selected vector solutions}  $\vPsi_k$ are {\it uniquely} defined by\footnote{The singular part of  $\Psi^{(sing)}$  is uniquely determined by the normalization, but not   $\Psi^{(sing)}$ itself, because  the analytic additive term  $\hbox{reg}(\lambda-u_k^0)$ is an arbitrary  linear combination of the remaining $n-1$ independent analytic solutions.}
\be
\label{20agosto2020-19}
\vPsi_k(\lambda):=\vpsi_k(\lambda)(\lambda-u_k^0)^{-\lp_k-1} \quad \hbox{for } \lp_k\not\in\mathbb{Z};\quad\quad 
 \vPsi_k(\lambda):=\vpsi_k(\lambda) \quad  \hbox{for  }  \lp_k\in\mathbb{Z}.
\ee
In case $\lp_k\in\mathbb{N}$, depending on the system, it may exceptionally happen that $$
\vPsi_k:=\vpsi_k\equiv 0 .$$  
  \end{itemize}
   
  \bre
 {\rm 
 Suppose $\lambda_k^\prime\in \mathbb{Z}$. In particular, if $\lp_k\leq -2$, suppose we are in the case when $\vec{\Psi}_k^{(sing)}$ is not identically zero. Then
$$ 
\vec{\Psi}_k(\lambda) =\frac{1}{2\pi i} \left( \vec{\Psi}_k^{(sing)}(l_k(\lambda))-\vec{\Psi}_k^{(sing)}(\lambda)\right),\quad \lambda\in\mathcal{P}_{\eta^{(0)}},
$$
is the difference of two singular solutions defined on $\mathcal{P}_{\eta^{(0)}}$.  
Here, in the notation of Remark \ref{3aprile2021-1}, the   function $ \vec{\Psi}_k^{(sing)}(l_k(\lambda))$     is the value at $ \lambda\in\mathcal{P}_{\eta^{(0)}}$  of the analytic continuation of $\vec{\Psi}_k^{(sing)}(\lambda)$ when passing from a prefixed component of $\mathcal{R}^\prime$, in this case   $\mathcal{P}_{\eta^{(0)}}$, to the component  associated with the  sequence $\{l_k\}$ of only one element. Namely,   the analytic continuation for  a small loop $(\lambda-u_k^{(0)})\longmapsto (\lambda-u_k^{(0)})e^{2\pi i}$.
}
  \ere

 \subsubsection*{Connection Coefficients}
  Above, the behaviour of $\vec{\Psi}_k(\lambda)$ has been described at $\lambda=u_k^0$. The behaviour at any point $\lambda=u_j^0$, for $j=1,...,n$, will be expressed by linear relations
\be
\label{30agosto2020-12}
 \vPsi_k(\lambda)=\vPsi_j^{(sing)}(\lambda) c_{jk}+\hbox{\rm reg}(\lambda-u_j^0).
\ee
$$ c_{jk}:=0,\quad \forall k=1,...,n,\quad \hbox{ when }  \vPsi_j^{(sing)}(\lambda)\equiv 0 \hbox{ (possible only if $\lp_j\in -\mathbb{ N}-2$)}.
$$
The above relations  {\it define} the {\bf connection coefficients} $c_{jk}$. 
 From the definition, we see that $
c_{kk}=1 $ for $ \lp_k\not\in \mathbb{ Z}$,  while $c_{kk}=0$ for $\lp_k\in \mathbb{Z}$. 
  In case $\lp_k\in \mathbb{ N}$, if it happens that $\vPsi_k\equiv 0$, then  $c_{jk}=0$ for any $j=1,..,n$.

\begin{figure}
\centerline{\includegraphics[width=0.6\textwidth]{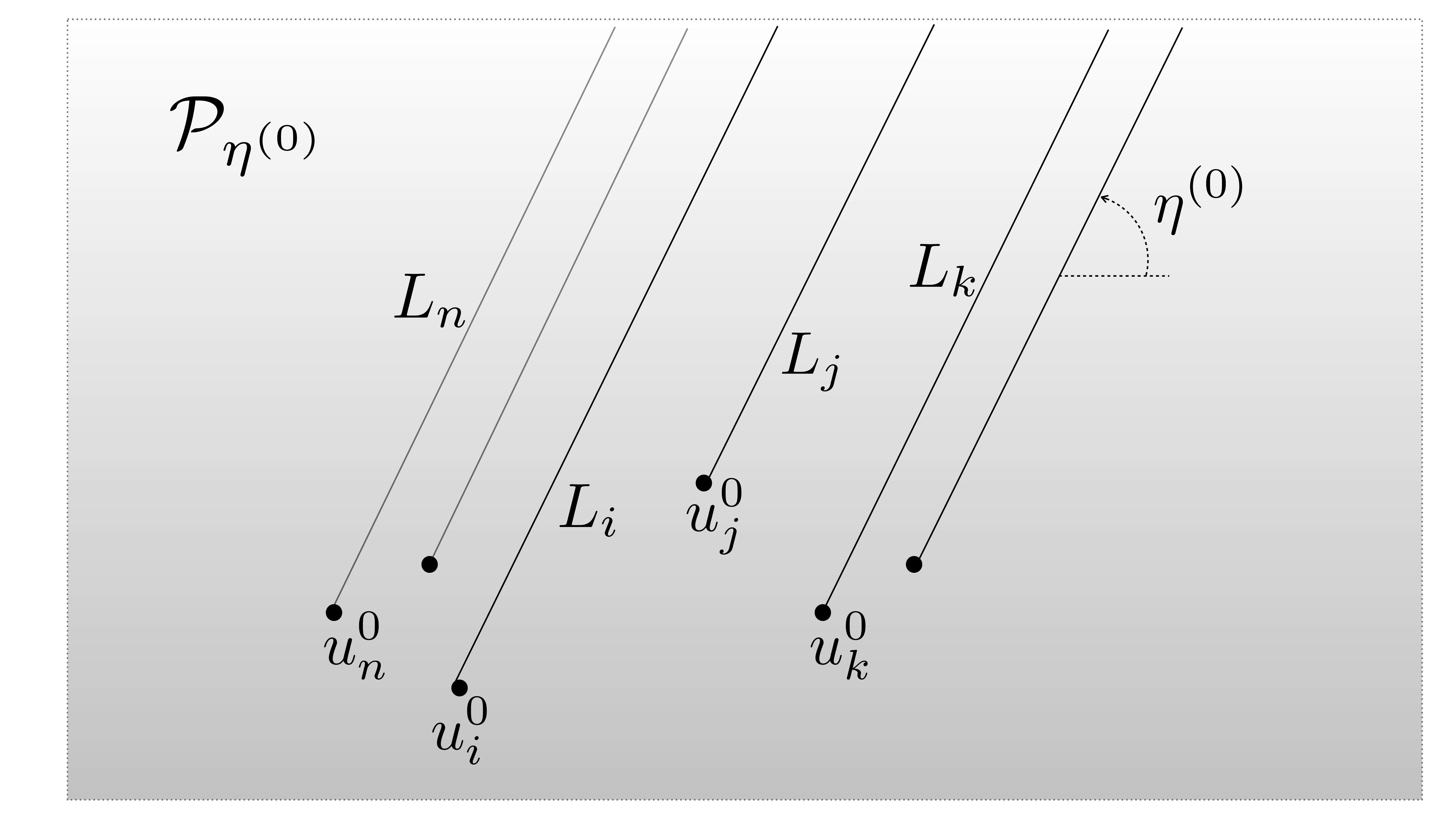}}
\caption{The poles $u_j^0$, $1\leq j \leq n$, of system (\ref{29agosto2020-1}) and plane $\mathcal{P}_{\eta^{(0)}}$ with branch cuts $L_j$.}
\label{cutss}
\end{figure}

\bpr[see \cite{BJL4} and propositions 3, 4 of \cite{guz2016} ]
\label{26novembre2020-2}
 If $A$ has no integer eigenvalues, then   \be
\label{29agosto2020-3}
\Psi(\lambda)=\Bigr[\vPsi_1(\lambda)~|~\cdots~|~\vPsi_n(\lambda)\Bigr],~~~~~\lambda\in{\cal P}_{\eta^{(0)}}
\ee 
(each $\vec{\Psi_k}$ occupies a column) is a fundamental matrix solution of \eqref{29agosto2020-1}.  Moreover,  the matrix $C:=(c_{jk})$ is invertible if and only if $A$ has no integer eigenvalues. 
 If  $A$ has integer eigenvalues and $\Psi$ is fundamental, then  some $\lp_k\in \mathbb{Z}$. 
\epr

\subsubsection*{Laplace transform and Stokes Matrices in terms of Connection Coefficients}

If $\eta^{(0)}$ is admissible in the $\lambda$-plane, with respect to the {\it fixed} and  pairwise distinct  $u_1^0,...,u_n^0$, then $$
\arg z=\tau^{(0)} := 3\pi/2-\eta^{(0)}$$
is an admissible direction \eqref{29agosto2020-2} in the $z$-plane for system \eqref{24nov2018-1} {\it at the fixed $u=u^0$}. We consider the Stokes rays of $\Lambda(u^0))$ as before. For some $\nu\in\mathbb{Z}$, a labelling  \eqref{22novembre2020-1} holds, so that 
\be
\label{8novembre2020-1}
\tau_\nu<\tau^{(0)}<\tau_{\nu+1} \quad \Longleftrightarrow \quad  \quad \eta_{\nu+1}<\eta^{(0)}<\eta_\nu ,\quad \quad \eta_\nu:=\frac{3\pi}{2}-\tau_\nu.
\ee 
In order to keep track of \eqref{8novembre2020-1}, we  label \eqref{29agosto2020-3} with $\nu$,
\be
\label{29agosto2020-7}
\Psi_\nu(\lambda)=\Bigr[\vPsi_1(\lambda~|\nu)~|~\cdots~|~\vPsi_n(\lambda~|\nu)\Bigr],\quad \lambda\in\mathcal{P}_{\eta^{(0)}}. \ee
The connection coefficients will be labelled accordingly as  $c^{(\nu)}_{jk}$. 
Also the singular vector solutions will be labelled   $\vec{\Psi}_k^{(sing)}(\lambda~|\nu)$,  $\lambda\in \mathcal{P}_{\eta^{(0)}}$ as above.

\vskip 0.2 cm 
The  relation between   solutions  $\vec{\Psi}_k(\lambda~|\nu)$ or $\vec{\Psi}_k^{(sing)}(\lambda~|\nu)$ and the columns of $
Y_\nu(z)$ is established in \cite{guz2016} for all values of $\lp_1,...,\lp_n$, and in  \cite{BJL4} for  non integer values only.  It  is given by  Laplace-type integrals (Proposition 8 of \cite{guz2016})
$$
\vec{Y}_k(z~|\nu)= \dfrac{1}{2\pi i } \int_{\gamma_k(\eta^{(0)})} e^{z\lambda} \vec{\Psi}_k^{(sing)}(\lambda~|\nu) d\lambda,  \hbox{ if } \lp_k\not\in\mathbb{Z}_{-};
\quad
\vec{Y}_k(z~|\nu)= \int_{L_k(\eta^{(0)})} e^{z\lambda} \vec{\Psi}_k(\lambda~|\nu) d\lambda,   \hbox{ if }  \lp_k\in\mathbb{Z}_{-}.
$$
Here,      $\gamma_k(\eta^{(0)})$ is the path coming from $\infty$ along the left side of the oriented $L_k(\eta^{(0)})$, encircling $u_k^0$ with a small loop excluding all the other poles, and going back to $\infty$ along the right side of $L_k(\eta^{(0)})$. 
    
    \vskip 0.2 cm 
    The same as \eqref{29agosto2020-7}  can be defined for the cut-plane  $\mathcal{P}_{\eta^\prime}$, with an admissible   direction $\eta^\prime$  satisfying   $$\eta_{\nu+k\mu^{(0)}+1}<\eta^\prime<\eta_{\mu^{(0)}+k\mu^{(0)}},\quad k\in\mathbb{Z},
    $$  and will be denoted by $\Psi_{\nu+k\mu^{(0)}}(\lambda)$, and analogously for the vectors  $\vec{\Psi}_k(\lambda~|\nu+k\mu^{(0)})$ and $ \vec{\Psi}_k^{(sing)}(\lambda~|\nu+k\mu^{(0)})$.  From the Laplace transforms of  $\vec{\Psi}_k(\lambda~|\nu+k\mu^{(0)})$ or $ \vec{\Psi}_k^{(sing)}(\lambda~|\nu+k\mu^{(0)})$, with the paths of integration  $\gamma_k(\eta^\prime)$ or $L_k(\eta^\prime)$, we receive $Y_{\nu+k\mu^{(0)}}(z)$.

\vskip 0.2 cm 
 Introduce in $\{1,2,...,n\}$ the ordering $\prec$ given by 
 $$
 j\prec k ~~\Longleftrightarrow~~ \Re(z(u_j^0-u_k^0))<0 ~\hbox{ for } \arg z=\tau^{(0)},~~~~~i\neq j,~~i,j\in\{1,...,n\}.
 $$
The following important results, proved in  theorem 1 of  \cite{guz2016} for all values of $\lp_1,...,\lp_n$, and  in the seminal paper  \cite{BJL4} in the generic  case  $\lp_1,...,\lp_n\not \in\mathbb{Z}$, establishes the relation between Stokes matrices and connection coefficients.\footnote{The key point is the fact  that $\vec{\Psi}_k^{(sing)}$ in \eqref{15agosto2020-1},
   or equivalently $\vec{\Psi}_k$ for $\lp_1,...,\lp_n\not\in\mathbb{Z}$,  can be substituted by another 
   set of vector solutions, denoted in \cite{guz2016}  by $\vec{\Psi}_k^*(\lambda, u~|\nu)$ and in \cite{BJL4} by $Y^*_k$. The effect of the change of the branch cut from $\eta_{\nu+1}<\eta<\eta_{\nu}$ to  $\eta_{\nu+\mu+1}<\eta^\prime<\eta_{\nu+\mu}$, namely from $\eta$ to $\eta^\prime=\eta-\pi$,  yields a 
     linear relation  
     $$\vec{\Psi}_k^*(\lambda, u~|\nu+\mu)=\vec{\Psi}_k^*(\lambda, u~|\nu)C^+_\nu
     ,\quad \lambda\in \mathcal{P}_\eta\cap\mathcal{P}_{\eta -\pi},
     $$  where the connection matrix $C_\nu^+$ is  expressed in terms of the connection coefficients $c_{jk}^{(\nu)}=c_{jk}^{(\nu)}(\eta)$ 
    associated with  $\vec{\Psi}_k^{(sing)}(\lambda, u~|\nu)$.  
    The same can be done for the change of branch cut from $\eta_{\nu+\mu+1}<\eta^\prime<\eta_{\nu+\mu}$ to $\eta_{\nu+2\mu+1}<\eta^{\prime\prime}<\eta_{\nu+2\mu}$ (namely,  $\eta^\prime=\eta-\pi$ and $\eta^{\prime\prime}=\eta-2\pi$) yielding a relation 
    $$\vec{\Psi}_k^*(\lambda, u~|\nu+2\mu)=\vec{\Psi}_k^*(\lambda, u~|\nu+\mu)C^-_\nu,\quad \lambda\in\mathcal{P}_{\eta-\pi}\cap \mathcal{P}_{\eta-2\pi}.
    $$
    Substituting these relations in the Laplace integrals, one proves Theorem \ref{20agosto2020-1}, being  $\mathbb{S}_\nu=C_\nu^+$ and $\mathbb{S}_{\nu+\mu}^{-1}=C_\nu^{-}$. See \cite{BJL4} and \cite{guz2016}}

\bth  
\label{20agosto2020-1}
Let $u=u^0$ be fixed so that $\Lambda(u^0)$ has pairwise distinct eigenvalues.  Let $\eta^{(0)}$  and  $\tau^{(0)} = 3\pi/2-\eta^{(0)}$ be admissible  for $u^0$ in the $\lambda$-plane and  $z$-plane respectively. Suppose that the labelling of Stokes rays is \eqref{22novembre2020-1} and  \eqref{8novembre2020-1}. 
Then, the Stokes  matrices of system \eqref{24nov2018-1} at $u=u^0$ are given in terms of the connection coefficients $c_{jk}^{(\nu)}$ of system  \eqref{29agosto2020-1}, according to the following  formulae
$$
\bigl(\mathbb{S}_\nu\bigr)_{jk}
=
\left\{
\begin{array}{cc}
e^{2\pi i \lp_k}\alpha_k ~c^{(\nu)}_{jk}& ~~~\hbox{ for } j\prec k,
\\
\\
             1   & ~~~\hbox{ for } j =k,
\\
\\
             0   &  ~~~\hbox{ for } j\succ k,
\end{array}
\right.
~~~~~~~~
\bigl( \mathbb{S}_{\nu+\mu^{(0)}}^{-1}\bigr)_{jk}
=
\left\{
\begin{array}{cc}
             0   &  ~~~\hbox{ for } j\prec k,
\\
\\
             1   &~~~ \hbox{ for } j =k,
\\
\\
        -e^{2\pi i (\lp_k-\lp_j)}\alpha_k~c^{(\nu)}_{jk} & 
                                         ~~~\hbox{ for } j\succ k.
\end{array}
\right.
$$
where, 
$$
\alpha_k:=(e^{-2\pi i \lp_k}-1)\quad \hbox{ if }\lp_k\not\in \mathbb{Z};\quad\quad \alpha_k:=2\pi i  \quad\hbox{ if }\lp_k\in \mathbb{Z}.
$$
\eth
\noindent
$\Box$

In the above discussion, the differential systems do not depend on parameters ($u$ is fixed). 
The purpose of the present paper is to extend the description of  Background 2 to the case depending on deformation parameters and include coalescences  in $\mathbb{D}(u^c)$,  and then  to obtain Theorem \ref{18agosto2020-8}  of Background 1 in terms of an isomonodromic Laplace transform.

\section{Equivalence of the Isomonodromy Deformation Equations for  \eqref{24nov2018-1}  and \eqref{03}}
\label{19agosto2020-1}

The first step in our construction is Proposition \ref{11agosto2020-5} below, establishing  the equivalence between strong isomonodromy deformations of systems \eqref{24nov2018-1} and \eqref{03}, for $u$ varying in a $\tau$-cell of $\mathbb{D}(u^c)$. In the specific case of Frobenius manifolds, this fact can be deduced from   Chapter 5 of \cite{Dub2}. Here we establish the equivalence in general terms. 

  According to Theorem \ref{9agosto2020-1}, system \eqref{24nov2018-1} is strongly isomonodromic in a polydisc   $\mathbb{D}(u^0)$ contained in a $\tau$-cell of $\mathbb{D}(u^c)$ if and only if \footnote{As already mentioned when stating Theorem \ref{9agosto2020-1}, equations $dA=[\omega_i(u),A]$ and $\omega_i(u)=[F_1,E_i]$ for $i=1,...,n$ are exactly the the Frobenius integrability conditions of \eqref{28nov2018-1}  when   \eqref{24nov2018-1}  is strongly isomomodromic \cite{CDG}.} 
\be
\label{10agosto2020-1}
dA=\sum_{j=1}^n [\omega_j(u),A]~du_j,\quad \quad \omega_j(u)=[F_1(u),E_j], \hbox{ given in \eqref{22agosto2020-6}}.
\ee

 On the other hand, system \eqref{03} is strongly isomonodromic in  $\mathbb{D}(u^0)$ by definition  (\cite{Guzz-SIGMA}, Appendix A), when fundamental matrix solutions in Levelt form at each pole $\lambda=u_j$, $j=1,...,n$, have {\it constant monodromy exponents} and are related to each other  by {\it constant connection matrices} (not to be confused with the connection coefficients). From  \cite{Bolibruch,Bolibruch1,Guzz-SIGMA}, the necessary and sufficient condition for  the deformation to be strongly isomonodromic (this can also be taken as the definition)  is    that \eqref{03} is the $\lambda$-component of a Frobenius integrable Pfaffian system with the following structure
\be
\label{11agosto2020-1}
d\Psi=P(\lambda,u) \Psi,\quad\quad P(\lambda,u)=\sum_{k=1}^n\frac{B_k(u)}{\lambda-u_k}d(\lambda-u_k) +\sum_{k=1}^n \gamma_k(u) du_k.
\ee
The integrability condition $dP=P \wedge P$ is  the non-normalized Schlesinger system (see Appendix A and   \cite{Bolibruch0,Bolibruch,Bolibruch1,Guzz-SIGMA,Guzz-notes,YT}) 
\begin{align} 
\label{10agosto2020-3}
&\partial_i\gamma_k-\partial_k\gamma_i=\gamma_i\gamma_k-\gamma_k\gamma_i,
\\
\label{10agosto2020-4}
&\partial_iB_k=\frac{[B_i,B_k]}{u_i-u_k}+[\gamma_i,B_k], \quad i\neq k
\\
\label{10agosto2020-5}
  & \partial_iB_i= -\sum_{k\neq i} \frac{[B_i,B_k]}{u_i-u_k}+[\gamma_i,B_i]
  \end{align}

  \bpr
  \label{11agosto2020-5}
The system  \eqref{10agosto2020-1} is equivalent to  \eqref{10agosto2020-3}-\eqref{10agosto2020-5} if and only if 
   $$\gamma_j(u)\equiv \omega_j(u)  \hbox{ as in \eqref{9agosto2020-2} and \eqref{22agosto2020-6}},\quad\quad j=1,...,n.
$$
 Namely,  \eqref{24nov2018-1} is strongly isomonodromic in a polydisc on $\mathbb{D}(u^0)$ contained in a $\tau$-cell if and only if \eqref{03} is strongly isomonodromic.
  \epr

  \begin{proof} See Appendix B. \end{proof}
  
  \section{Schlesinger System on $\mathbb{D}(u^c)$ and Vanishing Conditions}
  \label{24agosto2020-3}
In this section, Proposition \ref{8novembre2020-2}, we  holomorphically  extend to $\mathbb{D}(u^c)$ the non-normalized Schlesinger system associated with \eqref{03}, when certain  vanishing conditions \eqref{30agosto2020-3} are satisfied. This is the second step to obtain the results of \cite{CDG} by Laplace transform. 
  
  \ble
  \label{12agosto2020-2-2}
  Let  $A(u)$ be holomorphic on $\mathbb{D}(u^c)$ and $B_j(u):=-E_j(A(u)+I)$, $j=1,...,n$. 
  
  i) The vanishing  relations 
  \be
   \label{12agosto2020-1}
   [B_i(u),B_j(u)]\longrightarrow 0, \quad \hbox{ for  $u_i-u_j\to 0$ in $\mathbb{D}(u^c)$}.
   \ee
    hold 
   if and only if
   \be
   \label{30agosto2020-3}
    \bigl(A(u)\bigr)_{ij}\longrightarrow 0, \quad \hbox{ for  $u_i-u_j\to 0$ in $\mathbb{D}(u^c)$}. 
   \ee
   
 ii) The matrices $\omega_k(u)=[F_1(u),E_k]$ are holomorphic on $\mathbb{D}(u^c)$ if and only if \eqref{30agosto2020-3} holds.
   \ele

  \begin{proof} Let $u^*\in \Delta$, so that for some $i\neq j$ it occurs that $u_i^*=u_j^*$. Since 
  \be
  \label{22agosto2020-5}
  B_j=-E_j(A+I)=\begin{pmatrix} 0   &&&0&&&0
\\
                        \vdots &&&\vdots&&&\vdots
                        \\
 -A_{j1}&\cdots& -A_{j,j-1}&-\lp_j-1&-A_{j,j+1}& \cdots&-A_{jn}
 \\
\vdots&&&\vdots&&&\vdots
\\
0&&&0&&&0
\end{pmatrix}.
\ee
it is an elementary computation to check the equivalence between the relation  $[B_i(u^*),B_j(u^*)]=0$ and the relation $(A(u^*))_{ij}=0$. Since  $[F_1(u),E_k]$ is \eqref{22agosto2020-6}, the statement on its analyticity  is straightforward. 
\end{proof}

\bpr
  \label{8novembre2020-2}
  Consider a Frobenius integrable Pfaffian system \eqref{11agosto2020-1} on $\mathbb{D}(u^0)$ 
   with 
   \be
   \label{8novembre2020-6}
   B_j(u)=-E_j(A(u)+I)\quad \hbox{ and }\quad \gamma_j(u)\equiv\omega_j(u)=[F_1(u),E_j] \hbox{ in \eqref{22agosto2020-6}}.
   \ee
    Assume that $A(u)$ is holomorphic on the whole $\mathbb{D}(u^c)$. Then,  system  \eqref{11agosto2020-1}  is Frobenius integrable on the whole $\mathbb{D}(u^c)$ with holomorphic matrix coefficients if and only if the vanishing conditions \eqref{30agosto2020-3} hold.  
  \epr
  
  \begin{proof}
  If system \eqref{11agosto2020-1} is integrable on $\mathbb{D}(u^c)$ with holomorphic coefficients $B_k$ and $\gamma_k=\omega_k$, then analyticity of $\omega_k$ with structure  \eqref{22agosto2020-6}  implies  that  \eqref {30agosto2020-3} must hold, so that \eqref{12agosto2020-1} holds. Notice that by \eqref{12agosto2020-1}, the r.h.sides  of \eqref{10agosto2020-4}-\eqref{10agosto2020-5} are  holomorphic on $\mathbb{D}(u^c)$. 
    Conversely, suppose that \eqref{12agosto2020-1}-\eqref{30agosto2020-3} hold. By Proposition \ref{11agosto2020-5}, \eqref{10agosto2020-3}-\eqref{10agosto2020-4}-\eqref{10agosto2020-5} 
  are equivalent in $ \mathbb{D}(u^0)$ to 
\be
\label{5aprile2021-2}
  dA=\sum_{j=1}^n [\omega_j(u),A]~du_j,
  \quad u\in\mathbb{D}(u^0).
  \ee
  Now, the l.h.s  is well defined and holomorphic   on $\mathbb{D}(u^c)$, because   so is $A(u)$. The r.h.s. is also analytic on $\mathbb{D}(u^c)$, because  of  \eqref{30agosto2020-3}. Hence, the first part of the proof of Proposition \ref{11agosto2020-5} in Appendix B works in the whole $\mathbb{D}(u^c)$, and so the Pfaffian system \eqref{11agosto2020-1} is integrable there. 
  \end{proof}
  
  For completeness, we also state the
   following 
   \bpr
    Let system \eqref{11agosto2020-1}  with coefficients  \eqref{8novembre2020-6} be integrable on 
    $\mathbb{D}(u^c)\backslash\Delta$ and let the $B_k(u)$  be holomorphic on $\mathbb{D}(u^c)$. Then the vanishing conditions \eqref{12agosto2020-1}-\eqref{30agosto2020-3},  hold and the $\omega_k$ are holomorphic on $\mathbb{D}(u^c)$. 
   \epr
   
   \begin{proof}
   Analogously to the proof of Proposition  \ref{11agosto2020-5},  we see that \eqref{10agosto2020-3}-\eqref{10agosto2020-4}-\eqref{10agosto2020-5} on $\mathbb{D}(u^c)\backslash\Delta$ 
  are equivalent to 
\be
\label{5aprile2021-2-bis}
  dA=\sum_{j=1}^n [\omega_j(u),A]~du_j,
  \quad u\in\mathbb{D}(u^c)\backslash\Delta.
\ee
By holomorphy of $A(u)$ on $\mathbb{D}(u^c)$, the r.h.s is well defined, so that also the l.h.s. must be holomorphic on $\mathbb{D}(u^c)$.   From \eqref{5aprile2021-2-bis} we proceed as in Remark \ref{6aprile2021-2}, concluding that   $ A_{ij} =O(u_i-u_j)\to 0$ holomorphically for $u_i-u_j\to 0$.     The  proof can be done also with an argument similar to Remark \ref{7aprile2021-2}.  \end{proof}

  \section{Selected Vector solutions   depending on parameters $u\in \mathbb{D}(u^c)$, Theorem \ref{30agosto2020-5}}
  \label{30agosto2020-1}
  
   In this section we state one main result of the paper,  Theorem \ref{30agosto2020-5}  below,   introducing the  isomonodromic analogue of the selected and singular  vector  solutions \eqref{20agosto2020-19} and \eqref{20agosto2020-16}.  This is the  third  step required to obtain the results of \cite{CDG} by Laplace transform.

  Preliminarily, we   characterise the radius $\epsilon_0>0$ of  $\mathbb{D}(u^c)$ in \eqref{22novembre2020-5}.    The coalescence point $ 
  u^c=(u_1^c,...,u_n^c)$  contains $s<n$ distinct values, say  $\lambda_1,...,\lambda_s$, with algebraic multiplicities $p_1$, ..., $p_s$ respectively ($p_1+\dots+p_s=n$). Suppose that  $\arg z =\tau$ is a direction   {\it admissible  at $u^c$}, as defined in \eqref{20agosto2020-4}, and  let 
 $$\eta={3\pi/2}-\tau$$ 
 be the corresponding admissible direction in the $\lambda$-plane, where we draw parallel half lines $\mathcal{L}_1=\mathcal{L}_1(\eta)$, ..., $\mathcal{L}_s=\mathcal{L}_s(\eta)$  issuing from $\lambda_1$, ..., $\lambda_s$ respectively, with direction $\eta$,  as in  figure \ref{2settembre2020-10}.  Let 
 $$ 
 2\delta_{\alpha \beta }:=\hbox{ distance between $\mathcal{L}_\alpha$ and $\mathcal{L}_\beta$, for $1\leq \alpha \neq \beta\leq s$} 
 $$
 In formulae, $2\delta_{\alpha \beta}=\min_{\rho>0}|\lambda_\alpha-\lambda_\beta +\rho e^{\sqrt{-1}( 3\pi/2-\tau)}|$. 
Then, we require that 
\be
\label{30agosto2020-6}
\epsilon_0<\min_{1\leq \alpha \neq \beta \leq n} \delta_{\alpha \beta}.
\ee
The bound \eqref{30agosto2020-6} was introduced in \cite{CDG} in order to prove Theorem \ref{18agosto2020-8} in Background 1. It implies properties of the Stokes rays as $u$ varies in $\mathbb{D}(u^c)$,  described later in Section \ref{20agosto2020-10}. Let 
$$ 
\mathbb{D}_\alpha:= \{\lambda\in\mathbb{C}~|~|\lambda-\lambda_\alpha|\leq \epsilon_0\},\quad \alpha=1,...,s,
$$
be the disc centered a $\lambda_\alpha$ and radius $\epsilon_0$. 
  If $u_j$ is such that $u_j^c=\lambda_\alpha$,  the  bound  \eqref{30agosto2020-6} implies that $u_j$ remains in $\mathbb{D}_\alpha$ as 
   $u$ varies in $\mathbb{D}(u^c)$. Clearly, $\mathbb{D}_\alpha\cap \mathbb{D}_\beta=\emptyset$. 

  \begin{figure}
\centerline{\includegraphics[width=0.8\textwidth]{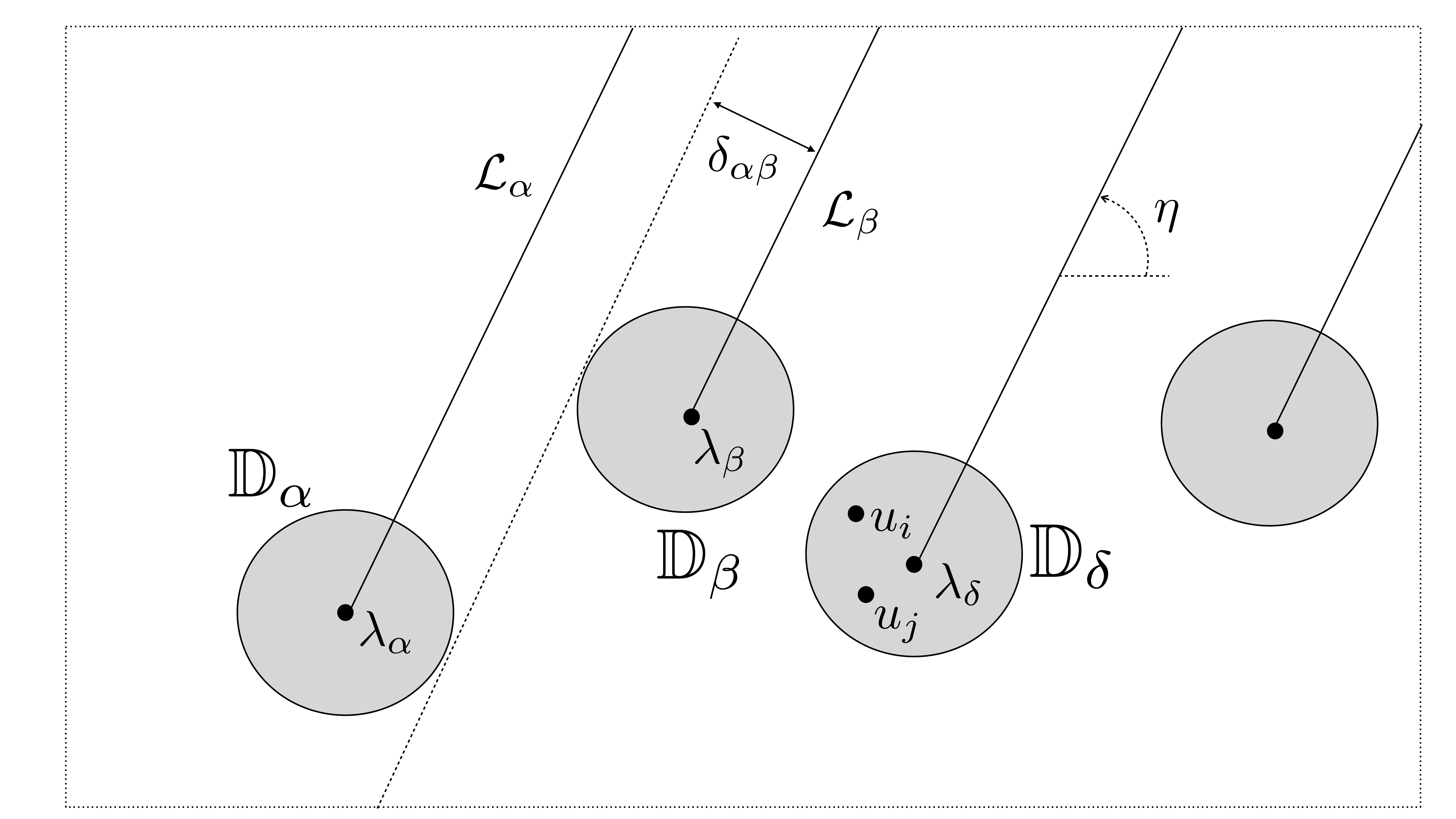}}
\caption{The figure represents the half lines $\mathcal{L}_\alpha$, $\mathcal{L}_\beta$, etc,  for $\alpha,\beta,...\in\{1,...,s\}$, in direction $\eta=3\pi/2-\tau$, the discs centred at the coordinates $\lambda_1,\dots,\lambda_s$ of the coalescence point $u^c$, and  the distances $\delta_{\alpha\beta}$. Also  two points  $u_i,u_j$ are represented, such that $u_i^c=u_j^c=\lambda_\delta$ for some $\delta\in \{1,...,s\}$. {\bf Important:} now $\eta$ refers to $u^c$, differently from Section \ref{19agosto2020-3} and figure \ref{cutss}. }
\label{2settembre2020-10}
\end{figure}

  \vskip 0.2 cm 
  
 The Stokes rays of $\Lambda(u^c)$ can be labeled as in \eqref{15novembre2020-1}. For a certain  $\nu\in\mathbb{Z}$ we have 
  \be
\label{23agosto2020-16}
\eta_{\nu+1}<\eta<\eta_\nu \quad \Longleftrightarrow \quad \tau_\nu <\tau <\tau_{\nu+1},\quad \quad \eta_\nu=\frac{3\pi}{2}-\tau_\nu.
\ee  
  For each $u\in\mathbb{D}(u^c)$, let  $\mathcal{P}_{\eta}=\mathcal{P}_{\eta}(u)$   be the $\lambda$-plane with  branch cuts $L_1=L_1(\eta)$, ..., $L_n=L_n(\eta)$ issuing from $u_1,...,u_n$  and the choice of the logarithms 
 $ \ln(\lambda- u_k)=\ln|\lambda- u_k|+i\arg(\lambda-u_k)$, given by 
 $$ \eta-2\pi <\arg(\lambda-u_k)<\eta,\quad \quad k=1,...,n.$$ 
   We define the domain (notation $\hat{\times}$  inspired by \cite{JMU})
  $$ 
  \mathcal{P}_\eta(u)\hat{\times}\mathbb{D}(u^c):=\{(\lambda,u)~|~ u\in \mathbb{D}(u^c),~\lambda\in \mathcal{P}_\eta(u)\}.
  $$

     According to Proposition \ref{8novembre2020-2}, for a Pfaffian system \eqref{11agosto2020-1} with coefficients \eqref{8novembre2020-6}, defined on a polydisc $\mathbb{D}(u^0)$  contained in a $\tau$-cell of $\mathbb{D}(u^c)$, if   $A(u)$ is holomorphic on  $\mathbb{D}(u^c)$, then the vanishing conditions  \eqref{30agosto2020-3}
     $$
      \bigl(A(u)\bigr)_{ij}\longrightarrow 0, \quad \hbox{ for  $u_i-u_j\to 0$ in $\mathbb{D}(u^c)$}.
      $$ are equivalent to Frobenius integrability on the whole $\mathbb{D}(u^c)$. With this in mind, we state the following 
  
  \bth
\label{30agosto2020-5}  Consider a  Pfaffian system,   
 \be
  \label{19agosto2020-5}
  d\Psi=P(\lambda,u) \Psi,\quad\quad P(z,u)=\sum_{k=1}^n\frac{B_k(u)}{\lambda-u_k}d(\lambda-u_k) +\sum_{k=1}^n \omega_k(u) du_k.
\ee
Frobenius integrable  on  $\mathbb{D}(u^c)$, with matrix coefficients \eqref{22agosto2020-6}, and   $A(u)$  holomorphic on  $\mathbb{D}(u^c)$.  Let the radius $\epsilon_0$ be as in \eqref{30agosto2020-6}.  Then,  two classes of vector solutions,  holomorphic on  $\mathcal{P}_\eta(u)\hat{\times}\mathbb{D}(u^c)$, exist as follows.

  \vskip 0.2 cm
 The {\bf  selected solution:} $\vPsi_1(\lambda,u~|\nu),~...~,~\vPsi_n(\lambda,u~|\nu)$. Each $\vPsi_k(\lambda,u~|\nu)$ is uniquely identified   by the  local behaviour below for $\lambda\in \mathbb{D}_\alpha$,  where $\alpha$ is  such that  $u_k^c=\lambda_\alpha$. The label $\nu$ keeps track of \eqref{23agosto2020-16}.
  
  \begin{itemize}

 \item  \underline{For $\lp_k\in\mathbb{C}\backslash \mathbb{Z}$ or  $\lp_k\in\mathbb{Z}_{-}=\{-1,-2,...\}$}, 
 \be
 \label{18ottobre2020-1} 
   \vPsi_k(\lambda,u~|\nu)=\vec{\psi}_k(\lambda,u~|\nu) (\lambda-u_k)^{-\lp_k-1},\quad k=1,...,n,
   \ee
   where $\vec{\psi}_k(\lambda,u~|\nu)$ is  holomorphic on $ \mathbb{D}_\alpha\times \mathbb{D}(u^c)$ and is represented by  a   uniformly convergent Taylor expansion  with  holomorphic on  $ \mathbb{D}(u^c)$ coefficients:
   \be
   \label{31agosto2020-1}
   \vec{\psi}_k(\lambda,u~|\nu)=f_k\vec{e}_k+
   \sum_{l=1}^\infty \vec{b}_l^{~(k)}(u) (\lambda-u_k)^l ,\quad \hbox{ for } \lambda\to u_k,
   \ee
 The following  normalization uniquely identifies $\vec{\Psi}_k$. 
  \be
  \label{19novembre2020-1}
f_k=
\left\{
\begin{array}{ccc} \Gamma(\lp_k+1),&  \lp_k\in\mathbb{C}\backslash \mathbb{Z},
\\
\noalign{\medskip}
\dfrac{(-1)^{\lp_k}}{(-\lp_k-1)!}, & \lp_k \in \mathbb{Z}_{-},
\end{array}
\right.
\ee

   \item  \underline{For  $\lp_k\in\mathbb{N}=\{0,1,2,...\}$},  
     \be
   \label{31agosto2020-1-bis}
   \vPsi_k(\lambda,u~|\nu)=\sum_{l=0}^\infty \vec{d}_l^{~(k)}(u)(\lambda-u_k)^l  ,\quad \hbox{ for } \lambda\to u_k,
   \ee
   is holomorphic on $ \mathbb{D}_\alpha\times \mathbb{D}(u^c)$, 
  the Taylor expansion being  uniformly convergent  with holomorphic   coefficients $
   \vec{d}_l^{~(k)}(u)$.  It is uniquely identified by the normalization  \eqref{17ottobre2020-2} of the singular solution \eqref{17ottobre2020-1} below. 
   Depending on the specific Pfaffian system \footnote{See the comment to  \eqref{19novembre2020-3} below.}, it may happen that identically $$\vPsi_k(\lambda,u~|\nu)\equiv 0.
   $$ 

      \end{itemize}
      
 The isolated  singularities of  $\vPsi_k(\lambda,u~|\nu)$, if any, are located at  $\lambda=u_j$ with  $u_j^c=\lambda_\beta$, $\beta\neq \alpha$, and   at $\lambda=u_k$  only in  case  $\lp_k\in\mathbb{C}\backslash \mathbb{Z}$.  For  $i\neq j$ such that $u_i^c=u_j^c$,  $\vec{\Psi}_i(\lambda,u~|\nu)$ and $\vec{\Psi}_j(\lambda,u~|\nu)$ are either linearly independent, or at least one of them is identically zero (identity to zero  may occur  only for    $\lp_i$ or   $\lp_j$ belonging to $\mathbb{N}$)

  \vskip 0.3 cm 

The   {\bf singular solutions:} $ 
  \vPsi_1^{(sing)}(\lambda,u~|\nu),~...~,~\vPsi_n^{(sing)}(\lambda,u~|\nu)
  $. Each    $\vPsi_k^{(sing)}(\lambda,u~|\nu)$ is a solution with an isolated singularity  at $\lambda=u_k$, whose singular behaviour is  uniquely characterized as follows.\footnote{The solution here defined   is not  uniquely identified by the singular behaviour if $\lp_k\in \mathbb{Z}_{-}$, see Remark \ref{7aprile2021-1}.}  Let $\mathbb{D}_\alpha$ be identified by  $\lambda_\alpha=u_k^c$.  \begin{itemize}

   \item  \underline{For $\lp_k\in\mathbb{C}\backslash\mathbb{Z}$} {\rm [algebraic or logarithmic branch-point]},
$$ 
\vec{\Psi}_k^{(sing)}(\lambda,u~|\nu):= \vec{\Psi}_k(\lambda,u~|\nu)=\vec{\psi}_k(\lambda,u~|\nu) (\lambda-u_k)^{-\lp_k-1}.
$$
\item \underline{For $\lp_k\in\mathbb{Z}_{-}$} {\rm [logarithmic branch-point]},
\begin{align}
\label{18ottobre2020-2}
\vec{\Psi}_k^{(sing)}(\lambda,u~|\nu)&=\vec{\Psi}_k(\lambda,u~|\nu) \ln (\lambda-u_k)+\sum_{m\neq k}^* r_m\vec{\Psi}_m(\lambda,u~|\nu) \ln (\lambda-u_m)+ \vec{\phi}_k(\lambda,u~|\nu) ,
\\
\label{18ottobre2020-3}
&\underset{\lambda\to u_k} =\vec{\Psi}_k(\lambda,u~|\nu) \ln (\lambda-u_k)+\hbox{\rm reg}(\lambda-u_k),\quad \quad r_m\in\mathbb{C},
\end{align}
where $\sum_{m\neq k}^*$ is  over all $m$ such that $u_m\in \mathbb{D}_\alpha$ and $\lp_m\in\mathbb{Z}_{-}$.  The vector  function $ \vec{\phi}_k(\lambda,u~|\nu)$ is holomorphic in $\mathbb{D}_\alpha\times \mathbb{D}(u^c)$. 

In particular, \underline{for $\lp_k\leq -2$}, depending on the  system,  it may happen that     there is no   solution with singularity  in $\mathbb{D}_\alpha$,  so that $$\vec{\Psi}_k^{(sing)}(\lambda,u~|\nu):=0.$$

\item \underline{For $\lp_k\in\mathbb{N}$} {\rm [logarithmic branch-point and pole]}, 
\be
\label{17ottobre2020-1}
\vec{\Psi}_k^{(sing)}(\lambda,u~|\nu)= \vec{\Psi}_k(\lambda,u~|\nu)\ln(\lambda-u_k)+\frac{\vec{\psi}_k(\lambda,u~|\nu)}{(\lambda-u_k)^{\lp_k+1}},
\ee
where $\vec{\psi}_k(\lambda,u~|\nu)$ is holomorphic in $\mathbb{D}_\alpha\times \mathbb{D}(u^c)$,  
\be
\label{17ottobre2020-2}
 \vec{\psi}_k(\lambda,u~|\nu)= \Gamma(\lp_k+1)\vec{e}_k+
   \sum_{l=1}^\infty \vec{b}_l^{~(k)}(u) (\lambda-u_k)^l ,\quad \hbox{ for } \lambda\to u_i,
   \ee
     the Taylor expansion being uniformly convergent  and the coefficients $\vec{b}_l^{~(k)}(u)$ holomorphic on  $ \mathbb{D}(u^c)$.     \end{itemize}

      Let $i,j$ be such that $u_i^c=u_j^c$. Then $\vec{\Psi}_i^{(sing)}(\lambda,u~|\nu)$ and $\vec{\Psi}_j^{(sing)}(\lambda,u~|\nu)$ are either linearly independent, or at least one of them is identically zero (identity to zero  can be realized only for    $\lp_i\leq -2 $ or   $\lp_j\leq -2$.)

   \eth
   \begin{proof} See Section \ref{19novembre2020-4}. \end{proof}
 
 \bre 
 {\rm 
  Of the coefficients of \eqref{17ottobre2020-2}, only $b_0^{(i)}(u)$, $b_1^{(i)}(u)$, ..., $b_{\lp_k}^{(k)}(u)$ will be useful later. 
 }
 \ere
 
 \bre
 \label{7aprile2021-1} 
 {\rm
For $\lp_k\not \in \mathbb{Z}_{-}$, the singular solution $\vec{\Psi}_k^{(sing)}$  is unique, identified by its singular behaviour at $\lambda=u_k$ and   the normalization  \eqref{31agosto2020-1}-\eqref{19novembre2020-1} when $\lp_k\in\mathbb{C}\backslash  \mathbb{Z}$, or
 by the normalization  \eqref{17ottobre2020-2} when $\lp_k\in \mathbb{N}$.  For $\lp_k\in \mathbb{Z}_{-}$,  a singular 
 solution  in \eqref{18ottobre2020-2}  is not  unique, but its singular behaviour \eqref{18ottobre2020-3} at $\lambda=u_k$
  is uniquely fixed by the normalization \eqref{31agosto2020-1}-\eqref{19novembre2020-1}. There is a freedom  
  due to the choice of the coefficients $r_m$ and of $\vec{\phi}_k$ in \eqref{18ottobre2020-2}. See also Remark \ref{19novembre2020-6}. }
 \ere
 
   The singular behaviour of $\vec{\Psi}_k$   at $\lambda=u_j$ is expressed by connection coefficients.  
 
  \begin{definition}
  \label{18ottobre2020-5}
    The connection coefficients are defined by 
   \be
   \label{31agosto2020-5}
\vec{\Psi}_k(\lambda,u~|\nu)\underset{\lambda\to u_j}=\vec{\Psi}_j^{(sing)}(\lambda,u~|\nu) ~c_{jk}^{(\nu)}~+\hbox{\rm reg}(\lambda-u_j),\quad\quad \lambda\in\mathcal{P}_\eta, 
\ee
and by 
\be
\label{9ottobre2020-1}
c_{jk}^{(\nu)}:=0,~\forall k=1,...,n,\quad \hbox{when $\vec{\Psi}^{(sing)}_j\equiv 0$, possibly occurring for $\lp_j\in -\mathbb{N}-2$}.
\ee
\end{definition}

The uniqueness of the singular behaviour of  $\vec{\Psi}_j^{(sing)}$ at $\lambda=u_j$ implies  that the $c_{jk}$ are {\it uniquely defined}. 
From the definition, we  see that 

$\bullet$ If $\lp_k\not\in\mathbb{Z}$, $c^{(\nu)}_{kk}=1$.

$\bullet$ If $\lp_k\in\mathbb{Z}$, $c^{(\nu)}_{kk}=0$.

$\bullet$ If $\lp_k\in\mathbb{N}$ and $\vec{\Psi}_k(\lambda,u~|\nu)\equiv 0$, then $c^{(\nu)}_{1k}=c^{(\nu)}_{2k}=\dots=c^{(\nu)}_{nk}=0$.

$\bullet$ If $\lp_j\in-\mathbb{N}-2$ and $\vec{\Psi}_j^{(sing)}(\lambda,u~|\nu)\equiv 0$, then $c^{(\nu)}_{j1}=c^{(\nu)}_{j2}=\dots=c^{(\nu)}_{jn}=0$.

\bpr
\label{13ottobre2020-1}
 The coefficients in \eqref{31agosto2020-5}-\eqref{9ottobre2020-1} are 
 {\bf isomonodromic connection coefficients}, namely they are independent of $u\in\mathbb{D}(u^c)$. They  satisfy the vanishing relations 
\be
\label{31agosto2020-4}
c_{jk}^{(\nu)}=0 \quad \hbox{ for $j\neq k$ such that } u^c_j=u^c_k.
\ee
\epr
   
 \begin{proof} See Section \ref{22novembre2020-7}.  \end{proof}

 \section{Proof of Theorem \ref{30agosto2020-5}}
 \label{19novembre2020-4}

 {\bf Remark on notations:} Throughout this section, we  work with functions $f=f(\lambda,u|~\nu)$ defined on $\mathcal{P}_\eta(u)\hat{\times}\mathbb{D}(u^c)$. For simplicity we omit  $\nu$ and write $f=f(\lambda,u)$. Similarly, we write   $c_{jk}$ in place of $c_{jk}^{(\nu)}$.  
 
 \subsection{Fundamental matrix solution of the Pfaffian System}
 
    Without loss of generality, we order the eigenvalues so  that\footnote{ In this way,  $\mathbb{D}(u^c)=\mathbb{D}_1^{\times p_1}\times\cdots \times \mathbb{D}_s^{\times p_s}$, where $\mathbb{D}_\alpha=\{x \in\mathbb{C}~|~ |x-\lambda_\alpha|\leq \epsilon_0\}$, $\alpha=1,...,s$.} 
\begin{align}
\label{2settembre2020-8}
&  u_1^c=\dots =u_{p_1}^c=\lambda_1;
\quad
  u^c_{p_1+1}=\dots=u^c_{p_1+p_2}=\lambda_2;
\\
\noalign{\medskip}
  \label{2settembre2020-9}
&u^c_{p_1+p_2+1}=\dots=u^c_{p_1+p_2+p_3}=\lambda_3; \quad ..... \hbox{ up to }  \quad  u^c_{p_1+\dots+p_{s-1}+1}=\dots=u^c_{p_1+\dots+p_{s-1}+p_s}=\lambda_s.
  \end{align}
  We  analyse first the coalescence of  $u_1$, ..., $u_{p_1}$ to $\lambda_1$. Other cases are analogous.   
   We change variables $(u_1,...,u_n,\lambda)\mapsto (x_1,...,x_{n+1})$  as follows
  $$ 
  x_{n+1}=\lambda-\lambda_1,\quad\quad 
  x_j=\left\{\begin{array}{cc}
  \lambda-u_j,& 1\leq j \leq p_1;
  \\
 \noalign{\medskip}
  u_j-\lambda_1,& p_1+1\leq j \leq n.
  \end{array}
  \right.
  $$
  The inverse transformation is 
  $$
  \lambda=x_{n+1}+\lambda_1,\quad\quad u_j=\left\{\begin{array}{c}
  x_{n+1}-x_j+\lambda_1, \quad 1\leq j \leq p_1,
  \\
  \noalign{\medskip}
  x_j+\lambda_1, \quad \quad \quad p_1+1\leq j \leq n.
  \end{array}
  \right.
    $$ 
Let 
    $$
  x:=(\underbrace{x_1,...,x_{p_1}}_{p_1},
        \underbrace{x_{p_1+1},....,x_n}_{n-p_1},~x_{n+1}) \equiv  (\underbrace{x_1,...,x_{p_1}}_{p_1},
         \boldsymbol{x}^\prime,x_{n+1}),
        $$
        where 
        $
  \boldsymbol{x}^\prime:=(x_{p_1+1},....,x_n)
  $. 
  We are interested in the behaviour of solutions for 
   $$ 
 x \longrightarrow (\underbrace{0,~0,...,0}_{p_1},~\boldsymbol{x}^\prime,0),  $$
 corresponding to 
 $$
 u_1\to \lambda_1,~\dots,~u_{p_1}\to \lambda_1,\quad\hbox{ and } \lambda\to \lambda_1$$
 namely  $u_i-u_j\to 0$ , $i\neq j$  and $\lambda-u_i\to 0$, for  $i, j\in \{1,...,p_1\}$.    
  The Pfaffian system \eqref{19agosto2020-5} in variables $x$, with Fuchsian singularities at  $x_1=0, \dots,x_{p_1}=0$, becomes
\be
\label{13agosto2020-3}
 d \Psi=P(x)\Psi,\quad \quad  P(x)=\sum_{j=1}^{p_1} \frac{P_j(x)}{x_j} dx_j + \sum_{j=p_1+1}^{n+1} \widehat{P}_j(x)dx_j
\ee
where 
   \begin{align*}
  & \frac{P_j(x)}{x_j}=\frac{B_j(x)}{x_j}-\omega_j(x),\quad 1\leq j \leq p_1,
 &  
 \widehat{P}_j(x)= \frac{B_j(x)}{x_j-x_{n+1}}+\omega_j(x),\quad  p_1+1 \leq j \leq n,
  \\
 &  \widehat{P}_{n+1}(x)=\sum_{j=p_1+1}^n\frac{B_j(x)}{x_{n+1}-x_j}+\sum_{j=1}^{p_1}\omega_j(x)
\end{align*}

\noindent
The Pfaffian system  is assumed integrable with  holomorphic in   $\mathbb{D}(u^c)$ coefficients, therefore 
 $P_1(x),...,P_{p_1}(x)$ and $\widehat{P}_{p_1+1}(x),...,\widehat{P}_{n+1}(x)$  are holomorphic at     $(\underbrace{0,\dots,0}_{p_1},~\boldsymbol{x}^\prime,0)$, for  $\boldsymbol{x}^\prime$ varying  as $u_{p_1+1},\dots,u_n$ vary in $\mathbb{D}(u^c)$.

 \bre
 \label{7aprile2021-2}
 {\rm  The  commutation relations \eqref{12agosto2020-1} at 
 $u=(\underbrace{\lambda_1,\dots,\lambda_1}_{p_1},\boldsymbol{u}^\prime)$, where  $\boldsymbol{u}^\prime:=(u_{p_1+1},\dots,u_n)$, 
  are 
 \be
\label{11agosto2020-3}
 [B_i(\lambda_1,\dots,\lambda_1,\boldsymbol{u}^\prime),B_j(\lambda_1,\dots,\lambda_1,\boldsymbol{u}^\prime)]=0,\quad1\leq i \neq j \leq p_1.
 \ee 
 They also follow from  the integrability condition  
 $ 
 dP(x)=P(x)\wedge P(x)$  
 of \eqref{13agosto2020-3}, which implies 
 $$\frac{\partial }{\partial x_i} \left(\frac{P_j}{x_j}\right)-\frac{\partial }{\partial x_j} \left(\frac{P_i}{x_i}\right)-\frac{P_iP_j-P_j P_i}{x_ix_j}=0,\quad\quad 1\leq i\neq j \leq p_1.
 $$
 Let $\hat{\boldsymbol{k}}=(k_1,...,k_{p_1})$, and write $\hat{\boldsymbol{l}}\leq \hat{\boldsymbol{k}}$ if $k_i\leq l_i$ for all $i\in \{1,...,p_1\}$. 
The Taylor convergent series 
  $ 
 P_i(x)= \sum_{k_1+\cdots +k_{p_1}\geq 0} P_{i,\hat{\boldsymbol{k}}}(\boldsymbol{x}^\prime,x_{n+1}) x_1^{k_1}\cdots x_{p_1}^{k_{p_1}},
 $
 has  coefficients $P_{i,\widehat{\boldsymbol{k}}}(\boldsymbol{x}^\prime,x_{n+1})$
 holomorphic of $\boldsymbol{x}^\prime$, $x_{n+1}$. 
 The integrability condition becomes  \cite{YT}
\be
\label{11agosto2020-2}
 k_jP_{i,\hat{\boldsymbol{k}}}-k_iP_{j,\hat{\boldsymbol{k}}}+\sum_{\boldsymbol{0}\leq \hat{\boldsymbol{l}}\leq \hat{\boldsymbol{k}}}[P_{i,\hat{\boldsymbol{l}}},P_{j,\hat{\boldsymbol{k}}-\hat{\boldsymbol{l}}}]=0,\quad 1\leq i\neq j \leq p_1.
 \ee
 In particular,
 $ 
 P_{i,\hat{\boldsymbol{0}}}(\boldsymbol{x}^\prime,x_{n+1})=B_i(\underbrace{\lambda_1,\dots,\lambda_1}_{p_1},\boldsymbol{u}^\prime)
$ for $\hat{\boldsymbol{k}}=\hat{\boldsymbol{0}}$, 
 so that \eqref{11agosto2020-2} reduces to \eqref{11agosto2020-3}. $\Box$
}
\ere
 
  Let us define Jordan matrices 
 \be
 \label{19settembre2020-4}
\widehat{T}^{(j)}=\hbox{\rm diag}(0,\dots,0,\underbrace{{-1-\lp_j}}_{\hbox{position $j$}},0,\dots,0), \quad \quad \hbox{ for  $\lp_j\neq -1$}.
\ee
 \be
 \label{19settembre2020-4-bis}
\widehat{T}^{(j)}:=J^{(j)}:=\begin{pmatrix} 
0  &&0~\cdots&&0
\\
\vdots &\ddots&& & \vdots 
\\
\noalign{\medskip}
0  &  ~\cdots~& 0~\cdots&r_{m_j}^{(j)}&0 
\\
\noalign{\medskip}
\vdots &&&\ddots&  \vdots 
\\
0 &&0~\cdots&&0
\end{pmatrix}\quad \longleftarrow \hbox{ row $j$}, \quad \hbox{for $\lp_j=-1$},
\ee 
$$ 
 r_{m_j}^{(j)}:=1,\quad \hbox{ is the only non-zero entry in position $(j,m_j)$, with  $m_j\geq p_1+1$}.
 $$  
   
 \ble
 \label{17settembre2020-1}
 Under the assumptions of Theorem \ref{30agosto2020-5}, for every $j\in\{1,...,n\}$ there exists a holomorphically   invertible matrix $G^{(j)}(u)$ on $\mathbb{D}(u^c)$ reducing $B_j(u)$ to constant Jordan form. 
Moreover, $B_1(u^c),...,B_{p_1}(u^c)$  are simultaneously reducible to   $\widehat{T}^{(1)}, ..., \widehat{T}^{(p_1)}$ respectively. 
 \ele

\vskip 0.3 cm 
\begin{proof} 
 For every $j\in\{1,...,n\}$, the Schlesinger system  \eqref{10agosto2020-3}-\eqref{10agosto2020-5} implies the  Frobenius integrability (on  $\mathbb{D}(u^c)$) of the  the linear Pfaffian system (see  Corollary \ref{18novembre2020-2},  Appendix A)
\be
  \label{18novembre2020-1} 
\frac{\partial G^{(j)}}{\partial u_k} =\left(\frac{B_k}{u_k-u_j}+\gamma_k\right)G^{(j)}, \quad k\neq j,\quad\quad
 \frac{\partial G^{(j)}}{\partial u_j}=-\sum_{k\neq j} \left(\frac{B_k}{u_k-u_j}+\gamma_k\right) G^{(j)}
\ee
From \eqref{10agosto2020-4}-\eqref{10agosto2020-5} and the above  we receive   $ 
\partial_k \bigl((G^{(j)})^{-1} B_j G^{(j)}\bigr)=0$, $ k=1,...,n$,  for a  holomorphic on $\mathbb{D}(u^c) $ fundamental matrix solution $G^{(j)}(u)$. 
 Thus, up to $G^{(j)}\mapsto G^{(j)}\mathcal{G}^{(j)}$, $\mathcal{G}^{(j)}\in GL(n,\mathbb{C})$, we can choose   $G^{(j)}(u)$ which  puts $B_j$ in constant Jordan form.
   If we consider each $B_j$ separately,  now for $j\in \{1,...,p_1\}$, it is straightforward that the Jordan forms are the matrices  $\widehat{T}^{(j)}$.\footnote{It is also elementary to find a holomorphic $G^{(j)}$ explicitly. For example, if all $B_j(u)$ are diagonalizable (i.e $\lp_j\neq -1$), an elementary computation shows that
$(G^{(j)}(u))^{-1} B_j(u) G^{(j)}(u)=\hbox{\rm diag}(0,\dots,0,-1-\lp_j,0,\dots,0)$, $j=1,2,...,n,$, where the columns of $G^{(j)}$ are as follows: 
$$ \hbox{$j$-th column  is multiple of  } \vec{e}_j\in\mathbb{C}^n;\quad\quad 
\hbox{$l$-th column, $l\neq j$,   is multiple of } \vec{e_l}-\frac{A_{jl}(u)}{\lp_j+1}\vec{e_j}.
$$
}
 An elementary computation  shows that $B_1(u^c),...,B_{p_1}(u^c)$ are actually  reducible to $\widehat{T}^{(1)},....,\widehat{T}^{(1)}$ simultaneously,\footnote{For example, in case of the previous footnote, the simultaneous reduction to Jordan form  at $u^*=(\lambda_1,...,\lambda_1, \boldsymbol{u}^\prime)$, where $\boldsymbol{u}^\prime=(u_{p_1+1},...,u_n)$), is realized by  the product $G^{(1)}(u^*)\cdots G^{(p_1)}(u^*)$, which depends holomorphically on $\boldsymbol{u}^\prime$} because   only the $j$-th row of $B_j(u^c)$ is non-zero, and  by \eqref{12agosto2020-1}  the first $p_1$ entries of this row are  zero, except for the $(j,j)$-entry equal to  $-\lp_j-1$. Namely,
 $$ 
B_j(u^c)=\begin{pmatrix} 
0  ~0&&&& \cdots & 0 
\\
\vdots &&&& & \vdots 
\\
{\bf 0 ~0}&-\lp_j-1 & {\bf 0 ~0}&~ -A_{j,p_1+1}^{(j)}(u^c)& \cdots & -A_{j,n}(u^c) 
\\
\vdots &&&&&  
\\
0~0 &&&& \cdots & 0
\end{pmatrix}\quad \longleftarrow \hbox{ row $j$}.
$$
\end{proof}
 
 \bre
 {\rm 
 As in Lemma \ref{17settembre2020-1},   $B_1(u^c),...,B_{p_1}(u^c)$ are reducible simultaneously to their respective Jordan forms,  $B_{p_1+1}(u^c),...,B_{p_1+p_2}(u^c)$ are reducible simultaneously to their  respective Jordan forms, and so on up to $B_{p_1+...+p_{s-1}+1}(u^c),...,B_{p_1+...+p_s}(u^c)$.
 }
 \ere
  

For short, let $\boldsymbol{p}_1:=(1,...,p_1)$. Without loss of generality, we   label $u_1,...,u_{p_1}$ so that 
$$ 
\lp_j\in\mathbb{C}\backslash\mathbb{Z}, \quad \hbox{ for } 1\leq j \leq q_1, \quad\quad \lp_j\in\mathbb{Z}, \quad \hbox{ for } q_1+1\leq j \leq p_1.
$$
If all $\lp_j\in\mathbb{Z}$, then $q_1=0$,  if all $\lp_j\not \in\mathbb{Z}$, then $q_1=p_1$.    The first and fundamental  step to achieve Theorem  \ref{30agosto2020-5}  is the following

\bth
 \label{30agosto2020-9}
 In the assumptions of Theorem \ref{30agosto2020-5}, the Pfaffian system \eqref{19agosto2020-5} admits the fundamental matrix solution 
 \be
 \label{19settembre2020-3-old}
    \Psi^{(\boldsymbol{p}_1)}(\lambda,u)=G^{(\boldsymbol{p}_1)} U^{(\boldsymbol{p}_1)}(\lambda,u)\cdot  \prod_{l=1}^{p_1} (\lambda-u_l)^{\widehat{T}^{(l)}}\cdot\prod_{j=q_1+1}^{p_1}(\lambda-u_j)^{\widehat{R}^{(j)}},\quad (\lambda,u)\in \mathcal{P}(u)\hat{\times}\mathbb{D}(u^c),
   \ee
    where  $G^{(\boldsymbol{p}_1)}$ is a constant invertible matrix simultaneously reducing $B_1(u^c), ..., B_{p_1}(u^c)$ to $\widehat{T}^{(1)}, ..., \widehat{T}^{(p_1)}$ as  in \eqref{19settembre2020-4}-\eqref{19settembre2020-4-bis}. The matrix function  $U^{(\boldsymbol{p}_1)}(\lambda,u)$ is  holomorphic in $\mathbb{D}_1\times \mathbb{D}(u^c)$ with convergent expansion 
    \begin{align*}
    &
   U^{(\boldsymbol{p}_1)}(\lambda,u)=  I+
   \\
   &
   +\sum_{\boldsymbol{k}>0,~k_{1}+...+k_{p_1}\geq 0}\Bigl[U^{(\boldsymbol{p}_1)}_{\boldsymbol{k}} \cdot 
 (u_{p_1+1}-u_{p_1+1}^c)^{k_{p_1+1}}\cdots (u_n-u_n^c)^{k_n} (\lambda-\lambda_1)^{k_{n+1}}\Bigr]  (\lambda-u_1)^{k_1}~\cdots ~(\lambda-u_{p_1})^{k_{p_1}},
    \end{align*}
and constant matrix coefficient $U^{(\boldsymbol{p}_1)}_{\boldsymbol{k}}$.   Here $\boldsymbol{k}:=(k_1,...,k_n,k_{n+1})$, $k_j\geq 0$, and $\boldsymbol{k}>0$ means that at least one $k_j>0$ ($j=1,...,n+1$).  The  exponents $\widehat{R}^{(q_1+1)},\dots,
 \widehat{R}^{(p_1)}$ are constant  nilpotent matrices. 
 \begin{itemize} 
 \item If $\lp_j=-1$,  
 \be
 \label{7ottobre2020-2}
 \widehat{R}^{(j)}=0.
 \ee 
 \item If $\lp_j\in \mathbb{N}=\{0,1,2,...\}$, only the entries $\widehat{R}^{(j)}_{mj}=:r_m^{(j)} $, for $m=1,...,n$ and $m \neq j$,  are possibly  non zero, namely 
 \be
 \label{28settembre2020-6}
 \widehat{R}^{(j)}=\left[
 \vec{0}~\left|~\cdots~\left|~\vec{0}~\left|~ \sum_{m\neq j, m=1}^n r^{(j)}_m \vec{e}_m ~\right|~\vec{0}~\right|~\cdots~\right|~\vec{0}
 \right],
 \ee
  where  only the $j$-th column is possibly non-zero.
  
 \item If $\lp_j\in -\mathbb{N}-2=\{-2,-3,...\}$, only the entries  $\widehat{R}^{(j)}_{jm}=:r_m^{(j)}$,  for $m=1,...,n$ and $m\neq j$,  are possibly non zero, namely
 \be
 \label{28settembre2020-7}
\widehat{R}^{(j)}=\begin{pmatrix} 
0 &\cdots &&&& \cdots & 0 
\\
\vdots &&&&& & \vdots 
\\
r_1^{(j)} & \cdots & r_{j-1}^{(j)} & 0 & r_{j+1}^{(j)}& \cdots & r_n^{(j)}
\\
\vdots &&&&& & \vdots 
\\
0 & \cdots&&&& \cdots & 0
\end{pmatrix}\quad \longleftarrow \hbox{ row $j$ is possibly non zero }. 
\ee
 
 \end{itemize}

The exponents 
 $\widehat{T}^{(l)}$ and $R^{(j)}$ satisfy the following commutation relations 
 \begin{align}
 \label{17settembre2020-3-bis-old}
& [\widehat{T}^{(i)},\widehat{T}^{(j)}]=0,   \quad i,j=1,...,p_1;
\\
\noalign{\medskip}
 \label{17settembre2020-3-old}
  & [\widehat{R}^{(j)},\widehat{R}^{(k)}]=0, \quad  [\widehat{T}^{(i)},\widehat{R}^{(j)}]=0, \quad  i=1,...,p_1,\quad i\neq j,\quad j,k=q_1+1,...,p_1.
 \end{align}
  By analytic continuation, $\Psi^{(\boldsymbol{p}_1)}(\lambda,u)$ defines an analyic function on the universal covering of $\mathcal{P}_\eta(u)\hat{\times} \mathbb{D}(u^c)$.  Another representation of \eqref{19settembre2020-3-old} will be given in \eqref{19settembre2020-3-new}.
\eth

\begin{proof}

 We  apply the results  of  \cite{YT}  at the point $x= x^c:=(\underbrace{0,~0,...,0}_{p_1},~\boldsymbol{x}_c^\prime,0)$,  with $  \boldsymbol{x}_c^\prime:=(x^c_{p_1+1},....,x_n^c)
 $, 
 corresponding to $u=u^c$ and $\lambda=\lambda_1$, 
 where
 $ 
 x^c_j= u_j^c-\lambda_1$, $ j=p_1+1, ..., n$. 
 By Theorem 7 of \cite{YT}, the Pfaffian system \eqref{13agosto2020-3} admits a fundamental matrix solution
 \be
 \label{20novembre2020-1}
 \Psi^{(\boldsymbol{p}_1)}(\lambda,u)=U_0~ U(x)~  Z(x),\quad \quad Z(x)= \prod_{j=1}^{p_1} x_l^{A_j}  \prod_{j=1}^{p_1} x_l^{Q_j} ,\quad \quad \det U_0\neq 0, 
 \ee
for certain   matrices   $A_j$ which are simultaneous  triangular forms of $B_1(u^c),..., B_{p_1}(u^c)$. While in  \cite{YT} a lower triangular form is considered, we equivalently use   the upper triangular one.  The matrices $Q_j$ will be described  below.  The matrix   $U(x)=V(x)\cdot W(x)$ has structure
  \begin{align*}
 &V(x)=
  I+\sum_{\boldsymbol{k}>0,~k_{p_1+1}+...+k_{n+1}>0}V_{\boldsymbol{k}}~
 x_1^{k_1}~\cdots ~x_{p_1}^{k_{p_1}}~(x_{p_1+1}-x^c_{p_1+1})^{k_{p_1+1}}\cdots (x_n-x^c_n)^{k_n}\cdot x_{n+1}^{k_{n+1}}
 \\
& W(x)=I+\sum_{k_{1}+...+k_{p_1}>0}W_{k_1,...,k_{p_1}}~
 x_1^{k_1}~\cdots ~x_{p_1}^{k_{p_1}}.
    \end{align*}
        The constant matrix coefficients $V_{\boldsymbol{k}}$, $W_{k_1,...,k_{p_1}}$ can be determined \cite{YT} from the constant  matrix coefficients $ P_{i,\boldsymbol{k}}$ in the Taylor expansion\footnote{   
          $$
 P_i(x)= \sum_{k_1+\cdots +k_{n+1}\geq 0} P_{i,\boldsymbol{k}} ~x_1^{k_1}\cdots  x_{p_1}^{k_{p_1}} \cdot (x_{p_1+1}-x_{p_1+1}^c)^{k_{p_1+1}}\cdots  (x_{n}-x_{n}^c)^{k_{n}}\cdot x_{n+1}^{k_{n+1}} .
 $$
  and analogous for $\widehat{P}_j(x)$} of the  $P_j(x)$ and $\widehat{P}_j(x)$. 
  Recall that 
    $ x_j= \lambda-u_j$, $1\leq j\leq p_1$, and $ x_{n+1}=\lambda-\lambda_1$. 
    Moreover, for $p_1+1\leq j \leq n$, we have 
    $ 
    x_j-x_j^c= (u_j-\lambda_1)-(u_j^c-\lambda_1)=u_j-u_j^c$. 
    Thus, restoring variables $(\lambda,u)$, we have 
    \begin{align*}
    &V(\lambda,u)=I+
   \\
   &
    +\sum_{
    k_{p_1+1}+...+k_{n+1}>0
   }\Bigl[V_{\boldsymbol{k}}
 (u_{p_1+1}-u_{p_1+1}^c)^{k_{p_1+1}}\cdot ...\cdot (u_n-u_n^c)^{k_n}\cdot (\lambda-\lambda_1)^{k_{n+1}}\Bigr]  (\lambda-u_1)^{k_1}~\cdots ~(\lambda-u_{p_1})^{k_{p_1}},
 \\
 \noalign{\medskip}
 & W(\lambda, u_1,...,u_{p_1})=I+\sum_{k_{1}+...+k_{p_1}>0}W_{k_1,...,k_{p_1}}~
(\lambda-u_1)^{k_1}~\cdot...\cdot ~(\lambda-u_{p_1})^{k_{p_1}}.
    \end{align*}
    Therefore, the matrices  appearing in the statement are $G^{(\boldsymbol{p}_1)}:=U_0$ and  $
  U^{(\boldsymbol{p}_1)}(\lambda,u):=  V(\lambda,u)W(\lambda,u)$, which is  holomorphic for $ (\lambda, u)\in\mathbb{D}_1\times \mathbb{D}(u^c)$. 

We show that the exponents $A_j$ and $Q_j$ are  respectively $\widehat{T}^{(j)}$ in  \eqref{19settembre2020-4}-\eqref{19settembre2020-4-bis} and   $\widehat{R}^{(j)}$  in \eqref{7ottobre2020-2}-\eqref{28settembre2020-6}-\eqref{28settembre2020-7}.     
According to \cite{YT}  (see theorems 2 and 5), the matrix function   $G^{(\boldsymbol{p}_1)}\cdot U^{(\boldsymbol{p}_1)}(\lambda,u)$ in \eqref{19settembre2020-3-old} provides the gauge transformation   
     $$ 
 \Psi=G^{(\boldsymbol{p}_1)}\cdot U^{(\boldsymbol{p}_1)}(\lambda,u)Z\underset{\hbox{in notation of \cite{YT}}}\equiv U_0U(x) Z,
 $$ which brings \eqref{13agosto2020-3} to the {\it reduced form} (being "reduced" is defined in  \cite{YT}) 
$$
dZ=\sum_{j=1}^{p_1}\frac{Q_j(x)}{x_j}~Z,\quad \quad  
Q_j(x)=A_j +\sum_{\widehat{\boldsymbol{k}}>0} Q_{\widehat{\boldsymbol{k}},j} x_1^{k_1}\cdots x_{p_1}^{k_{p_1}},
$$
where the notation $\widehat{\boldsymbol{k}}=(k_1,...,k_{p_1})>0$ means at least one $k_l>0$.
From \cite{YT}, we have the following.

 \noindent
$\bullet$  The $A_j$   are simultaneous triangular forms of $B_1(u^c),...,B_{p_1}(u^c)$.  Thus, by Lemma \ref{17settembre2020-1}, they can be taken to be
  $$ 
 A_j =\widehat{T}^{(j)} \hbox{ as in \eqref{19settembre2020-4}-\eqref{19settembre2020-4-bis},  $j=1,...,p_1$}. 
 $$
 \noindent
 $\bullet$ The   $Q_{\widehat{\boldsymbol{k}},j}$  satisfy  diag$(Q_{\widehat{\boldsymbol{k}},j})=0$, while   the  entry $(\alpha,\beta)$ for  $\alpha\neq \beta$ satisfies 
 $$ 
(Q_{\widehat{\boldsymbol{k}},j})_{\alpha\beta} \neq 0 \quad\hbox{ only if }  \quad (\widehat{T}^{(j)})_{\alpha\alpha}-(\widehat{T}^{(j)})_{\beta\beta}= k_j\geq 0,\quad \hbox{ for all $j=1,...,p_1$}.
$$
Taking into account the particular structure   \eqref{19settembre2020-4}-\eqref{19settembre2020-4-bis}, the above condition can be satisfied only for
$$
 \widehat{\boldsymbol{k}}=(\underbrace{0,...,0}_{q_1},\underbrace{0,...,0,k_j,0,...,0}_{p_1-q_1}),\quad k_j=|\lp_j+1|\geq 1 \hbox{ in position $j$},
 $$
because
\begin{align}
\label{20settembre2020-6}
(\widehat{T}^{(j)})_{\alpha\alpha}-(\widehat{T}^{(j)})_{\beta\beta}&=-\lp_j-1\geq 1 \quad \hbox{ when } \lp_j\in-\mathbb{N}-2\quad \hbox{ and } \alpha=j \quad (\beta\neq j),
\\ 
\noalign{\medskip}
\label{20settembre2020-7}
(\widehat{T}^{(j)})_{\alpha\alpha}-(\widehat{T}^{(j)})_{\beta\beta}&=\lp_j+1\geq 1 \quad \hbox{ when } \lp_j\in\mathbb{N}\quad \hbox{ and } \beta=j \quad (\alpha \neq j).
\end{align}
This can occur only for $j=q_1+1,...,p_1$.  Thus
\be
\label{7ottobre2020-1}
Q_{\widehat{\boldsymbol{k}},j}=0,\quad j=1,...,q_1,\quad\quad Q_{\widehat{\boldsymbol{k}},j}= \widehat{R}^{(j)} \hbox{ in  \eqref{7ottobre2020-2}-\eqref{28settembre2020-6}-\eqref{28settembre2020-7}},\quad j=q_1+1,...,p_1.
\ee
 In conclusion, the  reduced form turns out to be 
\be
\label{20settembre2020-3}
 dZ=\left[\sum_{j=1}^{p_1}\left(\frac{\widehat{T}^{(j)}+\widehat{R}^{(j)} x^{k_j}}{x_j}\right)\right]Z,\quad\quad \widehat{R}^{(1)}=\dots=\widehat{R}^{(q_1)}=0.
 \ee

Its integrability implies  the  commutation relations.  Indeed, the compatibility   $\partial_i\partial_j Z=\partial_j\partial_i Z$, $i\neq j$, holds  if and only if 
$$
\frac{[\widehat{T}^{(j)},\widehat{T}^{(i)}]}{x_ix_j}+[\widehat{R}^{(j)},\widehat{R}^{(i)}]x_i^{k_i-1}x_j^{k_j-1} +[\widehat{T}^{(j)},\widehat{R}^{(i)}]x_i^{k_i-2}+   [\widehat{R}^{(j)},\widehat{T}^{(i)}]x_j^{k_j-2}=0,
\quad\quad 1\leq i\neq j \leq p_1.
$$
Keeping into account that $\widehat{R}^{(1)}=\dots=\widehat{R}^{(q_1)}=0$, the above holds if and only if \eqref{17settembre2020-3-bis-old}-\eqref{17settembre2020-3-old} hold. 

\vskip 0.2 cm 
The last to be checked is that a fundamental matrix  of \eqref{20settembre2020-3} is $Z(x)$ in \eqref{20novembre2020-1}, namely 
$$ 
Z(x)= \prod_{l=1}^{p_1} x_l^{\widehat{T}^{(l)}}  \prod_{j=q_1+1}^{p_1} x_l^{\widehat{R}^{(j)}}.
$$ 
 It suffices to verify this by differentiating $Z(x)$, keeping into account  the commutation relations \eqref{17settembre2020-3-bis-old}-\eqref{17settembre2020-3-old}  and the formula $\partial_i x_i^M=(M/x_i)x_i^M$, for a constant matrix $M$.  
For $i=1,...,q_1$ we receive
$$\frac{\partial }{\partial x_i} Z(x) = \frac{\widehat{T}^{(i)}}{x_i} Z(x).
$$
For $i=q_1+1,...,p_1$ we receive
\begin{align*}
\frac{\partial }{\partial x_i} Z(x) &= \frac{T^{(i)}}{x_i} Z(x)+ \Bigl(\prod_{l=1}^{p_1} x_l^{\widehat{T}^{(l)}} \Bigr)\frac{\widehat{R}^{(i)}}{x_i} \Bigl(\prod_{j=q_1+1}^{p_1} x_l^{\widehat{R}^{(j)}}\Bigr)
\\
& 
=\frac{\widehat{T}^{(i)}}{x_i} Z(x)+  \Bigl(\prod_{l=1}^{i-1} x_l^{\widehat{T}^{(l)}} \Bigr) \frac{x_i^{\widehat{T}^{(i)}}\widehat{R}^{(i)}}{x_i}  \Bigl( \prod_{l=i+1}^{p_1} x_l^{\widehat{T}^{(l)}} \Bigr) \Bigl( \prod_{j=q_1+1}^{p_1} x_l^{\widehat{R}^{(j)}} \Bigr)=(**).
\end{align*}
Now, recalling that $k_i=|\lp_i+1|$ and \eqref{20settembre2020-6}-\eqref{20settembre2020-7}, we see that 
$
x_i^{\widehat{T}^{(i)}}\widehat{R}^{(i)}x_i^{-\widehat{T}^{(i)}}= \widehat{R}^{(i)}x_i^{k_i} 
$.
Therefore, 
$$
(**)=\frac{\widehat{T}^{(i)}}{x_i} Z(x)+ \frac{\widehat{R}^{(i)} x_i^{k_i}}{x_i}  \Bigl( \prod_{l=1}^{p_1} x_l^{\widehat{T}^{(l)}} \Bigr) \Bigl( \prod_{j=q_1+1}^{p_1} x_l^{\widehat{R}^{(j)}} \Bigr)= \frac{\widehat{T}^{(i)}+\widehat{R}^{(i)} x_i^{k_i}}{x_i} Z(x),
$$ as we wanted to prove.

Finally, the fact that $\Psi^{(\boldsymbol{p}_1)}(\lambda,u)$ has analytic continuation on the universal covering of $\mathcal{P}_\eta(u)\hat{\times} \mathbb{D}(u^c)$ follows from general results in the theory of linear Pfaffian systems \cite{Haraoka,IKSY,YT}. 

\end{proof}


It is convenient to introduce a slight change of  the exponents.  
Without loss in generality, we can label $u_1,...,u_{p_1}$ in such a way that, for some $q_1,c_1\geq 0$ integers,  the following ordering of eigenvalues of $A$ holds:
$$
\underline{\lp_1,~\dots,~\lp_{q_1}\in\mathbb{C}\backslash\mathbb{Z},\quad
\quad
 \lp_{q_1+1},~\dots,~\lp_{q_1+c_1}\in \mathbb{Z}_{-},
\quad\quad
\lp_{q_1+c_1+1},~\dots,~\lp_{p_1}\in \mathbb{N}.}
$$
Clearly, $0\leq q_1\leq p_1$, $0\leq c_1\leq p_1$ and $0\leq q_1+c_1\leq p_1$.  We define new exponents. 
\begin{itemize}
\item For $\lp_j\neq -1$, 
\be 
\label{28settembre2020-5-bis}
T^{(j)}:=\widehat{T}^{(j)},\quad j=1,...,p_1;\quad\quad\quad  R^{(j)}:=\widehat{R}^{(j)},\quad j=q_1+1,...,p_1.
\ee 
\item For $\lp_j= -1$ (so $j\in \{q_1+1,...,q_1+c_1\}$),
\be
\label{28settembre2020-5}
T^{(j)}:=0 ,\quad R^{(j)}:=\underbrace{J^{(j)}}_{\hbox{in \eqref{19settembre2020-4-bis}}}=\begin{pmatrix} 
0  &&0~\cdots&&0
\\
\vdots &\ddots&& & \vdots 
\\
\noalign{\medskip}
0  &  ~\cdots~& 0~\cdots&r_{m_j}^{(j)}&0 
\\
\noalign{\medskip}
\vdots &&&\ddots&  \vdots 
\\
0 &&0~\cdots&&0
\end{pmatrix} \longleftarrow \hbox{ row $j$}, \quad r_{m_j}^{(j)}=1.
\ee
Recall that $m_j\geq p_1+1$. 
\end{itemize}
This  new definitions  allow to treat together  the case  $\lp_j\in-\mathbb{N}-2$ and the case $\lp_j=-1$. 

\ble With the definition \eqref{28settembre2020-5-bis}-\eqref{28settembre2020-5}, the following relations hold.
 \begin{align}
 \label{17settembre2020-3-bis}
 & [T^{(i)},T^{(j)}]=0,   \quad i,j=1,...,p_1;
\\
\noalign{\medskip}
 \label{17settembre2020-3}
 & [R^{(j)},R^{(k)}]=0, \quad  [T^{(i)},R^{(j)}]=0,  
  \quad i=1,...,p_1,\quad i\neq j,\quad j,k=q_1+1,...,p_1,
 \end{align}
  \ele
  \begin{proof}
The equivalence between \eqref{17settembre2020-3-bis-old}-\eqref{17settembre2020-3-old} and \eqref{17settembre2020-3-bis}-\eqref{17settembre2020-3} is straightforward. 
\end{proof}

 \bcr 
 In Theorem \ref{30agosto2020-9}, the fundamental matrix solution  \eqref{19settembre2020-3-old} is 
  \be
 \label{19settembre2020-3-new}
    \Psi^{(\boldsymbol{p}_1)}(\lambda,u)=G^{(\boldsymbol{p}_1)}\cdot U^{(\boldsymbol{p}_1)}(\lambda,u)\cdot  \prod_{l=1}^{p_1} (\lambda-u_l)^{T^{(l)}}\cdot\prod_{j=q_1+1}^{p_1}(\lambda-u_j)^{R^{(j)}},
   \ee
   where the exponents are defined in \eqref{28settembre2020-5-bis}-\eqref{28settembre2020-5}.
   \ecr
   \begin{proof}
   It is an immediate consequence of the commutation relations being satisfied, that the representation \eqref{19settembre2020-3-old} for $\Psi^{(\boldsymbol{p}_1)}$ still holds with  the definition  \eqref{28settembre2020-5-bis}-\eqref{28settembre2020-5}. 
   \end{proof}

The commutation relations impose  a simplification on the structure of the matrices $R^{(j)}$. Let the new convention \eqref{28settembre2020-5-bis}-\eqref{28settembre2020-5} be used.  The relations  $[T^{(i)},R^{(j)}]=0$ for $i=1,...,p_1$ and $j=q_1+1,...,p_1$, $j\neq i$, imply the vanishing of the first $p_1$ non-trivial entries of $R^{(j)}$, so that  (by \eqref{28settembre2020-6}, \eqref{28settembre2020-7} and \eqref{28settembre2020-5}),
\be
\label{7ottobre2020-11}
 R^{(j)}=\left[
 \vec{0}~\left|~\cdots~\left|~\vec{0}~\left|~ \sum_{m=p_1+1}^n r^{(j)}_m \vec{e}_m ~\right|~\vec{0}~\right|~\cdots~\right|~\vec{0}
 \right],\quad\quad \lp_j\in\mathbb{N}.
 \ee
\be
\label{7ottobre2020-10} 
R^{(j)}=\begin{pmatrix} 
0 &\cdots &&&& \cdots & 0 
\\
\vdots &&&&& & \vdots 
\\
0 & \cdots &0 & 0 & r_{p_1+1}^{(j)}& \cdots & r_n^{(j)}
\\
\vdots &&&&& & \vdots 
\\
0 & \cdots&&&& \cdots & 0
\end{pmatrix}\quad \longleftarrow \hbox{ row $j$},\quad \quad    \lp_j\in \mathbb{Z}_{-};
\ee
 The relations $[R^{(j)},R^{(k)}]=0$ for either $j,k\in\{q_1+1,\dots,q_1+c_1\}$ or  $j,k\in \{q_1+c_1+1,\dots,p_1\}$ are automatically satisfied.  On the other hand, the commutators $[R^{(j)},R^{(k)}]=0$ for $j\in\{q_1+1,\dots,q_1+c_1\}$ and $k\in \{q_1+c_1+1,\dots,p_1\}$  imply the  further (quadratic) relations 
 \be
 \label{29settembre2020-1}
\sum_{m=p_1+1}^n r^{(j)}_m r^{(k)}_m=0.
\ee
In particular, if $\lp_j=-1$ and $R^{(j)}$ is \eqref{28settembre2020-5}, all the above conditions can be satisfied, provided that we take  $m_j\geq p_1+1$, as we have agreed from the beginning.

  \subsection{Selected Vector Solutions $\vec{\Psi}_i$}
  \label{11ottobre2020-2}
  
  \noindent
  {\bf Remark on notations.} For the sake of the proof, it is convenient to use a slightly different notation with respect to the statement of Theorem \ref{30agosto2020-5}. The identifications between objects in the proof and  objects in the statement is  $ 
\vec{\varphi}_i \longmapsto \vec{\psi}_i$, $r_i^{(m)}/r_k^{(i)}\longmapsto  r_m$ and $\vec{\varphi}_k/r_k^{(i)} \longmapsto \phi_i$. 

\vskip 0.3 cm 
We will construct  selected vector solutions of Theorem  \ref{30agosto2020-5} from suitable    linear combinations of columns of the fundamental matrix $\Psi^{(\boldsymbol{p}_1)}$ in \eqref{19settembre2020-3-new}.  The $i$-th column of an $n\times n$ matrix $M$ is $M\cdot \vec{e}_i$ (rows by columns multiplication), where $\vec{e}_i$ is the standard unit basic vector in $\mathbb{C}^n$.
From  \eqref{17settembre2020-3-bis}-\eqref{17settembre2020-3},    and \eqref{7ottobre2020-10}-\eqref{7ottobre2020-11}-\eqref{29settembre2020-1}, we receive 
  $$ 
\prod_{l=1}^{p_1} (\lambda-u_l)^{T^{(l)}}\cdot\prod_{j=q_1+1}^{p_1}(\lambda-u_j)^{R^{(j)}}\cdot \vec{e}_i=
$$
\be
\label{27settembre2020-1}
=
\left\{
\begin{array}{ccc}
(\lambda-u_i)^{-\lp_i-1} \vec{e}_i, & i=1,...,q_1+c_1, &  \lp_i\in\mathbb{C}\backslash\mathbb{N};
\\
\noalign{\medskip}
(\lambda-u_i)^{-\lp_i-1} \vec{e}_i+\left(\sum_{m=p_1+1}^n r^{(i)}_m \vec{e}_m \right) \ln(\lambda-u_i), & i=q_1+c_1+1,...,p_1,& \lp_i\in\mathbb{N};
\\
\noalign{\medskip}
\vec{e}_i+\sum_{m=q_1+1}^{q_1+c_1} \vec{e}_m r_i^{(m)} (\lambda-u_m)^{-\lp_m-1}\ln(\lambda-u_m), & i=p_1+1,...,n.
\end{array}
\right.
\ee
For $i=1,...,n$, let 
\be
\label{27settembre2020-2}
\vec{\varphi}_i(\lambda,u):=G^{(\boldsymbol{p}_1)} U(\lambda,u) \cdot \vec{e}_i,\quad \quad i=1,...,n,
\ee  
which is holomorphic for $(\lambda,u)\in\mathbb{D}_1\times \mathbb{D}(u^c)$. 
For  $i=1,...,p_1$, we define vector valued functions 
\be
\label{27settembre2020-3}
\vec{\Psi}_i(\lambda,u):=
\left\{
\begin{array}{ccc}
\vec{\varphi}_i(\lambda,u) (\lambda-u_i)^{-\lp_i-1},& i=1,...,q_1+c_1,& \lp_i\in\mathbb{C}\backslash\mathbb{N};
\\
\noalign{\medskip}
\sum_{k=p_1+1}^n r_k^{(i)} \vec{\varphi}_k(\lambda,u), & i=q_1+c_1+1,...,p_1, & \lp_i\in\mathbb{N}.
\end{array}
\right.
\ee
 Notice that for  $i=q_1+c_1+1,...,p_1$,   if   $r_k^{(i)}=0$ for all $k=p_1+1,...,n$, then $\vec{\Psi}_i(\lambda,u)$ is identically zero  
\be
\label{19novembre2020-3} 
\vec{\Psi}_i(\lambda,u)\equiv 0,\quad \lp_i\in\mathbb{N}, 
\ee
Hence,  the $i$-th column of $\Psi^{(\boldsymbol{p}_1)}(\lambda,u)$ is
\begin{align}
\label{25settembre2020-1}
& &\Psi^{(\boldsymbol{p}_1)}(\lambda,u)\cdot\vec{e}_i &=\vec{\Psi}_i(\lambda,u), &&i=1,...,q_1+c_1,
\\
\label{25settembre2020-2}
&&&
=\vec{\Psi}_i(\lambda,u)\ln(\lambda-u_i) +\frac{\vec{\varphi}_i(\lambda,u)}{(\lambda-u_i)^{\lp_i+1}}
, && i=q_1+c_1+1,...,p_1,
\\
\label{25settembre2020-3}
&&&=\varphi_i(\lambda,u)+\sum_{m=q_1+1}^{q_1+c_1} r_i^{(m)}\vec{\Psi}_m(\lambda,u) \ln (\lambda-u_m),&& i=p_1+1,...,n.
\end{align}

\bpr {\it  The vector functions \eqref{27settembre2020-3} coincide with  the following linear  combinations of columns of $\Psi^{(\boldsymbol{p}_1)}(\lambda,u)$, 
\be
\label{11ottobre2020-1}
\vec{\Psi}_i(\lambda,u)=
\left\{
\begin{array}{cc}
\Psi^{(\boldsymbol{p}_1)}(\lambda,u)\cdot \vec{e}_i, & i=1,...,q_1+c_1, \quad \hbox{ namely } \lp_i\in
 \mathbb{C}\backslash\mathbb{N};
 \\
 \noalign{\medskip}
\Psi^{(\boldsymbol{p}_1)}(\lambda,u)\cdot \sum_{k=p_1+1}^n r_k^{(i)} \vec{e}_k, & i=q_1+c_1+1,...,p_1, \quad \hbox{ namely } \lp_i\in
 \mathbb{N}.
 \end{array}
\right.
\ee
 As such, they are vector solutions (called {\bf selected}) of the Pfaffian system \eqref{19agosto2020-5}}. Those  $\vec{\Psi}_i(\lambda,u)$ which are not identically zero are linearly independent.   
\epr

\begin{proof} 
For  $i=1,...,q_1+c_1$,  \eqref{11ottobre2020-1} is just  \eqref {25settembre2020-1}, so it is a vector solution of  \eqref{19agosto2020-5}. In case $i=q_1+c_1+1,...,p_1$, we  claim that $\vec{\Psi}_i(\lambda,u)$  defined in \eqref{27settembre2020-3} coincides with   the following  linear combination 
$$
\vec{\Psi}_i(\lambda,u)= \sum_{k=p_1+1}^n r_k^{(i)}\left(
\Psi^{(\boldsymbol{p}_1)}(\lambda,u)\cdot\vec{e}_k\right), \quad  \quad i=q_1+c_1+1,...,p_1,
$$
of the  vector solutions \eqref{25settembre2020-3}. Indeed, 
\begin{align*}
\sum_{k=p_1+1}^n r_k^{(i)}\left(\Psi^{(\boldsymbol{p}_1)}(\lambda,u)\cdot\vec{e}_k\right)&= \sum_{k=p_1+1}^n r_k^{(i)}\left(
\varphi_k(\lambda,u)+\sum_{m=q_1+1}^{q_1+c_1} r_k^{(m)}\vec{\Psi}_m(\lambda,u) \ln (\lambda-u_m)\right)
\\
&\underset{\eqref{27settembre2020-3}}=\vec{\Psi}_i(\lambda,u) +\sum_{m=q_1+1}^{q_1+c_1}\left(\sum_{k=p_1+1}^n  r_k^{(i)}r_k^{(m)}\right)\vec{\Psi}_m(\lambda,u)\ln(\lambda-u_m).
\end{align*}
Now, it follows from \eqref{29settembre2020-1} that $\sum_{k=p_1+1}^n  r_k^{(i)}r_k^{(m)}=0$, so proving the claim and the expressions \eqref{11ottobre2020-1}. Linear independence follows from \eqref{11ottobre2020-1}. 
\end{proof}

 \subsection{Singular Solutions $\vec{\Psi}_i^{(sing)}$}
 \label{11ottobre2020-3}
 Using the previous results, we define  singular  vector solutions of the Pfaffian system.
 \begin{itemize}
 \item For \underline{$\lp_i\not\in\mathbb{Z}$}, i.e. $i=1,...,q_1$,
$$ \vec{\Psi}_i^{(sing)}(\lambda,u):= \vec{\Psi}_i(\lambda,u)~\equiv \Psi^{(\boldsymbol{p}_1)}(\lambda,u)\cdot \vec{e}_i
$$

\item For \underline{$\lp_i\in\mathbb{N}$}, i.e. $i=q_1+c_1+1,...,p_1$,
$$
\vec{\Psi}_i^{(sing)}(\lambda,u):= \vec{\Psi}_i(\lambda,u)\ln(\lambda-u_i)+\frac{\vec{\varphi}_i(\lambda,u)}{(\lambda-u_i)^{\lp_i+1}}~\equiv \Psi^{(\boldsymbol{p}_1)}(\lambda,u)\cdot \vec{e}_i.
$$

\item For \underline{$\lp_i\in \mathbb{Z}_{-}$}, i.e. $i=q_1+1,...,q_1+c_1$, we distinguish three subcases. 

\begin{itemize} 
\item[i)] If \underline{$\lp_i\leq -2$} and $r^{(i)}_k\neq 0$ for some $k\in \{p_1+1,...,n\}$, from \eqref{25settembre2020-3} (change notation $i\mapsto k$)
$$
\vec{\Psi}_i^{(sing)}(\lambda,u):=\frac{1}{r^{(i)}_k}\left\{\varphi_k(\lambda,u)+\sum_{m=q_1+1}^{q_1+c_1} r_k^{(m)}\vec{\Psi}_m(\lambda,u) \ln (\lambda-u_m)\right\}\equiv \frac{1}{r^{(i)}_k} \Psi^{(\boldsymbol{p}_1)}(\lambda,u)\cdot \vec{e}_k.
$$
\item[ii)]  If \underline{$\lp_i\leq -2$} and $r^{(i)}_k= 0$ for all $k\in \{p_1+1,...,n\}$,
$$ \vec{\Psi}_i^{(sing)}(\lambda,u):=0
$$
\item[iii)] If \underline{$\lp_i=-1$}, then  $r_{m_i}^{(i)}=1$ and in  i) above  we take  $k=m_i$, so that  
\begin{align*}
\vec{\Psi}_i^{(sing)}(\lambda,u)&:=\vec{\varphi}_{m_i}(\lambda,u)+ \vec{\Psi}_i(\lambda,u)\ln(\lambda-u_i)+\sum_{m\neq i, ~m=q_1+1}^{q_1+c_1} r_{m_i}^{(m)} \vec{\Psi}_m(\lambda,u) \ln(\lambda-u_m).
\\
& 
=\Psi^{(\boldsymbol{p}_1)}(\lambda,u)\cdot \vec{e}_{m_i},\quad \quad m_i\geq p_1+1. 
\end{align*}
\end{itemize}
The above $\vec{\Psi}_i^{(sing)}(\lambda,u)$ in i) and iii)  is singular at $u_i$, but possibly  also at $u_{q_1+1}, \dots, u_{q_1+c_1}$ corresponding to $\lp_m\in \mathbb{Z}_{- }$. By definition,  
 \be
 \label{25settembre2020-5}
  \vec{\Psi}_i^{(sing)}(\lambda,u)\underset{\lambda\to u_i}=\vec{\Psi}_i(\lambda,u) \ln (\lambda-u_i) +\hbox{\rm reg}(\lambda-u_i),\quad\quad i=q_1+1,...,q_1+c_1,
  \ee   
  
  \bre
  \label{19novembre2020-6}
  {\rm 
The definition  in i) contains the freedom of choosing  $k\in \{p_1+1,...,n\}$, which changes $\varphi_k(\lambda,u)$ and  the ratios ${r_k^{(m)}/r^{(i)}_k}$ (in formula \eqref{18ottobre2020-2}, $\varphi_k/r^{(i)}_k$ is  denoted by $\phi_i$ and  ${r_k^{(m)}/r^{(i)}_k}$ is $r_m$). Whatever is the choice of $k$, provided that  $r^{(i)}_k\neq 0$, the  behaviour at $\lambda=u_i$  of the corresponding $\vec{\Psi}_i^{(sing)}$ is always  \eqref{25settembre2020-5}, so it is uniquely fixed if we fix the normalization of $ \vec{\Psi}_i(\lambda,u)$. 
 }\ere

\end{itemize}

As a consequence of the above definitions  and Section \ref{11ottobre2020-2}, we receive  the  following 
\bpr
The $\vec{\Psi}_i^{(sing)}(\lambda,u)$ defined above, $i=1,...,p_1$, when not identically zero, are linearly independent. They  are represented as follows
$$ 
\vec{\Psi}_i^{(sing)}(\lambda,u)
=
\left\{
\begin{array}{ccc}
\Psi^{(\boldsymbol{p}_1)}(\lambda,u)\cdot \vec{e}_i, & \lp_i\in \mathbb{C}\backslash \mathbb{Z}_{-},&
\\
\noalign{\medskip}
\Psi^{(\boldsymbol{p}_1)}(\lambda,u)\cdot \dfrac{\vec{e}_k}{r_k^{(i)}}, & \lp_i\in \mathbb{Z}_{-},&\hbox{for some $k\in\{p_1+1,...,n\}$ such that $r_k^{(i)}\neq 0$}
\\
0, & \lp_i\in -\mathbb{N}-2,&\hbox{if $r_k^{(i)}= 0$ for all $k\in\{p_1+1,...,n\}$.}
\end{array}
\right.
$$ 
\epr
 
 \vskip 0.2 cm 
  
 \subsection{Expansions at $\lambda=u_i$, $i=1,...,p_1$ and completion of the proof.}

In order to proceed in the proof, and in view of the Laplace transform to come,  we need local behaviour at $\lambda=u_i$. 

\ble 
The following   Taylor expansion holds at $\lambda=u_i$, with coefficients $\vec{d}_l^{~(i)}(u)$ holomorphic on $\mathbb{D}(u^c)$,
$$ 
\vec{\Psi}_i(\lambda,u)
=
\sum_{l=0}^\infty \vec{d}_l^{~(i)}(u) (\lambda-u_i)^l,
\quad\quad \lp_i \in \mathbb{N}, ~\hbox{ namely }~ i=q_1+c_1+1,...,p_1.
$$ 

\ele 
\begin{proof}
By  \eqref{27settembre2020-3},  $\vec{\Psi}_i(\lambda,u)= 
G^{(\boldsymbol{p}_1)} U(\lambda,u) \cdot (\sum_{m=p_1+1}^n r^{(i)}_m \vec{e}_m )$, so it is holomorphic on $\mathbb{D}_1\times \mathbb{D}(u^c)$. From this we conclude.   \end{proof}

The coefficients $d_l^{(i)}(u)$ will be fixed by a chosen normalization for $\vec{\varphi}_i$  in \eqref{25settembre2020-2}, as in the following lemma. 

\ble The following   Taylor expansions hold at $\lambda=u_i$, uniformly convergent  for $u\in \mathbb{D}(u^c)$. 
$$
\left.
\begin{array}{cc} 
\hbox{$\lp_i\not \in \mathbb{N}$, i.e. $i=1,...,q_1+c_1$:} & \vec{\Psi}_i(\lambda,u) 
  \\
  \\
  \hbox{$\lp_i \in \mathbb{N}$, i.e.  $q_1+c_1+1,...,p_1$:} &\dfrac{\vec{\varphi}_i(\lambda,u)}{(\lambda-u_i)^{\lp_i+1}}
  \end{array}
\right\}
\underset{\lambda\to u_i}= \Bigl(f_i\vec{e}_i+\sum_{l=1}^\infty \vec{b}_l^{~(i)}(u)(\lambda-u_i)^l\Bigr)(\lambda-u_i)^{-\lp_i-1},
$$
with certain vector coefficients $\vec{b}_l^{~(i)}(u)$ holomorphic in $\mathbb{D}(u^c)$. In particular, the leading term is  {\rm constant}, and will be chosen as follows
\be
\label{26novembre2020-1}
f_i=
\left\{
\begin{array}{ccc} \Gamma(\lp_i+1),&  \lp_i\in\mathbb{C}\backslash \mathbb{Z},&i=1,...,q_1,
\\
\noalign{\medskip}
\dfrac{(-1)^{\lp_i}}{(-\lp_i-1)!}, & \lp_i \in \mathbb{Z}_{-},& i=q_1+1,...,q_1+c_1,
\\
\noalign{\medskip}
\lp_i!\equiv \Gamma(\lp_i+1), & \lp_i\in\mathbb{N}, & i=q_1+c_1+1,...,p_1.
\end{array}
\right.
\ee
 \ele
 
 \begin{proof}
   That the above convergent expansions must hold follows from the definitions. Work is required to prove that the leading term is $f_i\vec{e}_i$, with  $f_i\in\mathbb{C}\backslash\{0\}$. 
   
    From  definitions \eqref{27settembre2020-2}-\eqref{27settembre2020-3},  the leading term must  coincide with the leading term  of the expansion at $\lambda=u_i$ of the $i$-th   column $G^{(\boldsymbol{p}_1)} U(\lambda,u) \cdot \vec{e}_i$, for $i=1,...,p_1$. 
 To evaluate it,  observe that the solution   $\Psi^{(\boldsymbol{p}_1)}(\lambda,u)$,   restricted to  a  polydisc $\mathbb{D}(u^0)$  contained in a $\tau$-cell of $\mathbb{D}(u^c)$, is a fundamental matrix solution of the Fuchsian system \eqref{03} in the {\rm Levelt form} \eqref{28settembre2020-1} below at $\lambda=u_i$, $i=1,...,p_1$. 
 Indeed, by \eqref{17settembre2020-3} it can be written as 
$$ 
\Psi^{(\boldsymbol{p}_1)} (\lambda,u)=\Bigl\{G^{(\boldsymbol{p}_1)} U^{(\boldsymbol{p}_1)} (\lambda,u)  \prod_{
\begin{array}{c} l=1 \\ l\neq i\end{array}}^{p_1}
(\lambda-u_l)^{T^{(l)}}
\prod_{
\begin{array}{c} j=q_1+1 \\ j\neq i\end{array}}^{p_1}(\lambda-u_j)^{R^{(j)}}\Bigr\}~\cdot (\lambda-u_i)^{T^{(i)}} (\lambda-u_i)^{R^{(i)}},
$$
where it is understood that $R^{(i)}=0$ if $i=1,...,q_1$. We have
$$ U^{(\boldsymbol{p}_1)} (\lambda,u)=I+F_i(u)+O(\lambda-u_i), \quad \quad \lambda\to u_i, \quad\quad F_{i}(u):=U^{(\boldsymbol{p}_1)} (u_i,u),
$$ 
and $O(\lambda-u_i)$ represent vanishing terms  at $\lambda=u_i$, holomorphic in $\mathbb{D}_1\times \mathbb{D}(u^c)$. The expansion   at $\lambda=u_i$ of the factors $(\lambda-u_l)^{T^{(i)}}$ and $(\lambda-u_j)^{R^{(j)}}$, for  $l,j\neq i$, yields the 
 {\it Levelt form}
\be
\label{28settembre2020-1}
\Psi^{(\boldsymbol{p}_1)}(\lambda,u)\underset{\lambda\to u_i}= G^{(i;\boldsymbol{p}_1)}(u)\Bigl(I+O(\lambda-u_i)\Bigr)(\lambda-u_i)^{T^{(i)}} (\lambda-u_i)^{R^{(i)}},\quad i=1,...,p_1,
\ee
where $O(\lambda-u_i)$ are higher order terms, provided that $u\in \mathbb{D}(u^0)$ (they contain negative  powers $(u_i-u_k)^{-m}$), and
$$
G^{(i;\boldsymbol{p}_1)}(u):= G^{(\boldsymbol{p}_1)}(I+F_i(u)) \prod_{
\begin{array}{c} l=1 \\ l\neq i\end{array}}^{p_1}
(u_i-u_l)^{T^{(l)}}
\prod_{
\begin{array}{c} j=q_1+1 \\ j\neq i\end{array}}^{p_1}(u_i-u_j)^{R^{(j)}},\quad i=1,...,p_1
.
$$
The matrix $G^{(i;\boldsymbol{p}_1)}(u)$ is holomorphically invertible if restricted to a polydisc $\mathbb{D}(u^0)$ contained in a $\tau$-cell, but it is branched at the coalescence locus $\Delta$ on the whole $\mathbb{D}(u^c)$.

We reach our goal if we show that the $i$-th column  $G^{(i;\boldsymbol{p}_1)}(u)\cdot \vec{e}_i$ is constant in  $\mathbb{D}(u^c)$. First,  it follows from \eqref{28settembre2020-1}  and the standard isomonodromic theory of \cite{JMU} that $G^{(i;\boldsymbol{p}_1)}(u)$    holomorphically  in $\mathbb{D}(u^0)$  reduces $B_i(u)$   to the diagonal  form $T^{(i)}$, when $\lp_i\neq -1$,
$$ 
\Bigl(G^{(i;\boldsymbol{p}_1)}(u)\Bigr)^{-1} B_i(u) ~G^{(i;\boldsymbol{p}_1)}(u)=T^{(i)},
$$
or to non-diagonal Jordan form \eqref{28settembre2020-5}  when  $\lp_i =-1$
$$ 
\Bigl(G^{(i;\boldsymbol{p}_1)}(u)\Bigr)^{-1} B_i(u) ~G^{(i;\boldsymbol{p}_1)}(u)=R^{(i)}\equiv J^{(i)},\quad\quad \lp_i =-1.
$$
For this reason, the $i$-th row  is proportional to the eigenvector $\vec{e}_i$ of $B_i(u)$ relative to the eigenvalue $-\lp_i-1$. Namely, for some scalar function $f_i(u)$,
$$ 
   G^{(i;\boldsymbol{p}_1)}(u)\cdot \vec{e}_i= f_i(u) \vec{e}_i. 
	$$
	This is obvious for $\lp_i\neq -1$, namely for diagonalizable $B_i$. If $\lp_i=-1$, 
the eigenvalue $0$ of $B_i$ appearing in $J^{(i)}$ at entry $(i,i)$ is associated with the eigenvector 
$f_i(u) \vec{e}_i$. Moreover, for every invertible matrix $G=[*|\cdots|*|\vec{e}_i|*|\cdots|*]$, where  $\vec{e}_i$ 
 occupies the $k$-th column, then $G^{-1}B_i(u) G$ is zero eveywhere, except for the $k$-th row.
Now, since $R^{(i)}=J^{(i)}$ has only one non-zero entry
   on the $i$-th row,  it follows that  the eigenvector   $f_i(u) \vec{e}_i$ must occupy the $i$-th column of $G^{(i;\boldsymbol{p}_1)}(u)$.

\vskip 0.2 cm 
$\bullet$ {\it $f_i(u)$ is  holomorphic on $\mathbb{D}(u^c)$}. Indeed, by \eqref{27settembre2020-1}, 
$$
 \prod_{
\begin{array}{c} l=1 \\ l\neq i\end{array}}^{p_1}
(u_i-u_l)^{T^{(l)}}
\prod_{
\begin{array}{c} j=q_1+1 \\ j\neq i\end{array}}^{p_1}(u_i-u_j)^{R^{(j)}}\cdot  \vec{e}_i=\vec{e}_i. 
$$
Therefore 
$ 
f_i(u) \vec{e}_i \equiv G^{(i;\boldsymbol{p}_1)}(u)\cdot \vec{e}_i ~ = G^{(\boldsymbol{p}_1)}(I+F_i(u))\vec{e}_i
$. We conclude, because   $F_i(u)$ is holomorphic on $\mathbb{D}(u^c)$. 

\vskip 0.2 cm 
$\bullet$ {\it $f_i$ is constant} on $\mathbb{D}(u^c)$.  Indeed, since $\Psi^{(\boldsymbol{p}_1)}(\lambda,u)$  is an isomonodromic solution in $\mathbb{D}(u^0)$, the matrix $G^{(i;\boldsymbol{p}_1)}(u)$ must satisfy the Pfaffian system  (see Appendix A,  identify $G^{(i;\boldsymbol{p}_1)}$ with $G^{(i)}$ in Corollary \ref{18novembre2020-2})
\be
\label{22agosto2020-7}
   \dfrac{\partial  G^{(i;\boldsymbol{p}_1)}}{\partial u_j}=\left( \dfrac{B_j}{u_j-u_i} +\omega_j\right) G^{(i;\boldsymbol{p}_1)}, \quad j\neq i 
  ;
  \quad\quad
    \dfrac{\partial  G^{(i;\boldsymbol{p}_1)}}{\partial u_i}=\sum_{j\neq i} \left( \dfrac{B_j}{u_i-u_j}+\omega_j\right) G^{(i;\boldsymbol{p}_1)} . 
    \ee
   From  \eqref{22agosto2020-6} and \eqref{22agosto2020-5}, the $i$-th column of 
    $
     \dfrac{B_j}{u_j-u_i} +\omega_j$ is null. Hence,     $$ 
      \frac{\partial}{\partial u_j}\left(G^{(i;\boldsymbol{p}_1)}\cdot \vec{e}_i\right)=0, \quad \forall j\neq i .
      $$ 
      Moreover, summing the equations \eqref{22agosto2020-7}, we get $
      \sum_{j=1}^n \partial_j  G^{(i;\boldsymbol{p}_1)}=0$.  
     Thus,   $G^{(i;\boldsymbol{p}_1)}\cdot \vec{e}_i$ is constant on $\mathbb{D}(u^0)$, and being holomorphic on $\mathbb{D}(u^c)$, it is constant on  $\mathbb{D}(u^c)$.  
     The choice  \eqref{26novembre2020-1} will be made. 
\end{proof}

The above obtained expansions for the $\vec{\Psi}_i$ and $\vec{\Psi}_i^{(sing)}$ and $\vec{\varphi}_i$  prove Theorem \ref{30agosto2020-5}   for $i=1,...,p_1$, with some obvious identifications between objects in the proof and  objects in the statement, namely $ 
\vec{\varphi}_i \longmapsto \vec{\psi}_i$, $r_i^{(m)}/r_k^{(i)}\longmapsto  r_m$ and $\vec{\varphi}_k/r_k^{(i)} \longmapsto \phi_i$.

\subsection{Analogous proof for all coalescences}

 With the labelling \eqref{2settembre2020-8}-\eqref{2settembre2020-9}, the same strategy above  holds for every coalescence $$
 (u_{p_1+...+p_{\alpha-1}+1}, ..., u_{p_1+...+p_{\alpha}})\longrightarrow (\lambda_\alpha,...,\lambda_\alpha), \quad \alpha=1,...,s.
 $$
We find corresponding isomondromic fundamental matrices  for the Pfaffian system (with self-explaining notations)
 $$
  \Psi^{(\boldsymbol{p}_\alpha)}(\lambda,u)
 =G^{(\boldsymbol{p}_\alpha)}\cdot U^{(\boldsymbol{p}_\alpha)}(\lambda,u)\cdot \prod_{l=p_1+...+p_{\alpha-1}+1}^{p_1+...+p_{\alpha}}
  (\lambda-u_l)^{T^{(l)}} \prod_{j=(p_1+...+p_{\alpha-1}+1)+q_\alpha}^{p_1+...+p_{\alpha}}(\lambda-u_j)^{R^{(j)}}.
 $$
 where 
 $\boldsymbol{p}_\alpha=(p_1+...+p_{\alpha-1}+1,\dots,p_1+...+p_{\alpha})$. 
 Then, we proceed in the same way, constructing  the  solutions $\vec{\Psi}_i$ and  $\vec{\Psi}_i^{(sing)}$, with $p_1+...+p_{\alpha-1}+1\leq i \leq p_1+...+p_{\alpha}$. $\Box$

\subsection{Proof of Proposition  \ref{13ottobre2020-1}}
\label{22novembre2020-7}
\begin{proof} 
For simplicity, we omit $\nu$ in the connection coefficients $c_{jk}^{(\nu)}$  in \eqref{31agosto2020-5}-\eqref{9ottobre2020-1}.  
It follows from the very definitions of the $\vec{\Psi}_k$ and $\vec{\Psi}_j^{(sing)}$ that 
$$ 
c_{jk}=0\quad \hbox{ if $u_j^c=u_k^c$}.
$$ 

In order to prove independence of $u$, we express  the monodromy of  $$ 
\Psi(\lambda,u):=[\vec{\Psi}_1(\lambda,u)~|~\cdots~|\vec{\Psi}_n(\lambda,u)],
$$  in terms of the connection coefficients. From the definition, we have (using  the notations in the statement of Theorem \ref{30agosto2020-5})
\be
\label{11ottobre2020-5}
 \vec{\Psi}_k(\lambda,u) =
\left\{
\begin{array}{cc}
\vPsi_j(\lambda,u)c_{jk}+\hbox{reg}(\lambda-u_j),& \lp_j\not\in \mathbb{ Z}
\\
\\
 \vPsi_j(\lambda,u)\ln(\lambda-u_j)c_{jk}+\hbox{reg}(\lambda-u_j),& \lp_j\in \mathbb{Z}_{-}
\\
\\
 \left(\vPsi_j(\lambda,u)\ln(\lambda-u_j)+\dfrac{\psi_j(\lambda,u)}{ (\lambda-u_j)^{\lp_j+1}}\right)c_{jk}+\hbox{reg}(\lambda-u_j),& \lp_j\in \mathbb{ N}
\end{array}
\right.
\ee
For $u\not \in \Delta$ and a small loop $(\lambda-u_k)\mapsto (\lambda-u_k)e^{2\pi i}$ we obtain from Theorem \ref{30agosto2020-5}
$$ 
\vec{\Psi}_k(\lambda,u) \longmapsto \vec{\Psi}_k(\lambda,u) e^{-2\pi i \lp_k}, \quad \hbox{ which includes also the case $\lp_k\in \mathbb{Z}$, with  $e^{-2\pi i \lp_k}=1$}.
$$
For a small loop $(\lambda-u_j)\mapsto (\lambda-u_j)e^{2\pi i}$, $j\neq k$, from Theorem \ref{30agosto2020-5} and  \eqref{11ottobre2020-5}  we obtain  
\begin{align*}
&\vec{\Psi}_k \longmapsto \vec{\Psi}_j e^{-2\pi i \lp_j} c_{jk} +\underbrace{\hbox{\rm reg}(\lambda-u_j)}_{\vec{\Psi}_k-\vec{\Psi}_jc_{jk}}=\vec{\Psi}_k+(e^{-2\pi i \lp_j}-1)c_{jk} \vec{\Psi}_j & \hbox{for $\lp_j\not\in\mathbb{Z}$}
\\
& \vec{\Psi}_k \longmapsto \vPsi_j\Bigl(\ln(\lambda-u_j)+2\pi i\Bigr)c_{jk}+\hbox{reg}(\lambda-u_j)= \vec{\Psi}_k+2\pi i c_{jk} \vec{\Psi}_j, &\hbox{for $\lp_j\in\mathbb{Z}_{-}$}
\\
\noalign{\medskip}
& \vec{\Psi}_k \longmapsto \left(\vPsi_j\Bigl(\ln(\lambda-u_j)+2\pi i\Bigr)+\dfrac{\psi_j(\lambda,u)}{ (\lambda-u_j)^{\lp_j+1}}\right)c_{jk}+\hbox{reg}(\lambda-u_j)= \vec{\Psi}_k+ 2\pi i c_{jk} \vec{\Psi}_j, & \hbox{ for $\lp_j\in\mathbb{N}$}.
\end{align*}
Therefore, for $u\not \in \Delta$ and a small loop $\gamma_k:(\lambda-u_k)\mapsto (\lambda-u_k)e^{2\pi i}$  not encircling other points $u_j$ (we denote the loop by  $\lambda\mapsto \gamma_k\lambda$),  we receive 
$$ 
\Psi(\lambda,u)\longmapsto  \Psi(\gamma_k\lambda,u)=\Psi(\lambda,u) M_k(u),
$$ 
where
$$ 
(M_k)_{jj}=1\quad j\neq k,\quad  (M_k)_{kk}=e^{-2\pi i \lp_k};\quad \quad
(M_k)_{kj}= \alpha_kc_{kj},\quad j\neq k;\quad\quad (M_k)_{ij}=0 \hbox{ otherwise}.
$$
and 
           $$ 
\alpha_k:=(e^{-2\pi i \lp_k}-1), \quad \hbox{ if }\lp_k\not\in \mathbb{Z};\quad\quad\quad 
\alpha_k:=2\pi i ,\quad \hbox{ if }\lp_k\in \mathbb{Z}. 
$$
We proceed by first analyzing the generic case, and then the general case. 

{\bf Generic case.} Suppose that $A(u)$ has {\it no integer eigenvalues} (recall that eigenvalues do not depend on $u$).
Let us fix $u$ in a $\tau$-cell. By  Proposition \ref{26novembre2020-2},  $\Psi(\lambda,u)$ is a fundamental matrix solution of \eqref{03} for the fixed $u$, and $C=(c_{jk})$ is invertible.  Thus 
$$ 
M_k(u)= \Psi(\gamma_k\lambda,u)\Psi(\lambda,u)^{-1}.
$$ 
The above makes sense  for every $u$ in the considered $\tau$-cell, being $\Psi(\lambda,u)$  invertible at such an $u$. But $\Psi(\lambda,u)$ and $\Psi(\gamma_k\lambda,u)$ are holomorphic on $\mathcal{P}_\eta(u)\hat{\times} \mathbb{D}(u^c)$, so that  the matrix $M_k(u)$ is holomorphic on the $\tau$-cell. Repeating the above argument for another $\tau$-cell, we conclude that  $M_k(u)$ is holomorphic on each $\tau$-cell.  Now, on a $\tau$-cell, we have 
$$ 
d\Psi(\gamma_k\lambda,u)=P(\lambda,u)\Psi(\gamma_k\lambda,u)=P(\lambda,u)\Psi(\lambda,u) M_k,
$$ 
and at the same time 
$$ 
d\Psi(\gamma_k\lambda,u)=d\Bigl(\Psi(\lambda,u) M_k\Bigr)= d\Psi(\lambda,u)~M_k+\Psi(\lambda,u)~dM_k=
P(\lambda,u)\Psi(\lambda,u) M_k+\Psi(\lambda,u)~dM_k.
$$
The two expressions are equal if and only if $dM_k=0$, because $\Psi(\lambda,u)$ is invertible on a $\tau$-cell.  Recall that $\tau$-cells are disconnected from each other, so that  {\it separately on each cell, $M_k$ is constant, and so  the connection coefficients are constant separately on each cell}.

We further suppose that {\it none of the $\lp_j$ is integer}. In this case, $\vec{\Psi}^{(sing)}_j=\vec{\Psi}_j$ for all $j=1,...,n$, so that from \eqref{11ottobre2020-5} for $u_k^c\neq u_j^c$ (otherwise $c_{jk}=0$ and there is nothing to prove) 
\be
\label{11ottobre2020-9}
\vec{\Psi}_k(\lambda,u)\underset{\lambda\to u_j}= \vec{\Psi}_j(\lambda,u)c_{jk}+ \hbox{\rm reg}(\lambda-u_j). 
\ee
Using the labelling   \eqref{2settembre2020-8}-\eqref{2settembre2020-9}, from the proof of Theorem \ref{30agosto2020-9}  we have the fundamental matrix solution
 $$ 
 \Psi^{(\boldsymbol{p}_1)}(\lambda,u)=
 \Bigl[
 \vec{\Psi}_1(\lambda,u)~\Bigr|~\cdots~\Bigr|\vec{\Psi}_{p_1}(\lambda,u)~\Bigr|~\vec{\varphi}_{p_1+1}^{~(1)}(\lambda,u)~\Bigr|~\cdots ~\Bigr|~\vec{\varphi}_{n}^{~(1)}(\lambda,u)
 \Bigr]
$$
and in general  at each $\lambda_\alpha$, $\alpha=1,...,s$ (with $\sum_{j=1}^{\alpha-1}p_j=0$ for $\alpha=1$) we have 
\begin{align*} 
 \Psi^{(\boldsymbol{p}_\alpha)}(\lambda,u)=& \\
=& \Bigl[
\vec{\varphi}_{1}^{~(\alpha)}(\lambda,u)~\Bigl|~\cdots 
~\Bigl|~
\vec{\varphi}_{\sum_{j=1}^{\alpha-1}p_j}^{~(\alpha)}(\lambda,u)
~\Bigl|
 \vec{\Psi}_{\sum_{j=1}^{\alpha-1}p_j+1}(\lambda,u)~\Bigl|~ \vec{\Psi}_{\sum_{j=1}^{\alpha-1}p_j+2}(\lambda,u)~\Bigl|~\cdots~\Bigr|\vec{\Psi}_{\sum_{j=1}^{\alpha}p_j}(\lambda,u)~\Bigr|
 \\
 &\Bigr|~\vec{\varphi}_{\sum_{j=1}^{\alpha}p_j+1}^{~(\alpha)}(\lambda,u)~\Bigr|~\cdots ~|~\vec{\varphi}_{n}^{~(\alpha)}(\lambda,u)
 \Bigr]
\end{align*}
where 
$$
\vec{\Psi}_m(\lambda,u)=\vec{\psi}_m(\lambda,u)(\lambda-u_m)^{-\lp_m-1}, \quad m=\sum_{j=1}^{\alpha-1}p_j+1,~ \dots ~,\sum_{j=1}^{\alpha}p_j, 
$$
and the $\vec{\psi}_m(\lambda,u)$ and  $\vec{\varphi}_{r}^{~(\alpha)}(\lambda,u)$   are holomorphic functions  in the corresponding  $\mathbb{D}_\alpha\times \mathbb{D}(u^c)$. 
The above allows us to explicitly rewrite   \eqref{11ottobre2020-9}, for $j$ such that   $u_j^c=\lambda_\alpha$,  as 
\be
\label{23ottobre2020-1}
\vec{\Psi}_k(\lambda,u)=
\sum_{m=p_{1}+\cdots+p_{\alpha-1}+1}^{p_{1}+\cdots+p_{\alpha}} c_{mk}~ \vec{\Psi}_m(\lambda,u)+
\sum_{r\not \in \{p_{1}+\cdots+p_{\alpha-1}+1,...,p_{1}+\cdots+p_{\alpha}\}} h_r \vec{\varphi}_r^{~(\alpha)}(\lambda,u),
\ee
for suitable constant coefficients  $h_r$. Here one of the $c_{mk}$ is $c_{jk}$ of  \eqref{11ottobre2020-9}.

Each  $u_m$, with $m=p_{1}+\cdots+p_{\alpha-1}+1, ..., p_{1}+\cdots+p_{\alpha}$, varies in  $\mathbb{D}_\alpha$. 
Firstly, we can fix $\lambda=\lambda_\alpha$ in \eqref{23ottobre2020-1},  consider the branch cut $\mathcal{L}_\alpha$ from $\lambda_\alpha$ to infinity in direction $\eta$ 
(see Figure \ref{2settembre2020-10}), and let $u$ vary in such a way that 
each 
 $u_{p_{1}+\cdots+p_{\alpha-1}+1}$, ..., $u_{p_{1}+\cdots+p_{\alpha}}$
  varies in  $\mathbb{D}_\alpha\backslash \mathcal{L}_\alpha$, 
  so that in the r.h.s. of \eqref{23ottobre2020-1} all the 
  $\vec{\Psi}_m(\lambda_\alpha,u)$
   and  
   $\vec{\varphi}_r^{~(\alpha)}(\lambda_\alpha,u)$ 
    are holomorphic with respect to  $u$,
     provided that $u_m\neq \lambda_\alpha$.    
If $u$ varies, with the constraint that the $u_m$'s  must remain in $\mathbb{D}_\alpha\backslash \mathcal{L}_\alpha$,  every  $\tau$-cell of $\mathbb{D}(u^c)$ can be reached  starting from an initial point in one specific cell. This proves,
  by $u$-analytic continuation of  \eqref{23ottobre2020-1} with fixed  $\lambda=\lambda_\alpha$, that the coefficients $c_{mk}$ are 
  constant\footnote{Recall that $\mathbb{D}(u^c)=\bigtimes_{\beta=1}^s\mathbb{D}_\beta^{\times p_\beta}$.}  
  in $(\mathbb{D}_\alpha\backslash \mathcal{L}_\alpha)^{\times p_\alpha}\times\left( \bigtimes_{\beta\neq \alpha} \mathbb{D}_\beta^{\times p_\beta}\right)\subset \mathbb{D}(u^c)$.
  
   Now, we can slightly  vary $\eta$ in $\eta_{\nu+1}<\eta<\eta_\nu$, so that the cut $\mathcal{L}_\alpha$ is irrelevant\footnote{The crossing locus $X(\tau)$, $\tau=3\pi/2-\eta$, is as arbitrary as is the choice of $\tau$ in the range $\tau_\nu<\tau<\tau_{\nu+1}$.}. Thus,  the $c_{mk}$ are constant on 
   $\left\{ u\in \mathbb{D}(u^c) ~|~ u_{p_1+...+p_{\alpha-1}+1}\neq \lambda_\alpha,\dots, u_{p_1+...+p_{\alpha}}\neq \lambda_\alpha\right\}$. 

Finally, we  fix another value   $\lambda=\lambda^*\in \mathbb{D}_\alpha$ in \eqref{23ottobre2020-1}, and repeat the above discussion with  cuts $\mathcal{L}_\alpha$ issuing from $\lambda^*$,  so that all  the $c_{mk}$ are constant on 
   $\left\{ u\in \mathbb{D}(u^c) ~|~ u_{p_1+...+p_{\alpha-1}+1}\neq \lambda^*,\dots, u_{p_1+...+p_{\alpha}}\neq \lambda^*\right\}$. This proves constancy of the $c_{mk}$, $m$ associated with $\lambda_\alpha$, on the whole $\mathbb{D}(u^c)$. Then, we repeat this for all $\alpha=1,..,s$, proving constancy of the $c_{jk}$ for all $j=1,...,n$.   Hence,  Proposition  \ref{13ottobre2020-1} is  proved in the generic case. 
\vskip 0.2 cm 

{\bf General case}  of any $A(u)$.  
 If some of the diagonal entries $\lp_1,...,\lp_n$ of $A$ are integers, or some eigenvalues are integers, there exists a sufficiently small $\gamma_0>0$ such that, for any $0<\gamma<\gamma_0$, $A-\gamma I$ has diagonal non-integer  entries $\lp_1-\gamma,...,\lp_n-\gamma$  and no integer eigenvalues. Take such a  $\gamma_0$, and for any    $0<\gamma<\gamma_0$ consider 
 \be
\label{02bis}
(\Lambda-\lambda)\frac{d}{ d\lambda} (~{}_\gamma\Psi)=  \Bigl((A(u)-\gamma I) +I\Bigr)~{}_\gamma\Psi.
\ee 
namely 
\be
\label{03bis}
 \frac{d}{ d\lambda}({}_\gamma\Psi)=\sum_{k=1}^n \frac{B_k[\gamma](u) }{ \lambda-u_k}{}_\gamma\Psi,~~~~~B_k[\gamma](u):=-E_k\Bigl(A(u)+(1-\gamma)I\Bigr).
\ee

\ble The above system \eqref{03bis} is strongly  isomonodromic in $\mathbb{D}(u^0)$ contained in a $\tau$-cell, and  $\lambda$-component of the integrable  Pfaffian system 
\be
\label{24ottobre2020-1}
 d{}_\gamma\Psi=P_{[\gamma]}(\lambda,u) {}_\gamma\Psi, \quad\quad P_{[\gamma]}(\lambda,u)
 =
 \sum_{k=1}^n \frac{B_k[\gamma](u)}{\lambda -u_k} d(\lambda-u_k) +\sum_{j=1}^n [F_1(u),E_j]du_j.
 \ee
 where $F_1(u)$ is defined as in \eqref{26nov2018-3}, $(F_1)_{i j}=\frac{A_{ij}}{u_j-u_i}$, $i\neq j$,  and $[F_1(u),E_j]$  is   \eqref{22agosto2020-6}. 

\ele

\begin{proof}

We do a gauge transformation 
\be
\label{2settembre2020-2}
{}_\gamma Y(z):=z^{-\gamma} Y(z),\quad\quad \gamma\in\mathbb{C},
\ee
 which transforms \eqref{24nov2018-1} into 
\be 
\label{01bis}
\frac{d({}_\gamma Y)}{dz}=\left( \Lambda+\frac{A-\gamma I}{ z} \right)~ {}_\gamma Y
\ee
 For $u\in\mathbb{D}(u^0)$ contained in a $\tau$-cell, we  write the unique formal solution 
 \be
\label{2settembre2020-3}
 {}_\gamma Y_F(z,u)=z^{-\gamma}Y_F(z,u),
 \ee
 where $Y_F(z,u)$  is \eqref{15agosto2020-9}, so that  
 $$ 
 {}_\gamma Y_F(z,u)=F(z,u) z^{B-\gamma I} e^{\Lambda z},\quad\quad B-\gamma I=\hbox{\rm diag}(A-\gamma)=
 {\rm diag}(\lp_1-\gamma,~...~,\lp_n-\gamma).
 $$
The crucial point is that  $F(z,u)$ is the same as \eqref{15agosto2020-8}, so all the $F_k(u)$ are independent of $\gamma$. The fundamental matrix solutions $${}_\gamma Y_\nu(z,u):=z^{-\gamma} Y_\nu(z,u),
$$ 
are uniquely defined by their asymptotics $ {}_\gamma Y_F(z,u)$ in $\mathcal{S}_\nu(\mathbb{D}(u^0))$. Their Stokes matrices do not depend on $\gamma$  because 
 $$ 
 {}_\gamma Y_{\nu+(k+1)\mu}(z,u)= {}_\gamma Y_{\nu+k\mu}(z,u) \mathbb{S}_{\nu+k\mu}\quad
 \Longleftrightarrow 
 \quad 
 Y_{\nu+(k+1)\mu}(z,u)=  Y_{\nu+k\mu}(z,u) \mathbb{S}_{\nu+k\mu}.
$$
The system \eqref{01bis}  is thus strongly isomonodromic. 
By Proposition \ref{11agosto2020-5} we conclude. 
\end{proof}

\bcr 
Let the assumptions of Theorem \ref{30agosto2020-5} hold. Then Theorem \ref{30agosto2020-5}    holds also for  \eqref{24ottobre2020-1}.

\ecr

By Theorem \ref{30agosto2020-5} applied to \eqref{24ottobre2020-1}, we receive independent vector solutions ${}_\gamma\vec{\Psi}_k(\lambda,u)\equiv {}_\gamma\vec{\Psi}^{(sing)}_k(\lambda,u)$, $k=1,...,n$, which form a fundamental matrix 
$$ 
{}_\gamma\Psi(\lambda,u):=[{}_\gamma\vec{\Psi}_1(\lambda,u)~|~\cdots~|{}_\gamma\vec{\Psi}_n(\lambda,u)].
$$ 
For system  \eqref{24ottobre2020-1} the results  already proved  in the generic case  hold.  
Therefore, the  connection coefficients   $c_{jk}^{(\nu)}[\gamma]$ defined  by 
\be 
 \label{15agosto2020-7}
{}_\gamma \vec{\Psi}_k(\lambda,u~|\nu)={}_\gamma\vec{\Psi}_j(\lambda,u~|\nu) ~c_{jk}^{(\nu)}[\gamma]~+\hbox{\rm reg}(\lambda-u_j),\quad\quad \lambda\in\mathcal{P}_\eta,
 \ee
 are constant on $\mathbb{D}(u^c)$. 
They depend on $\gamma$, but not on $u\in\mathbb{D}(u^c)$.

\bre
{\rm It is explained   in section 8  of \cite{guz2016}   what is the relation between   $\vec{\Psi}^{(sing)}_k$ and $ {}_\gamma \vec{\Psi}_k$,  by means of their primitives, and that in general both $\lim_{\gamma\to 0}  {}_\gamma \vec{\Psi}_k$ and   $\lim_{\gamma\to 0}c_{jk}^{(\nu)}[\gamma] $ are divergent. 
}
\ere

 Now, we invoke Proposition 10 of \cite{guz2016}, which holds with no assumptions on eigenvalues and diagonal entries of $A(u)$.\footnote{The proof  in \cite{guz2016} is  laborious, because it is necessary to take into account all possible values of the diagonal entries  $\lp_k$ of $A$, including integer values. In \cite{BJL4} the proof is given only for non-integer values.} This result, adapted to our case, reads as follows. 

\bpr 
\label{20agosto2020-14}
Let $u$ be fixed in a $\tau$-cell. 
Let $\gamma_0>0$ be small enough such that for any $0<\gamma<\gamma_0$ the matrix $A-\gamma I$ has no integer eigenvalues, and its   diagonal part  has no integer entries.\footnote{Recall that eigenvalues and diagonal entries do not depend on $u$, in the isomonodromic case.}
 Let $c_{jk}^{(\nu)}$ be the  connection coefficients  of the Fuchsian system \eqref{03} at the fixed $u$, as in  Definition \ref{18ottobre2020-5}. 
 Let   $c_{jk}^{(\nu)}[\gamma]$  be the connection coefficients  in  \eqref{15agosto2020-7}. Let 
$$
\alpha_k:=\left\{\begin{array}{cc}
e^{-2\pi i \lp_k}-1,& \lp_k\not\in\mathbb{Z}
\\
\noalign{\medskip}
2\pi i, & \lp_k\in\mathbb{Z}
\end{array}
\right.;
~~~~~~~~\alpha_k[\gamma]:=e^{-2\pi i (\lp_k-\gamma)}-1
$$ 
Then, the following equalities hold
\be
\label{31agosto2020-7}
\alpha_k c_{jk}^{(\nu)}= 
e^{-2\pi i \gamma}\alpha_k[\gamma]~c_{jk}^{(\nu)}[\gamma],\quad  \hbox{ if }k\succ j;\quad\quad 
\alpha_k c_{jk}^{(\nu)}= \alpha_k[\gamma]~c_{jk}^{(\nu)}[\gamma] , \quad  \hbox{ if } k\prec j;
\ee
where the ordering relation $j \prec k$  means, for the fixed $u$,  that $\Re(z(u_j-u_k))<0 $ for $\arg z=\tau =3\pi/2-\eta$ satisfying \eqref{23agosto2020-16}. 
\epr

 We use Proposition \ref{20agosto2020-14} to conclude the proof of Proposition \ref{13ottobre2020-1} in the general case.  Indeed, the proposition  is already proved in the generic case, so it holds  for the $c_{jk}^{(\nu)}[\gamma]$. Therefore, they  are  constant on the whole $\mathbb{D}(u^c)$.  Equalities \eqref{31agosto2020-7} hold at any fixed $u$ in $\tau$-cell, so that each  $c_{jk}^{(\nu)}$ is constant on a $\tau$-cell, and  such consta nt is the same in each $\tau$-cell.  With a slight variation of $\eta$ in $(\eta_{\nu+1},\eta_\nu)$, equalities \eqref{31agosto2020-7}  hold also at the crossing locus $X(\tau)$. They analytically extend at $\Delta$.
 
 \end{proof}

    
    \section{Laplace Transform in $\mathbb{D}(u^c)$, Theorem \ref{20agosto2020-5}}
 \label{20agosto2020-10}   
 By means of the Laplace transform with deformation parameters, we prove     points {\it (I1),(I2), (I3), (II1), (II2)} and {\it (II5)}   of  Theorem \ref{18agosto2020-8}. 
    Stokes matrices will be expressed in terms of the isomonodromic connection coefficients satisfying Proposition \ref{13ottobre2020-1}. The result is in Theorem \ref{20agosto2020-5} below, which is the last step of our construction.

  Let $\tau$ be the chosen direction in the $z$-plane admissible at $u^c$, and $\eta=3\pi/2-\tau$ in the $\lambda$-plane.  The Stokes rays of $\Lambda(u^c)$ will be labelled as in  \eqref{15novembre2020-1}, so that  \eqref{23agosto2020-16}  holds for  a certain  $\nu\in\mathbb{Z}$. We define the sectors
\be
\label{20agosto2020-12}
\mathcal{S}_\nu=\{z\in\mathcal{R}(\mathbb{C}\backslash\{0\})~\hbox{ such that }~\tau_\nu-\pi<\arg z <\tau_{\nu+1}\}.
\ee

If $u$ only varies in $\mathbb{D}(u^0)$ contained in a ${\tau}$-cell, then none of the Stokes rays  associated with $\Lambda(u)$ crosses $\arg z={\tau}$ mod $\pi$.
 If $u$ varies in $\mathbb{D}(u^c)$,  some Stokes rays associated with $\Lambda(u)$ necessarily cross $\arg z={\tau}$ mod $\pi$ (see Section \ref{20agosto2020-15}). 
 Consider the subset 
 of the set of Stokes rays  satisfying  $\Re(z(u_j-u_k))=0$, $z\in\mathcal{R}$, 
 associated with pairs $(u_j,u_k)$   such that  $u_j\in\mathbb{D}_\alpha$ and $u_k\in \mathbb{D}_\beta$,   $\alpha\neq \beta$, namely $u_j^c\neq u_k^c$. Following \cite{CDG}, we denote this subset by $\mathfrak{R}(u)$. 
 If $u$ varies in $\mathbb{D}(u^c)$ and $\epsilon_0$ satisfies \eqref{30agosto2020-6},   the rays in $\mathfrak{R}(u)$ continuously rotate, but    never cross  the admissible rays $\arg z={\tau}+h\pi$, where 
 \be
 \label{15novembre2020-2}
 \tau_{\nu+h\mu}<{\tau}+h\pi<\tau_{\nu+h\mu+1}, \quad \quad h\in \mathbb{Z},
 \ee
 The above allows to define $\widehat{\mathcal{S}}_{\nu+h\mu}(u)$ to be the unique sector containing $S\bigl({\tau}+(h-1)\pi,{\tau}+h\pi\bigr)$  and extending up to the nearest Stokes rays in $\mathfrak{R}(u)$. Then, let  
\be
\label{16ottobre2020-1}
\widehat{\mathcal{S}}_{\nu+h\mu}:=\bigcap_{u\in \mathbb{D}(u^c)} \widehat{\mathcal{S}}_{\nu+h\mu}(u). 
\ee
It has angular amplitude greater than $\pi$. The reason for the labeling is that 
$
\widehat{\mathcal{S}}_{\nu+h\mu}(u^c)=\mathcal{S}_{\nu+h\mu}$  in \eqref{20agosto2020-12}. 

Suppose that $u$ is fixed in a $\tau$-cell.  Let   $$
Y_{\nu+h\mu}(z,u):=\Bigl[\vec{Y}_1(z,u~|\nu+h\mu)~\Bigr| ~\dots ~\Bigr| ~ \vec{Y}_n(z,u~|\nu+h\mu)\Bigr],
$$
 be  defined by 
\begin{align}
\label{15agosto2020-1}
&
\vec{Y}_k(z,u~|\nu+h\mu):= \dfrac{1}{2\pi i } \int_{\gamma_k(\eta-h\pi)} e^{z\lambda} \vec{\Psi}_k^{(sing)}(\lambda,u~|\nu+h\mu) d\lambda,& \hbox{ for $\lp_k\not\in\mathbb{Z}_{-}$,}
\\
\noalign{\medskip}
\label{15agosto2020-1-bis}
&\vec{Y}_k(z,u~|\nu+h\mu):= \int_{L_k(\eta-h\pi)} e^{z\lambda} \vec{\Psi}_k(\lambda,u~|\nu+h\mu) d\lambda,  & \hbox{  for $\lp_k\in\mathbb{Z}_{-}$}.
\end{align}
In the $\lambda$-plane, the  admissible directions  $\eta-h\pi$  correspond to ${\tau}+h\pi$, with 
 \be
 \label{15novembre2020-3}
   \eta_{\nu+h\mu+1}<\eta-h\pi< \eta_{\nu+h\mu}.
   \ee
Here, $\vec{\Psi}_k(\lambda,u~|\nu+h\mu)$, $\vec{\Psi}_k^{(sing)}(\lambda,u~|\nu+h\mu)$ are the vector solutions of Theorem \ref{30agosto2020-5} for $\lambda\in \mathcal{P}_{\eta-h\pi}(u)$, with $u$  {\it fixed} in a $\tau$-cell.    
    $L_k(\eta-h\pi)$ is the cut in direction $\eta-h\pi$,  issuing from $u_k$ and  {\it oriented} from $u_k$ to $\infty$, and  $\gamma_k(\eta-h\pi)$ is the path coming from $\infty$ along the left side of $L_k(\eta-h\pi)$, encircling $u_k$ with a small loop excluding all the other poles, and going back to $\infty$ along the right side of $L_k(\eta-h\pi)$.  Here ``right'' and ``left'' refer to the orientation of $L_k(\eta-h\pi)$. The label $\nu+h\mu$ keeps track of  \eqref{23agosto2020-16} and \eqref{15novembre2020-2}-\eqref{15novembre2020-3}.

        \bth
     \label{20agosto2020-5}
     
      Let the assumptions of Theorem \ref{30agosto2020-5} hold.

  \begin{itemize}
  \item[{\bf 1)}]   The  matrices $Y_{\nu+h\mu}(z,u)$, obtained by Laplace transform \eqref{15agosto2020-1}-\eqref{15agosto2020-1-bis} at a fixed $u\in \mathbb{D}(u^0)$  contained in a $\tau$-cell, define  holomorphic matrix valued  functions of $(\lambda,u)\in \mathcal{R}(\mathbb{C}\backslash\{0\})\times \mathbb{D}(u^c)$, which are fundamental matrix solutions of  \eqref{24nov2018-1}.
   
  \item[{\bf 2)}]  They have structure 
$$ 
Y_{\nu+h\mu}(z,u)=\widehat{Y}_{\nu+h\mu}(z,u) z^{B} e^{z\Lambda},\quad\quad
  B=\hbox{\rm diag}(\lp_1,...,\lp_n),
  $$
  with asymptotic behaviour, uniform in $u\in\mathbb{D}(u^c)$,  
  $$
  \widehat{Y}_{\nu+h\mu}(z,u) \sim F(z,u)=I+\sum_{l=1}^\infty \frac{F_l(u)}{z^l}, \quad \quad z\to\infty \hbox{ in } \widehat{S}_{\nu+h\mu},
  $$
  given by the formal solution 
  $ 
  Y_F(z,u)=F(z,u)z^{B} e^{z\Lambda}
  $. 
The coefficients $ F_l(u)$  are holomorphic in $\mathbb{D}(u^c)$. The explicit expression of their columns is  \eqref{15novembre2020-5}, \eqref{17ottobre2020-3},  \eqref{17ottobre2020-4} (or \eqref{15novembre2020-7})  and  \eqref{15novembre2020-6}.

  \item[{\bf 3)}] Stokes matrices defined by
  \be
  \label{26novembre2020-3}
  Y_{\nu+(h+1)\mu}(z,u)=Y_{\nu+h\mu}(z,u) \mathbb{S}_{\nu+h\mu},
  \quad
  \quad
  z\in \widehat{\mathcal{S}}_{\nu+h\mu}\cap \widehat{\mathcal{S}}_{\nu+(h+1)\mu},
\ee
 are constant in the whole $\mathbb{D}(u^c)$ and satisfy 
\be
\label{16ottobre2020-6}
(\mathbb{S}_{\nu+h\mu})_{ab}=(\mathbb{S}_{\nu+h\mu})_{ba}=0 \quad\hbox{ for $a\neq b$ such that $u_a^c=u_b^c$}.
\ee

\item[{\bf 4)}] The following representation  in terms of the constant connection coefficients $c_{jk}^{(\nu)}$ of Proposition \ref{13ottobre2020-1} holds on $\mathbb{D}(u^c)$:
\be
\label{16ottobre2020-5}
\begin{array}{c}
(\mathbb{S}_\nu)_{jk}=\quad\quad\quad\quad\quad\quad\quad
\\
\\
\left\{
\begin{array}{cc}
e^{2\pi i \lp_k}\alpha_k ~c_{jk}^{(\nu)},&\hbox{ $j\prec k$,  $u_j^c\neq u_k^c$},
\\
\noalign{\medskip}
             1   &  j =k,
\\
\noalign{\medskip}
             0   &  \hbox{ $j\succ k$, $u_j^c\neq u_k^c$},
             \\
\noalign{\medskip}
             0   &  \hbox{  $j\neq  k$ ,  $u_j^c= u_k^c$},
\end{array}
\right.;
\end{array}
\quad
\quad
\begin{array}{c}
(\mathbb{S}_{\nu+\mu}^{-1})_{jk} 
= \quad\quad
\\
\\
\left\{
\begin{array}{cc}
0   &  \hbox{ $j\neq  k$,  $u_j^c= u_k^c$},
\\
\noalign{\medskip}
             0   &  \hbox{  $j\prec k$, $u_j^c\neq u_k^c$},
\\
\noalign{\medskip}
             1   & j =k,
\\
\noalign{\medskip}
         -e^{2\pi i (\lp_k-\lp_j)}\alpha_k~c_{jk}^{(\nu)} & 
                                         \hbox{  $j\succ k$,  $u_j^c\neq u_k^c$},
\end{array}
\right.
\end{array} 
\ee
\end{itemize}
where the relation  $j\prec k$ is  defined for $j\neq k$ such that $u_j^c\neq u_k^c$ and  means  that $\Re(z(u_j^c-u_k^c))<0$ when $\arg z=\tau$. 
\eth 

\bre{\rm 
The above \eqref{16ottobre2020-5}  generalises Theorem \ref{20agosto2020-1} in presence of  isomonodromic deformation parameters, including coalescences. Notice that the ordering relation $\prec$ here is referred to $u^c$, while in Theorem \ref{20agosto2020-1} it refers to $u^0$.  
}
\ere

\begin{proof} We use  the labelling \eqref{2settembre2020-8}-\eqref{2settembre2020-9} for $u^c$. 
{~}
\vskip 0.2 cm
\noindent
{\bf a)} \underline{Case $\lp_k\not\in \mathbb{Z}$}.

$\bullet$ {\bf Construction of $ \vec{Y}_k(z,u~|\nu)$}.  We have 
$
\vPsi_k^{(sing)}(\lambda,u|~\nu)=\vec{\Psi}_k(\lambda,u|~\nu)$ and \eqref{15agosto2020-1} is    \be
   \label{31agosto2020-3}
    \vec{Y}_k(z,u~|\nu):= \dfrac{1}{2\pi i } \int_{\gamma_k(\eta)} e^{z\lambda} \vec{\Psi}_k(\lambda,u~|\nu) d\lambda
   \ee
   Since $ \vec{\Psi}_k(\lambda,u~|\nu)$ grows at infinity no faster than some power of $\lambda$, the integral converges in a sector of amplitude at most $\pi$.  Now, $\vec{\Psi}_k(\lambda,u~|\nu) $ satisfies 
  Theorem \ref{30agosto2020-5}, hence  
  if  $u$ varies in $\mathbb{D}(u^c)$ the following facts hold.
  \begin{itemize}
  \item[1.]  $\vec{\Psi}_k(\lambda,u~|\nu)$ is branched at $\lambda=u_k$ and possibly at other poles $u_l$ such that  $u_l^c\neq u_k^c$. 
  
  \item[2.] 
  $\vec{\Psi}_k(\lambda,u~|\nu)$ is holomorphic at all $\lambda= u_j$ such that $u_j^c = u_k^c$, $j\neq k$.
  
  \end{itemize}
  It follows from 1. and 2. that  the path of integration can be modified:  for $\alpha$ such that $u_k^c=\lambda_\alpha$, we have
  \be
  \label{14agosto2020-1}
   \vec{Y}_k(z,u~|\nu)= \dfrac{1}{2\pi i } \int_{\Gamma_\alpha(\eta)} e^{z\lambda} \vec{\Psi}_k(\lambda,u~|\nu) d\lambda
   ,
   \ee
  where  $\Gamma_\alpha(\eta)$ is the path which  comes from $\infty$ in direction $\eta-\pi$, encircles $\lambda_\alpha$ along $\partial \mathbb{D}_\alpha$ anti-clockwise and goes to $\infty$ in direction $\eta$. This path encloses all the $u_j$ such that $u_j^c=\lambda_\alpha$, end excludes the others. See figure \ref{2settembre2020-11}. 
We conclude that  $u$ can vary in $\mathbb{D}(u^c)$ and the integral \eqref{14agosto2020-1} converges for $z$ in the sector
  $$ \mathcal{S}(\eta):=\Bigl\{
  z\in\mathcal{R}(\mathbb{C}\backslash\{0\})~\hbox{ such that } \frac{\pi}{2}-\eta <\arg z <\frac{3\pi}{2}-\eta
  \Bigr\},
  $$
defining $\vec{Y}_k(z,u~|\nu)$ as a holomorphic function of $(z,u)\in  \mathcal{S}(\eta)\times \mathbb{D}(u^c)$. 
   Now,   if  $u$ varies in $\mathbb{D}(u^c)$ and $\epsilon_0$ satisfies \eqref{30agosto2020-6}  none of the vectors $
    u_i-u_j $ such that $ u_i^c=\lambda_\alpha$ and $ u_j^c=\lambda_\beta$, $ 1\leq \alpha \neq \beta \leq s$, 
cross a direction $\eta$ mod $\pi$, for every $\eta_{\nu+1}<\eta<\eta_\nu$. Due to 1. and 2. above, a vector function $\vec{\Psi}_k(\lambda,u~|\nu)$ is well defined in $\mathcal{P}_\eta$ and $\mathcal{P}_{\tilde{\eta}}$ for any  $\eta_{\nu+1}<\eta<\tilde{\eta}< \eta_\nu$, and so on  $\mathcal{P}_\eta \cup \mathcal{P}_{\tilde{\eta}}$. Therefore, the integral in \eqref{14agosto2020-1} satisfies 
$$ 
 \dfrac{1}{2\pi i } \int_{\Gamma_\alpha(\eta)}e^{z\lambda} \vec{\Psi}_k(\lambda,u~|\nu) d\lambda= \dfrac{1}{2\pi i } \int_{\Gamma_\alpha(\tilde{\eta})}e^{z\lambda} \vec{\Psi}_k(\lambda,u~|\nu) d\lambda,
 \quad \quad 
 z\in \mathcal{S}(\eta)\cap \mathcal{S}(\tilde{\eta}),
 $$ 
 namely one is the analytic continuation of the other, so defining the function $\vec{Y}_k(z,u~|\nu)$ as analytic on $\widehat{\mathcal{S}}_\nu\times \mathbb{D}(u^c)$, where  $$
 \widehat{\mathcal{S}}_\nu:=\bigcup_{\eta_{\nu+1}<\eta<\eta_\nu}\mathcal{S}(\eta)
 $$ 
coincides with \eqref{16ottobre2020-1} (with $h=0$). 
  \begin{figure}
\centerline{\includegraphics[width=0.8\textwidth]{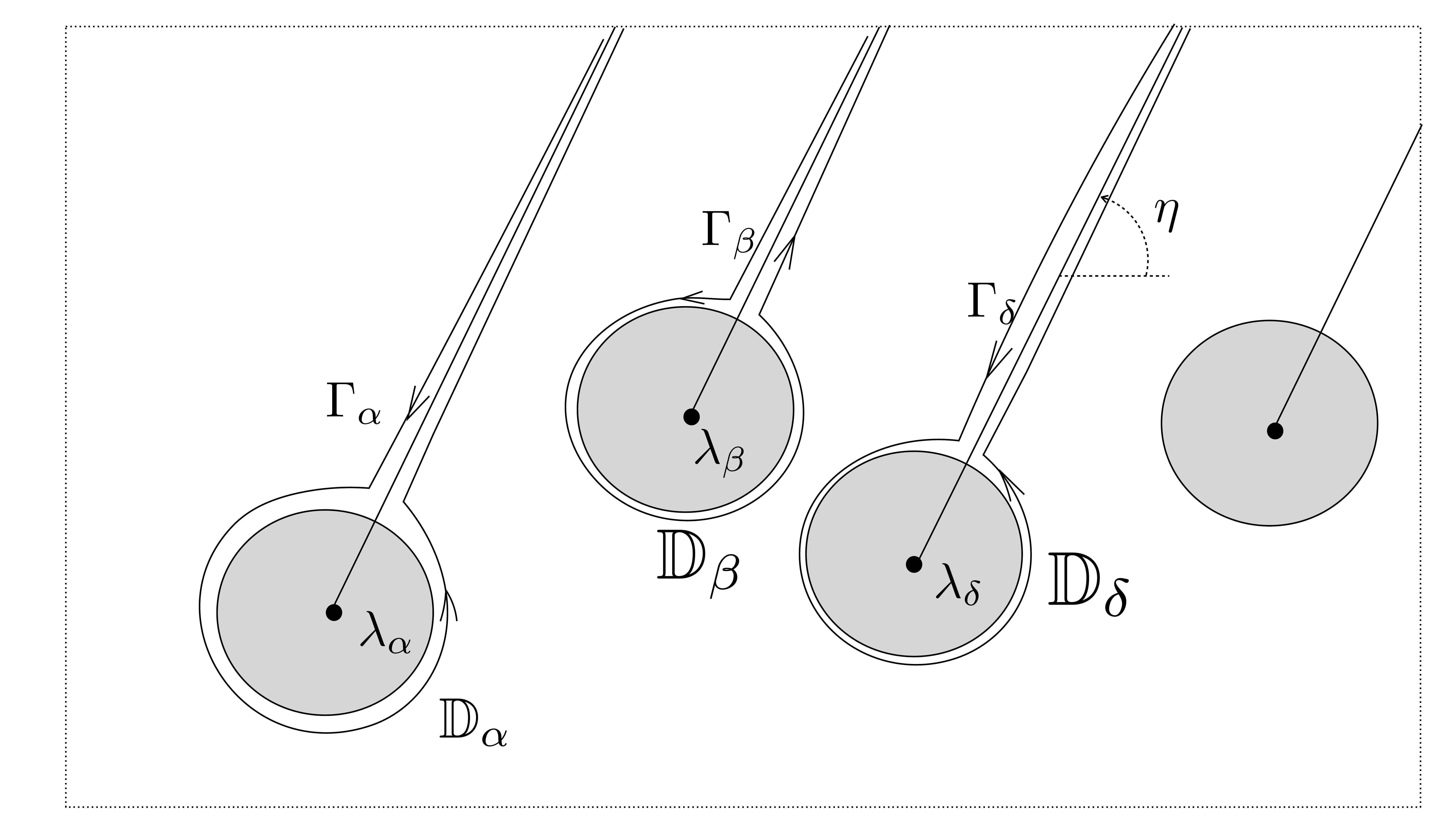}}
\caption{The paths of integration $\Gamma_\alpha$, $\Gamma_\beta$, etc $\alpha,\beta, ...\in\{1,...,s\}$.}
\label{2settembre2020-11}
\end{figure}
Finally, notice that $e^{\lambda z}(\lambda-\Lambda) \vec{\Psi}_k(\lambda,u~|\nu)\Bigl|_{\Gamma(\alpha)}=0$, due to the exponential factor.  By \eqref{31agosto2020-2}, the vector solutions $\vec{Y}_k(z,u~|\nu)$  satisfy  system  \eqref{24nov2018-1}.

\vskip 0.2 cm 
  $\bullet$ {\bf Asymptotic behaviour}. From \eqref{18ottobre2020-1}-\eqref{31agosto2020-1}, we write  \eqref{14agosto2020-1} as 
$$
 \vec{Y}_k(z,u~|\nu)= \dfrac{1}{2\pi i } \int_{\Gamma_\alpha(\eta)} e^{z\lambda} \Bigl(\Gamma(\lp_j+1)\vec{e}_j+\sum_{l\geq 1} \vec{b}_l^{(k)}(u)(\lambda-u_k)^l\Bigr)(\lambda-u_k)^{-\lp_k-1},
 d\lambda.
 $$
 with holomorphic $\vec{b}_l^{(k)}(u)$ on $\mathbb{D}(u^c)$. 
 We split the series as $\sum_{l\geq 1}=\sum_{l=1}^{\mathcal{N}}+\sum_{l\geq \mathcal{N}+1}$, and recall the standard formula (see \cite{Doe})
  $$ 
  \int_{\Gamma_\alpha(\eta)}(\lambda-\lambda_k)^a e^{z\lambda} d\lambda= \int_{\gamma_k(\eta)}(\lambda-\lambda_k)^a e^{z\lambda} d\lambda=\frac{z^{-a-1} e^{\lambda_k z}}{\Gamma(-a)}
  $$
 so that 
 $$ 
 \vec{Y}_k(z,u~|\nu)=\left( \vec{e}_k+\sum_{l=1}^{\mathcal{N}} \frac{\vec{b}_l^{(k)}(u)}{\Gamma(\lp_k+1-l)}z^{-l}+R_{\mathcal{N}}(z)\right)~z^{\lp_k}e^{\lambda_k z},
  $$
  with remainder $$ 
  R_{\mathcal{N}}(z)
=
\oint_{\Gamma_0(\eta)}\sum_{l\geq \mathcal{N}} \frac{\vec{b}_l^{(k)}(u) }{z^l} ~e^x  x^{l-\lp_k-1} ~dx~=O(z^{-\mathcal{N}+1}).
  $$
  The integral is along a path $\Gamma_0(\eta)$, coming from $\infty$ along the left part of the half line oriented from $0$ to $\infty$ in direction $\eta+\arg z$,  going around $0$, and back to $\infty $ along the right part. The estimate $O(z^{-\mathcal{N}+1})$ is standard. We conclude that 
$$
 \vec{Y}_k(z,u~|\nu) \left(z^{\lp_k}e^{u_k z}\right)^{-1}
  \sim  \vec{e}_k+\sum_{l=1}^{\infty} \frac{\vec{b}_l^{(k)}(u)}{\Gamma(\lp_k+1-l)}z^{-l}\equiv \vec{e}_k+\sum_{l=1}^{\infty} \vec{f}_l^{~(k)}(u) z^{-l},\quad \quad z\to \infty \hbox{ in }  \widehat{\mathcal{S}}_\nu
  $$
  with 
  \be
  \label{15novembre2020-5}
   \vec{f}_l^{~(k)}(u):= \frac{\vec{b}_l^{(k)}(u)}{\Gamma(\lp_k+1-l)}.
  \ee

\vskip 0.2 cm 
\noindent
{\bf b)} \underline{Case $\lp_k\in\mathbb{N}=\{0,1,2,...\}$}.

$\bullet$ {\bf Construction of $ \vec{Y}_k(z,u~|\nu)$}.  Definition \eqref{15agosto2020-1} is
\begin{align*}
\vec{Y}_k(z,u~|\nu)&:= \frac{1}{2\pi i} \int_{\gamma_k(\eta)} e^{z\lambda} \vec{\Psi}_k^{(sing)}(\lambda,u~|\nu) d\lambda
\\
&
\underset{\eqref{17ottobre2020-1}}=
 \frac{1}{2\pi i} \int_{\gamma_k(\eta)} e^{z\lambda} \left(
 \frac{\vec{\psi}_k(\lambda,u~|\nu)}{(\lambda-u_k)^{\lp_k+1}}+ \vec{\Psi}_k(\lambda,u~|\nu)\ln(\lambda-u_k) 
 \right) d\lambda.
\end{align*}
The same facts 1. and 2. of the previous case apply to $ \vec{\Psi}_k(\lambda,u~|\nu)$ and $\vec{\psi}_k(\lambda,u~|\nu)$ and allow to rewrite 
\begin{align*}
\vec{Y}_k(z,u~|\nu)
&=
 \frac{1}{2\pi i} \int_{\Gamma_\alpha(\eta)} e^{z\lambda} \left(
 \frac{\vec{\psi}_k(\lambda,u~|\nu)}{(\lambda-u_k)^{\lp_k+1}}+ \vec{\Psi}_k(\lambda,u~|\nu)\ln(\lambda-u_k) 
 \right) d\lambda
 \\
 &
 =
 \frac{1}{2\pi i} \int_{\Gamma_\alpha(\eta)} e^{z\lambda} 
 \vec{\Psi}_k^{(sing)}(\lambda,u~|\nu) d\lambda.
 \end{align*}
 We conclude that  $\vec{Y}_k(z,u~|\nu)$ is analytic on $\widehat{\mathcal{S}}_\nu\times \mathbb{D}(u^c)$. Moreover,   $e^{\lambda z}(\lambda-\Lambda) \vec{\Psi}_k^{(sing)}(\lambda,u~|\nu)\Bigl|_{\Gamma(\alpha)}=0$, due to the exponential factor.  By \eqref{31agosto2020-2}, the vector solution $\vec{Y}_k(z,u~|\nu)$  satisfies the system  \eqref{24nov2018-1}.
  
  \vskip 0.2 cm 
  $\bullet$ {\bf Asymptotic behaviour}. By \eqref{31agosto2020-1-bis} and \eqref{17ottobre2020-2}, and the fact that $\vec{\psi}_k$ has no singularities at $u_j\in \mathbb{D}_\alpha$, $j\neq k$,  so that the terms $ \sum_{l\geq 1+\lp_k} \vec{b}_l^{~(k)}(u) (\lambda-u_k)^l$ in $\vec{\psi}_k(\lambda,u~|\nu)$ do not contribute to the integration,  we can write
  $$ 
\vec{Y}_k(z,u~|\nu)=  \frac{1}{2\pi i} \int_{\Gamma_\alpha(\eta)} \left(
  \frac{ \lp_k!\vec{e}_k+
   \sum_{l=1}^{\lp_k} \vec{b}_l^{~(k)}(u) (\lambda-u_k)^l}{(\lambda-u_k)^{\lp_k+1}}+
  \sum_{l=0}^\infty \vec{d}_l^{~(k)}(u)(\lambda-u_k)^l~\ln(\lambda-u_k)
  \right)e^{z\lambda}~d\lambda .
  $$
 By Cauchy formula 
\begin{align*}
  \frac{1}{2\pi i} \int_{\Gamma_\alpha(\eta)} \left(
  \frac{ \lp_k!\vec{e}_k+
   \sum_{l=1}^{\lp_k} \vec{b}_l^{~(k)}(u) (\lambda-u_k)^l}{(\lambda-u_k)^{\lp_k+1}}\right)e^{z\lambda}~d\lambda
  & =
   \frac{1}{\lp_k!} \frac{d^{\lp_k}}{d\lambda^{\lp_k}}\left.\left[ \left(
\lp_k!\vec{e}_k+
   \sum_{l=1}^{\lp_k} \vec{b}_l^{~(k)}(u) (\lambda-u_k)^l\right)
   e^{z\lambda}
   \right] \right|_{\lambda=u_k}
 \\
 &
 =
 z^{\lp_k} e^{u_k z}\left(\vec{e}_k+ \sum_{l=1}^{\lp_k}  \vec{f}_l^{~(k)}(u) \frac{1}{z^l} \right),
\end{align*}
 where 
 \be
 \label{17ottobre2020-3}
   \vec{f}_l^{~(k)}(u):=\frac{\vec{b}_l^{(k)}(u)}{(\lp_k-l)!},\quad l=1,...,\lp_k.
 \ee
 In order to evaluate the terms with logarithm, we observe that for any function $g(\lambda)$ holomorphic   along $L_k(\eta)$, including $\lambda=u_k$,  we have 
 $$
 \int_{\gamma_k(\eta)}g(\lambda)\ln(\lambda-u_k)d\lambda= \int_{L_k(\eta)^{-}} g(\lambda)\ln(\lambda-u_k)_{-}d\lambda-
  \int_{L_k(\eta)^{+}} g(\lambda)\ln(\lambda-u_k)_{+}d\lambda,
  $$
  where $L_k(\eta)^{+}$ and $L_k(\eta)^{-}$ respectively are the left and right  parts of $L_k(\eta)$, oriented  from $0$ to $\infty$. Since 
  $\ln(\lambda-u_k)_{+}=\ln(\lambda-u_k)_{+}-2\pi i $, we conclude that 
\be
\label{18ottobre2020-4}
   \int_{\gamma_k(\eta)}g(\lambda)\ln(\lambda-u_k)d\lambda=  2\pi i \int_{L_k(\eta)} g(\lambda)d\lambda.
  \ee
    Keeping into account that the integral along $\Gamma_\alpha$ can be interchanged with that along $\gamma_k$, it follows that
\begin{align*}
  \frac{1}{2\pi i} \int_{\Gamma_\alpha(\eta)} 
 \vec{\Psi}_k(\lambda,u~|\nu)\ln(\lambda-u_k)
  e^{z\lambda}~d\lambda
  &
  = \int_{L_k(\eta)} 
 \vec{\Psi}_k(\lambda,u~|\nu)
  e^{z\lambda}~d\lambda
\\
  &
  =\int_{L_k(\eta)}
  \sum_{l=0}^\infty \vec{d}_l^{~(k)}(u)(\lambda-u_k)^l~
  e^{z\lambda}~d\lambda.
\end{align*}
We conclude, by the standard evaluation of the remainder analogous to $R_{\mathcal{N}}(z)$ considered before, and the variation of $\eta$ in  the range $(\eta_{\nu+1},\eta_\nu)$, that\footnote{Notice that, by abuse of notation, if $f(\lambda)e^{-u_k\lambda}\sim \sum_{0}^\infty c_l z^{-l}$ we write $f(\lambda)\sim e^{u_k\lambda} \sum_{0}^\infty c_lz^{-l}$.}
\begin{align*}
 \int_{L_k(\eta)} 
 \vec{\Psi}_k(\lambda,u~|\nu)
  e^{z\lambda}~d\lambda
  &
   \sim  e^{u_k z}\left( \sum_{l=0}^\infty (-1)^{l+1} l! ~  \vec{d}_l^{(k)}(u) ~z^{-l-1}\right),\quad \quad z\to\infty \hbox{ in } \widehat{\mathcal{S}}_\nu.
\\
&
 =z^{\lp_k}e^{u_k z}\left(\sum_{l=\lp_k+1}^\infty  \vec{f}_l^{~(k)}(u)~z^{-l}\right),
 \end{align*}
 where 
 \be
 \label{17ottobre2020-4} 
  \vec{f}_l^{~(k)}(u):=(-1)^{l-\lp_k}(l-\lp_k-1)!~\vec{d}_{l-\lp_k-1}^{~(k)}(u) ,\quad\quad  l\geq \lp_k+1.
 \ee
 
 In conclusion, we have the expansion 
 $$ 
 \vec{Y}_k(z,u~|\nu)\sim z^{\lp_k}e^{u_k z}\left(\vec{e}_k+\sum_{l=1}^\infty  \vec{f}_l^{~(k)}(u)~z^{-l}\right),   \quad \quad z\to\infty \hbox{ in } \widehat{\mathcal{S}}_\nu,
 $$ 
 with coefficients $ \vec{f}_l^{~(k)}(u)$  holomorphic in $\mathbb{D}(u^c)$ defined in \eqref{17ottobre2020-3}-\eqref{17ottobre2020-4}. Notice that, in exceptional cases, $\vec{\Psi}_k$ may be identically zero, so that  
 \be
 \label{15novembre2020-7}
 \hbox{$\vec{f}_l^{~(k)}=0$ for $ l\geq \lp_k+1$}.
 \ee

 \vskip 0.2 cm 
 \noindent
{\bf c)}  \underline{Case $\lp_k\in\mathbb{Z}_{-}=\{-1,-2,...\}$}
 \vskip 0.2 cm $\bullet$ {\bf Construction of $ \vec{Y}_k(z,u~|\nu)$}. Definition \eqref{15agosto2020-1-bis} is
 $$
 \vec{Y}_k(z,u~|\nu):= \int_{L_k(\eta)} e^{\lambda z} \vec{\Psi}_k(\lambda,u~|\nu) d\lambda\equiv  \int_{\mathcal{L}_\alpha(\eta)} e^{\lambda z} \vec{\Psi}_k(\lambda,u~|\nu) d\lambda
 .
 $$
 In the last equality, we have used the fact that $ \vec{\Psi}_k(\lambda,u~|\nu)$ is analytic in $\mathbb{D}_\alpha\times \mathbb{D}(u^c)$, where $\lambda_\alpha=u_k^c$.  
 
 We conclude analogously to previous cases that $ \vec{Y}_k(z,u~|\nu)$ is analytic in $\widehat{\mathcal{S}}_\nu\times \mathbb{D}(u^c)$. It is a solution of \eqref{24nov2018-1}, by \eqref{31agosto2020-2}, because $\vec{\Psi}_k(\lambda,u~|\nu)$ is analytic at $\lambda=u_k$ and behaves as in \eqref{18ottobre2020-1}-\eqref{31agosto2020-1}, so that 
 $$
 e^{\lambda z} (\lambda I - \Lambda(u)) \vec{\Psi}_k(\lambda,u~|\nu)\Bigr|_{\mathcal{L}_\alpha}=
  e^{\lambda z} (\lambda I - \Lambda(u)) \vec{\Psi}_k(\lambda,u~|\nu)\Bigr|_{L_k}=0-(u_kI-\Lambda(u)) \vec{\Psi}_k(\lambda,u_k~|\nu)=0.
  $$
   \vskip 0.2 cm 
  $\bullet$ {\bf Asymptotic behaviour}. We have, from  \eqref{18ottobre2020-1}-\eqref{31agosto2020-1}, 
  $$
  \vec{Y}_k(z,u~|\nu)= \int_{\mathcal{L}_\alpha(\eta)} e^{\lambda z} \left(
  \frac{(-1)^{\lp_k} ~\vec{e}_k}{(-\lp_k-1)!} (\lambda-u_k)^{-\lp_k-1}+\sum_{l\geq 1} \vec{b}_l^{(k)}(u)(\lambda-u_k)^{l-\lp_k-1}
  \right) d\lambda
  $$
  We integrate term by term in order to obtain the asymptotic expansion (the remainder for the truncated series is evaluate in standard way, as $R_{\mathcal{N}}(z)$ above). For the integration, we use 
$$
  \int_{L_k(\eta)} (\lambda-u_k)^m e^{\lambda z} d\lambda
  =\frac{e^{u_kz}}{z^{m+1}}\int_{+\infty e^{i\phi}}^0 x^m e^x dx  = \frac{e^{u_kz}}{z^{m+1}} m!~(-1)^{m+1}, \quad\quad \quad \frac{\pi}{2}<\phi<\frac{3\pi}{2}.
$$
We obtain, analogously to previous cases,  
$$ 
  \vec{Y}_k(z,u~|\nu)\sim z^{\lp_k}e^{u_k z}\left(\vec{e}_k + \sum_{l=1}^\infty \vec{f}_l^{~(k)}(u) z^{-l}\right),\quad \quad 
  z\to\infty \hbox{ in } \widehat{\mathcal{S}}_\nu,
 $$ 
 where the holomorphic in $\mathbb{D}(u^c)$ coefficients are 
\be
\label{15novembre2020-6}
\vec{f}_l^{~(k)}(u):= (-1)^{l-\lp_k} (l-\lp_k-1)!~\vec{b}_l^{(k)}(u).
\ee

\bre{\rm 
We cannot use  $\vec{\Psi}_k^{(sing)}(\lambda,u~|\nu)$ in \eqref{18ottobre2020-2}  to  define $ \vec{Y}_k(z,u~|\nu)$ if $u$ varies in the whole $\mathbb{D}(u^c)$. On the other hand, if $u$ is restricted to a $\tau$-cell, so that the eigenvalues $u_j$ are all distinct, by \eqref{18ottobre2020-4} we can  write 
$$  \vec{Y}_k(z,u~|\nu)= \int_{L_k(\eta)} e^{\lambda z} \vec{\Psi}_k(\lambda,u~|\nu) d\lambda\underset{\eqref{18ottobre2020-4}}=
\frac{1}{2\pi i } \int_{\gamma_k(u)} e^{\lambda z}\vec{\Psi}_k(\lambda,u~|\nu) \ln(\lambda-u_k)d\lambda.
$$
Then, we can use the local expansion \eqref{18ottobre2020-3} and the fact that $\int_{\gamma_k(u)} \hbox{\rm reg}(\lambda-u_k)d\lambda=0$, receiving
$$
 \vec{Y}_k(z,u~|\nu) =\frac{1}{2\pi i } \int_{\gamma_k(u)} e^{\lambda z}\vec{\Psi}_k^{(sing)}(\lambda,u~|\nu) d\lambda
$$
   }
   \ere

  \subsubsection*{Fundamental matrix solutions}

The vector solutions $\vec{Y}_k(z,u~|\nu) $ constructed above can be arranged as columns of the matrix 
$$ 
Y_\nu(z,u):=\Bigl[\vec{Y}_k(z,u~|\nu)~\Bigl|~\cdots~\Bigr|~\vec{Y}_n(z,u~|\nu)\Bigr], 
$$
which thus solves system \eqref{24nov2018-1}. From the general theory of differential systems,   it admits analytic continuation as analytic matrix valued function  on $\mathcal{R}(\mathbb{C}\backslash\{0\})\times \mathbb{D}(u^c)$. 
 Letting $B=\hbox{\rm diag} A=\hbox{\rm diag}(\lp_1,...,\lp_n)$, the asymptotic expansions obtained above are  summarized as
  $$ 
  Y_\nu(z,u~|\nu) ~z^{-B}e^{-\Lambda(u) z}\sim F(z,u)=I+\sum_{l=1}^\infty F_l(u) z^{-l},\quad \quad z\to \infty \hbox{ in }  \widehat{\mathcal{S}}_\nu
  ,
  $$
 $$
  F_l(u)=\Bigl[\vec{f}_l^{~(1)}(u)~|~\cdots~|~\vec{f}_l^{~(n)}(u)\Bigr].
  $$
    Therefore, the coefficients $F_l(u)$ of the formal  solution $Y_F(z,u)=F(z,u)z^{B}e^{\Lambda(u) z}$   are holomorphic in $\mathbb{D}(u^c)$.  Moreover, the leading term is the identity $I$, which implies that $
Y_\nu(z,u)$ is a fundamental matrix solution.

   \vskip 0.2 cm 
      
Consider now another direction $\eta$, satisfying  $\eta_{\nu+\mu+1}<\eta<\eta_{\nu+\mu}$. The above discussion can be repeated. We obtain a fundamental matrix solution  $Y_{\nu+\mu}(z,u)$ with canonical asymptotics $Y_F(z,u)$ in $\widehat{S}_{\nu+\mu}$. 
Again, for  $\eta$ satisfying $\eta_{\nu+2\mu+1}<\eta<\eta_{\nu+2\mu}$ we obtain the analogous result for $Y_{\nu+2\mu}(z,u)$ with canonical asymptotics in $\widehat{S}_{\nu+2\mu}$. 
This can be repeated for every $\nu+h\mu$, $h\in\mathbb{Z}$, obtaining the   fundamental matrix solutions  $Y_{\nu+h\mu}(z,u)$  with canonical asymptotics $Y_F(z,u)$  in $\widehat{S}_{\nu+h\mu}$.  
 So, Points {\bf 1)} and {\bf 2)} of Theorem \ref{20agosto2020-5}  are proved.

\vskip 0.2 cm 

 Stokes matrices are defined by \eqref{26novembre2020-3}. 
Thus, $ \mathbb{S}_{\nu+h\mu}(u)=Y_{\nu+h\mu}(z,u)^{-1}Y_{\nu+(h+1)\mu}(z,u)$ is holomorphic in $\mathbb{D}(u^c)$. 
 Let us consider the relations for $h=0,1$:
  \be
\label{15agosto2020-3}
  Y_{\nu+\mu}(z,u)=Y_\nu(z,u) \mathbb{S}_\nu(u), \quad \quad Y_{\nu+2\mu}(z,u)=Y_{\nu+\mu}(z,u) \mathbb{S}_{\nu+\mu}(u).
\ee
Let $u$ be fixed in a $\tau$-cell, so that  $\Lambda$ has distinct eigenvalues. From   Theorem \ref{20agosto2020-1}   at the fixed $u$ we receive 
     $$
\bigl( \mathbb{S}_\nu(u)\bigr)_{jk}
=
\left\{
\begin{array}{cc}
e^{2\pi i \lp_k}\alpha_k ~c_{jk}^{(\nu)}& ~~~\hbox{ for } j\prec k,
\\
\\
             1   & ~~~\hbox{ for } j =k,
\\
\\
             0   &  ~~~\hbox{ for } j\succ k,
\end{array}
\right.
~~~~~~~~
\bigl( \mathbb{S}_{\nu+\mu}^{-1}(u)\bigr)_{jk}
=
\left\{
\begin{array}{cc}
             0   &  ~~~\hbox{ for } j\prec k,
\\
\\
             1   &~~~ \hbox{ for } j =k,
\\
\\
        -e^{2\pi i (\lp_k-\lp_j)}\alpha_k~c_{jk}^{(\nu)} & 
                                         ~~~\hbox{ for } j\succ k.
\end{array}
\right.
$$
Here, for $j\neq k$ the ordering relation $j\prec k$ $\Longleftrightarrow$ $\Re (z(u_j-u_k))|_{\arg z=\tau}<0$ is well defined for every $u$ in the $\tau$-cell, because no Stokes rays $\Re(z(u_j-u_k))=0$ cross $\arg z=\tau$ as $u$ varies in the $\tau$-cell. 

The relation $j\prec k$ may change to $j\succ k$ when  passing from one $\tau$-cell to another  only for a pair $u_j$, $u_k$ such that $u_j^c=u_k^c$. This is due to the choice of $\epsilon_0$ as in \eqref{30agosto2020-6}. On the other hand,  $c_{jk}^{(\nu)}=0$ whenever $u_j^c=u_k^c$. This means that 
 \eqref{16ottobre2020-5} is true at every fixed $u$ in every $\tau$-cell, with ordering relation $j\prec k$     defined for $j\neq k$ such that $u_j^c\neq u_k^c$, namely $\Re(z(u_j^c-u_k^c))<0$ when $\arg z=\tau$.

Since the $\mathbb{S}_{\nu+h\mu}$ are holomorphic in $\mathbb{D}(u^c)$ and the $c_{jk}^{(\nu)}$ are constant in $\mathbb{D}(u^c)$, we conclude  that Stokes matrices are constant  in $\mathbb{D}(u^c)$ and hence  \eqref{16ottobre2020-5}  holds  in $\mathbb{D}(u^c)$. 
 The vanishing conditions \eqref{16ottobre2020-6} follow from the vanishing  conditions \eqref{31agosto2020-4} for the connection coefficients, plus the fact that we can generate all the  $\mathbb{S}_{\nu+h\mu}$ from the  formula $\mathbb{S}_{\nu+2\mu}=e^{-2\pi i B}\mathbb{S}_\nu e^{2\pi i B}$.

\end{proof}
\section{Non-Uniqueness at $u=u^c$ of the Formal Solution}
\label{16agosto2020-1}

By Laplace transform we prove Corollary  \ref{26nov2018-14} in Background 1, asserting that system \eqref{24nov2018-2} has  unique formal solution if and only if the constant diagonal entries of $A(u)$ satisfy the  {\it partial non-resonance}
$$
\lp_i-\lp_j\not \in \mathbb{Z}\backslash\{0\} \quad \hbox{ for every $ i\neq j$ such that $u_i^c=u_j^c$}
. 
$$
   Otherwise, the Laplace transform will be proved to generate a family of formal solutions $$
\mathring{Y}_F(z)=\Bigl(I+\sum_{l=1}^\infty \mathring{F}_l z^{-l}\Bigr)z^Be^{\Lambda(u^c) z},
$$
 whose {\it coefficients $\mathring{F}_l$ depend on a finite number of arbitrary parameters}. 

\vskip 0.2 cm 
     Due to the strategy of Section \ref{22novembre2020-7}, it will suffice to consider the generic case when all $\lp_1,...,\lp_n ~\not\in\mathbb{Z}$ and $A$ has no integer eigenvalues. Indeed, if this is not the case, the gauge transformation \eqref{2settembre2020-2} relates a formal solution ${}_\gamma Y_F$ to $Y_F$ at any point $u$, through \eqref{2settembre2020-3}, so that the coefficients $F_l$ of a formal expansion do not depend on $\gamma$. We are interested  in these coefficients.

      Consider system \eqref{03} under the assumptions that it is  (strongly) isomonodormic in $\mathbb{D}(u^c)$, so that $(A)_{ij}(u^c)=0$ for $u_i^c=u_j^c$. For simplicity, we order the eigenvalues as in \eqref{2settembre2020-8}-\eqref{2settembre2020-9}. Since $B_1(u)$, ..., $B_n(u)$ are holomorphic at $u^c$, system \eqref{03} at $u=u^c$ is 
      \be
      \label{16agosto2020-3} 
      \frac{d\Psi}{d\lambda} =
      \left(\frac{\sum_{j=1}^{p_1}B_j(u^c)}{\lambda-\lambda_1}+\frac{\sum_{j=p_1+1}^{p_1+p_2}B_j(u^c)}{\lambda-\lambda_2}
      +
      \cdots
      +
      \frac{\sum_{j=p_1+...+p_{s-1}+1}^n B_j(u^c)}{\lambda-\lambda_s}
      \right)\Psi
      \ee
      Let $G^{(\boldsymbol{p}_1)}$ be as in \eqref{19settembre2020-3-new}. The gauge transformation 
      $ 
      \Psi(\lambda)=G^{(\boldsymbol{p}_1)}\widetilde{\Psi}(\lambda)
      $ 
      yields
      \be
      \label{16agosto2020-2} 
      \frac{d\widetilde{\Psi}}{d\lambda}=\left(
      \frac{T^{(\boldsymbol{p}_1)}}{\lambda-\lambda_1}+\sum_{\alpha=2}^s \frac{D_\alpha^{(\boldsymbol{p}_1)}}{\lambda-\lambda_\alpha}
      \right)\widetilde{\Psi},
      \ee
      where
    $$
      T^{(\boldsymbol{p}_1)}:=T^{(1)}+...+T^{(p_1)}=\hbox{diag}(-\lp_1-1,~...,~-\lp_{p_1}-1,~\underbrace{0,~...~,~0}_{n-p_1}).
    $$
    and 
     $D_\alpha^{(\boldsymbol{p}_1)}:= {G^{(\boldsymbol{p}_1)}}^{-1} \cdot \sum_{j=p_1+...+p_{\alpha-1}+1}^{p_1+...+p_{\alpha}} B_j(u^c) \cdot G^{(\boldsymbol{p}_1)}$. 
      The matrix coefficient in system \eqref{16agosto2020-2} has convergent Taylor series at $\lambda=\lambda_1$ 
      $$
       \frac{d\widetilde{\Psi}}{d\lambda}=\frac{1}{\lambda-\lambda_1}\left(
       T^{(\boldsymbol{p}_1)}+
       \sum_{m=1}^\infty \mathfrak{D}_m(\lambda-\lambda_1)^m
       \right)\widetilde{\Psi},
      \quad\quad \mathfrak{D}_m= \sum_{\alpha=2}^s \frac{(-1)^{m+1}}{(\lambda_1-\lambda_\alpha)^m}D_\alpha^{(\boldsymbol{p}_1)}.
      $$ 
      
      We consider $\eta_{\nu+1}<\eta<\eta_\nu$ and $\lambda$ in the plane with branch cuts $\mathcal{L}_\alpha= \mathcal{L}_\alpha(\eta)$ issuing from $\lambda_1,...,\lambda_s$ to infinity in direction $\eta$, as in \eqref{23agosto2020-16}. 
      Close to the Fuchsian singularity $\lambda=\lambda_1$  a fundamental matrix solution to \eqref{16agosto2020-3}  has Levelt  form 
      \be
      \label{17agosto2020-1} 
     \mathring{\Psi}^{(\boldsymbol{p}_1)}(\lambda)=G^{(\boldsymbol{p}_1)}\Bigl(
      I+\sum_{l=1}^\infty \mathfrak{G}_l (\lambda-\lambda_1)^l
      \Bigr)(\lambda-\lambda_1)^{T^{(\boldsymbol{p}_1)}},
      \ee
      where  the matrix entries $ ( \mathfrak{G}_l)_{ij}$, $1\leq i\leq j\leq n$,  are recursively computed by the following formulae (see   Appendix C for an  explanation of  \eqref{17agosto2020-1}, or \cite{Guzz-notes,Wasow}).
      \begin{itemize}
      \item If $T^{(\boldsymbol{p}_1)}_{ii}-T^{(\boldsymbol{p}_1)}_{jj}= l$ positive integer,
       $ ( \mathfrak{G}_l)_{ij} $ is {\it arbitrary}.
       
      \item If $T^{(\boldsymbol{p}_1)}_{ii}-T^{(\boldsymbol{p}_1)}_{jj}\neq l$ (positive integer)
      $$ 
      ( \mathfrak{G}_l)_{ij} =\frac{1}{T^{(\boldsymbol{p}_1)}_{jj}-T^{(\boldsymbol{p}_1)}_{ii}+l}\left(
      \sum_{p=1}^{l-1} \mathfrak{D}_{l-p}\mathfrak{G}_l+\mathfrak{D}_l
      \right)_{ij}      \quad\hbox{(sum is zero for $l=1$)}.$$
     
      \end{itemize}
      
       Since we have assumed that all the $\lp_k$ are not integers,  the only possibility to have $T^{(\boldsymbol{p}_1)}_{ii}-T^{(\boldsymbol{p}_1)}_{jj}= l$ occurs for $1\leq i,j\leq p_1$, precisely 
        \be
        \label{26novembre2020-5}
        T^{(\boldsymbol{p}_1)}_{ii}-T^{(\boldsymbol{p}_1)}_{jj}=\lp_j-\lp_i=l.
        \ee
In this case,  \eqref{17agosto2020-1}  is a family depending on a finite number of   parameters due to the arbitrary $(\mathfrak{G}_l)_{ij}$. 
       Thus,  in the  first $p_1$ columns of a solution of type \eqref{17agosto2020-1} 
      $$
      \vec{\mathring{\Psi}}_j(\lambda~|\nu)=\Bigl(\Gamma(\lp_k+1)\vec{e}_k +\sum_{l=1}^\infty\mathring{b}_l^{(j)} (\lambda-\lambda_1)\Bigr)(\lambda-\lambda_1)^{-\lp_j-1},\quad \quad j=1,...,p_1.
      $$
      the vectors $\mathring{b}_l^{(j)}$ contain a finite number of parameters. The Laplace transform      $$
      \vec{\mathring{Y}}_j (z~|\nu)=\int_{\Gamma_1(\eta)}
  e^{z\lambda}       \vec{\mathring{\Psi}}_j(\lambda~|\nu) d\lambda ,\quad \quad j=1,...,p_1,
  $$
  yields  the first $p_1$ columns of a fundamental matrix solution of \eqref{24nov2018-2}. 
  Repeating the same computations of  Section \ref{20agosto2020-10}, we obtain, for $ j=1,...,p_1$, 
  $$ 
  \vec{\mathring{Y}}_j (z~|\nu) ~z^{-\lp_j}e^{-\lambda_1 z}\sim \vec{e}_j +\sum_{l=1}^\infty \frac{ \mathring{b}_l^{(j)} }{\Gamma(\lp_j+1-l)} \frac{1}{z^l}, \quad z\to \infty \hbox{ in } \mathcal{S}_\nu, 
  $$
   where $\mathcal{S}_\nu$ is given in \eqref{20agosto2020-12}. 
      We repeat the same construction at all $\lambda_1$, ..., $\lambda_s$. This yields a family of fundamental matrix solutions of \eqref{24nov2018-2}  
      $$ 
      \mathring{Y}_\nu(z) = \Bigl[
      \vec{\mathring{Y}}_1 (z~|\nu)~|~\cdots ~|~ \vec{\mathring{Y}}_n (z~|\nu)
      \Bigr],
      $$ 
      depending on a finite number of parameters,  with the behaviour for $z\to \infty$  in  ${S}_\nu$
      $$ 
       \mathring{Y}_\nu(z)\sim \mathring{Y}_F(z)= \Bigl(I+\sum_{l=1}^\infty \mathring{F}_lz^{-l} \Bigr) z^B e^{\Lambda(u^c) z};
       \quad\quad  \mathring{F}_l= \Bigl[
      \vec{\mathring{f}}_1^{~(l)} ~|~\cdots ~|~ \vec{\mathring{f}}_n^{~(l)}
      \Bigr],\quad\quad 
      \vec{\mathring{f}}_j^{~(l)}=\frac{\vec{\mathring{b}}_j^{~(l)}}{\Gamma(\lp_j+1-l)}.
      $$ 
     We conclude that the formal solution is not unique whenever a condition \eqref{26novembre2020-5} occurs. Only one element in the family satisfies 
      $ 
       \mathring{Y}_F(z)=Y_F(z,u^c)$. 
   \bre
   {\rm
   If we choose one formal solution $\mathring{Y}_F(z)$, then the corresponding  $\mathring{Y}_\nu(z)$ with asymptotic expansion $\mathring{Y}_F(z)$  in $\mathcal{S}_\nu$ is unique. For more details on the Stokes phenomenon at $u=u^c$  see \cite{CDG}. 
   }
   \ere

 \section{Appendix A. Non-normalized Schlesinger System}
    \label{23agosto2020-1}
    
    \ble
   The integrability condition $dP=P\wedge P$  of the  Pfaffian system \eqref{11agosto2020-1}  defined on a polydisc $\mathbb{D}(u^0)$ contained in a $\tau$-cell  is the non-normalized Schlesinger system  \eqref{10agosto2020-3}-\eqref{10agosto2020-5}.
    
    \ele
    \begin{proof}
For  every $i\in \{1,...,n\}$, the  Pfaffian system \eqref{11agosto2020-1}   can be rewritten as
$$
 P= \left(\frac{B_i}{\lambda-u_i}+\sum_{j\neq i}\frac{B_j}{\lambda-u_j}\right)d(\lambda-u_i)
 +\sum_{j\neq i} \left(\gamma_j-\frac{B_j}{\lambda-u_j}\right) d(u_j-u_i) +\sum_{j=1}^n \gamma_j(u) d\lambda.
 $$
 We study $\lambda-u_i\to 0$, while  $u_j-u_i\neq 0$ for $j\neq i$ in $\mathbb{D}(u^0)$. In  new variables
 $$ 
 \lambda=\lambda, \quad\quad y_i=\lambda-u_i, \quad y_j=u_j-u_i,\quad j\neq i,
 $$
 $P$ is rewritten in the following way (which defines the matrices $\mathcal{A}_j(y)$)
 \begin{align*}
P=& \left(\frac{B_i}{y_i}+\sum_{j\neq i}\frac{B_j}{y_i-y_j}\right)dy_i
 +\sum_{j\neq i} \left(\gamma_j-\frac{B_j}{y_i-y_j}\right) dy_j+\sum_{j=1}^n \gamma_j(y) d\lambda
 \\
 &=: \mathcal{A}_i(y) dy_i +\sum_{j\neq i} \mathcal{A}_j(y) dy_j +\sum_{j=1}^n \gamma_j(y) d\lambda. 
 \end{align*}
  The only singular term at $y_i=0$  is ${B_i/y_i}$ in $\mathcal{A}_i(y)$.  The components  relative to   $dy_1,...,dy_n$ of $dP=P\wedge P$ are
\be
\label{23agosto2020-2}
 \frac{\partial {\mathcal{A}}_l}{\partial y_k}+{\mathcal{A}}_l{\mathcal{A}}_k=\frac{\partial {\mathcal{A}}_k}{\partial y_l}+{\mathcal{A}}_k{\mathcal{A}}_l,\quad k\neq l,
\ee
For $k\neq i$ and $l=i$, from   \eqref{23agosto2020-2} we receive
$$ 
\frac{\partial}{\partial y_k} \left(\frac{B_i}{y_i}+\hbox{\rm reg}(y_i)\right)
+
\left(\frac{B_i}{y_i}+\hbox{\rm reg}(y_i)\right){\mathcal{A}}_k=
\frac{\partial {\mathcal{A}}_k}{\partial y_i} + {\mathcal{A}}_k\left( \frac{B_i}{y_i}+\hbox{\rm reg}(y_i)\right),
$$
where $\hbox{\rm reg}(y_i)$ stands for an analytic term at $y_i=0$. We expand the above in Taylor series at $y_i=0$.   The singular term (the residue at $y_i=0$) is 
\be
\label{23agosto2020-4} 
\frac{\partial B_i}{\partial y_k}=\bigl[\mathcal{A}_k|_{y_i=0},B_i\bigr]= \frac{[B_k,B_i]}{u_k-u_i}+[\gamma_k,B_i],\quad k\neq i.
 \ee
    The above gives the non-normalized Schlesinger equations \eqref{10agosto2020-4}-\eqref{10agosto2020-5}, because 
 \begin{align}
 \label{23agosto2020-5}
 &\frac{\partial B_i}{\partial y_k}=\frac{\partial B_i}{\partial (u_k-u_i)}= \frac{\partial u_k}{\partial (u_k-u_i)}\frac{\partial B_i }{\partial u_k}= \frac{\partial B_i}{\partial u_k},
 \\
 \noalign{\medskip}
 \label{23agosto2020-6}
 & \frac{\partial B_i}{\partial u_i}= \sum_{k\neq i} \frac{\partial(u_k-u_i)}{\partial u_i} \frac{\partial B_i}{\partial (u_k-u_i)}=
 -\sum_{k\neq i} \frac{\partial B_i}{\partial u_k}\quad\Longrightarrow \quad \sum_{k=1}^n \frac{\partial B_i}{\partial u_k}=0.
 \end{align}
 If we write the components of  $dP=P\wedge P$ referring to   $dy_l$ ad $d\lambda$, and we substitute into them  \eqref{23agosto2020-5}-\eqref{23agosto2020-6}, we receive  \eqref{10agosto2020-3}, namely 
 $ 
 \partial_l\gamma_k-\partial_k\gamma_l=\gamma_l\gamma_k-\gamma_k\gamma_l
 $. 
 \end{proof}
 
 \bcr
 \label{18novembre2020-2}
 A solution $B_i(u)$, $i=1,...,n$,   of \eqref{10agosto2020-3}-\eqref{10agosto2020-5} is holomorphically similar  to a constant Jordan form on $\mathbb{D}(u^0)$. The similarity is realized by a fundamental matrix solution $G^{(i)}(u)$  of the Pfaffian system \eqref{22novembre2020-6} below.  
 
 \ecr
 
 \begin{proof}
 We must show that there exists a holomorphically invertible $G^{(i)}(u)$  on $\mathbb{D}(u^0)$ such that 
 $(G^{(i)})^{-1} B_i G^{(i)}$ is a constant Jordan form. 
The conditions \eqref{23agosto2020-2} for $k,l\neq i$ can be evaluated at $y_i=0$, and  become
 $$ 
 \frac{\partial {\mathcal{A}}_l |_{y_i=0}}{\partial y_k}+{\mathcal{A}}_l |_{y_i=0}{\mathcal{A}}_k |_{y_i=0}=\frac{\partial {\mathcal{A}}_k|_{y_i=0}}{\partial y_l}+{\mathcal{A}}_k|_{y_i=0}{\mathcal{A}}_l|_{y_i=0},\quad \quad k\neq i, ~l\neq i ,~ k\neq l.
$$    
Hence, the following Pfaffian system is Frobenius integrable
\be 
\label{23agosto2020-7}
\frac{\partial G}{\partial y_k}= \mathcal{A}_k|_{y_i=0} ~G~\equiv \left(\frac{B_k}{u_k-u_i}+\gamma_k\right)G
,\quad k\neq i.
\ee
Using the chain rule as in \eqref{23agosto2020-5}, we receive \eqref{18novembre2020-1} 
\be
\label{22novembre2020-6}
\frac{\partial G}{\partial u_k} =\left(\frac{B_k}{u_k-u_i}+\gamma_k\right)G,\quad k\neq i, \quad 
\quad\quad 
 \frac{\partial G}{\partial u_i}=-\sum_{k\neq i} \left(\frac{B_k}{u_k-u_i}+\gamma_k\right) G
 \ee
Notice that  for both $\varphi(u)=B_i(u)$ and $\varphi(u)=G(u)$ we have 
\be
\label{23agosto2020-8} 
\sum_{k=1}^n \frac{\partial \varphi}{\partial u_k}=0\quad\Longrightarrow \quad \varphi(u)=\varphi(u_1-u_i,~\dots~,u_n-u_i).
\ee
We can take a solution $G(u)$ which holomorphically reduces  $B_i$ to Jordan form.  Indeed 
\begin{align*}
\hbox{ for $k\neq i$,  
 } \quad \frac{\partial}{\partial y_k} (G^{-1}B_iG)
&=
-G^{-1}\frac{\partial G}{\partial y_k} G^{-1}B_iG+G^{-1}\frac{\partial B_i}{\partial y_k} G + G^{-1} B_i \frac{\partial G}{\partial y_k} 
\\
& 
\underset{\eqref{23agosto2020-4},\eqref{23agosto2020-7}}=
-G^{-1}\mathcal{A}_k|_{y_i=0} B_i G 
+ G^{-1}\bigl[\mathcal{A}_k|_{y_i=0},~B_i\bigr]G
+G^{-1}B_i \mathcal{A}_k|_{y_i=0} G
\\
&
\quad\quad =0.
\end{align*}
Therefore, keeping into account \eqref{23agosto2020-8}, we see that $\mathcal{B}_i:=G^{-1}(u)B_i(u)G(u))$ is independent of $u$. Thus, there exists a constant matrix 
$\mathcal{G}$ such that $\mathcal{G}^{-1}\mathcal{B}_i \mathcal{G}$ is a constant Jordan form, and $G^{(i)}(u):=G(u) \mathcal{G}$ realises the holomorphic "Jordanization" .  The above arguments are standard, see for example \cite{Haraoka}.

\end{proof}

If the $B_i(u)$ are holomorphic on $\mathbb{D}(u^c)$ and the vanishing conditions \eqref{12agosto2020-1} hold, the coefficients of the Pfaffian system \eqref{22agosto2020-7} are holomorphic on $\mathbb{D}(u^c)$, so that $G^{(i)}(u)$ extends holomorphically there, and Corollary \ref{18novembre2020-2} holds on $\mathbb{D}(u^c)$.


\section{Appendix B. Proof of Proposition \ref{11agosto2020-5}}

  According to Theorem \ref{9agosto2020-1}, system \eqref{24nov2018-1} is strongly isomonodromic in   $\mathbb{D}(u^0)$ contained in a $\tau$-cell of $\mathbb{D}(u^c)$ if and only if \eqref{10agosto2020-1} holds. 
In this case $G^{(0)}$ in \eqref{25nov2018-7}  holomorphically reduces $A(u)$ to constant Jordan form and  satisfies
\be
\label{10agosto2020-2}
dG^{(0)}=\sum_{j=1}^n \omega_j(u) du_j ~G^{(0)}.
\ee

\begin{proof}[Proof of Proposition \ref{11agosto2020-5}]

  Suppose that  \eqref{24nov2018-1} is strongly isomonodromic, so that  \eqref{10agosto2020-1} holds.      Let $\mathcal{A}:=-A-I$, so that $E_k\mathcal{A}=B_k$, and \eqref{10agosto2020-1} are rewritten as $\partial_i \mathcal{A}=[\omega_i(u),\mathcal{A}]$. We multiply these equations to the left by $E_k$, with $k\neq i$. We receive 
    $$
    E_k\partial_i \mathcal{A} =E_k[\omega_i(u),\mathcal{A}].
    $$ 
    The l.h.s. is $ 
    E_k\partial_i \mathcal{A}=\partial_i B_k$. 
        The r.h.s.  is 
    $$ 
    E_k[ \omega_i,\mathcal{A}]=E_k  \omega_i\mathcal{A}-E_k\mathcal{A}  \omega_i=E_k \omega_i\mathcal{A}-B_k
     \omega_i
= \bigl(E_k \omega_i\mathcal{A}- \omega_i B_k\bigr)+[\omega_i,B_k].
$$
In conclusion 
$$ 
\partial_i B_k=\bigl(E_k \omega_i\mathcal{A}- \omega_i B_k\bigr)+[\omega_i,B_k],\quad i\neq k.
$$
The only terms we need to evaluate are 
   $$
   E_k \omega_i\mathcal{A}- \omega_i B_k=E_k[F_1,E_i]\mathcal{A}-[F_1,E_i]B_k=
   $$ 
   $$
   =E_kF_1E_i\mathcal{A}+E_iF_1B_k= E_kF_1E_i B_i+E_iF_1E_k B_k.
   $$
   In the second line we have used  $E_iE_k=E_kE_i=E_iB_k=0$, for $i\neq k$, and $E_i^2=E_i$. Now, observe that $E_kF_1E_i $ has zero entries, except for the entry $(k,i)$, which is $(F_1)_{ki}= (A)_{ki}/(u_i-u_k)$. This implies that 
   $$ 
   E_kF_1E_i B_i+E_iF_1E_k B_k= \frac{[B_i,B_k]}{u_i-u_k}. 
   $$ 
   In conclusion, we have proved that \eqref{10agosto2020-1} implies \eqref{10agosto2020-4}. On the other hand \eqref{10agosto2020-4}-\eqref{10agosto2020-5} are equivalent to the system given by  \eqref{10agosto2020-4} and the equations
   $$ 
   \partial_i \sum_kB_k= [ \omega_i,\sum_k B_k],\quad i=1,...,n.
   $$ 
   which are exactly \eqref{10agosto2020-1} if $B_k=E_k\mathcal{A}$.   
   Finally, notice that  \eqref{10agosto2020-3}, here with $\gamma_j=\omega_j$, is the integrability condition on $\mathbb{D}(u^0)$ of $dG=\sum_{j=1}^n \omega_j(u) du_j ~G$. 
On the other hand, it is a computation to see that \eqref{10agosto2020-1} implies  the the same  conditions.

   \vskip 0.2 cm

Conversely,   let system \eqref{03} be strongly isomonodromic, so that the integrability conditions \eqref{10agosto2020-3}-\eqref{10agosto2020-5} hold. 
  Firstly, we show that  \eqref{10agosto2020-4}-\eqref{10agosto2020-5} imply a Pfaffian system for $A$ of type  \eqref{10agosto2020-1}. To this end,  we sum \eqref{10agosto2020-4} and \eqref{10agosto2020-5}: 
   $$
   \sum_{k=1}^n \partial_i B_k= \sum_{k\neq i} \frac{[B_i,B_k]}{u_i-u_k}- \sum_{k\neq i} \frac{[B_i,B_k]}{u_i-u_k}+[\gamma_i, \sum_{k=1}^n B_k]=[\gamma_i, \sum_{k=1}^n B_k].
   $$ 
   Using $B_k=-E_k(A+I)$ and $\sum_k E_k=I$, the above becomes
 \be
 \label{18agosto2020-5}
   \partial_i A=[\gamma_i,A],\quad i=1,...,n.
   \ee
    Since $\gamma_1,...,\gamma_n$ satisfy  \eqref{10agosto2020-3}, it is directly verified  that \eqref{18agosto2020-5} is  Frobenius integrable.    
 Secondly,  we must show that we can choose
 $$
 \gamma_j:=\omega_j= \left(
 \frac{A_{ab}(\delta_{aj}-\delta_{bj})}{u_a-u_b}
 \right)_{a,b=1}^n
  \quad \hbox{ as in \eqref{22agosto2020-6}.}
 $$
 Substituting this choice into  \eqref{10agosto2020-3}, we see that if \eqref{18agosto2020-5} holds, in the form   $\partial_i A=[\omega_i,A]$, then \eqref{10agosto2020-3} are satisfied.\footnote{This is exactly what has been said before:  \eqref{10agosto2020-1} implies by computation   the integrability conditions of \eqref{10agosto2020-2},  namely  exactly equations \eqref{10agosto2020-3} with $\gamma_j=\omega_j$.}
  Now,  since   \eqref{18agosto2020-5} follows from \eqref{10agosto2020-4}-\eqref{10agosto2020-5} with  matrices  $B_k=-E_k(A+I)$, we conclude that  \eqref{10agosto2020-4}-\eqref{10agosto2020-5}  and the choice $\gamma_j=\omega_j$  guarantee that both \eqref{10agosto2020-3} and \eqref{10agosto2020-1}  are satisfied. 

\end{proof}
 
 The Schlesinger system can be used to show that  there is a  fundamental matrix solution $G^{(0)}$ of 
  \eqref{10agosto2020-2}  that   holomorphically  reduces $A$ to constant Jordan form on $\mathbb{D}(u^0)$. Since \eqref{10agosto2020-3} is the integrability condition on $\mathbb{D}(u^0)$ of the linear Pfaffian system
 \be
 \label{5aprile2021-1} dG=\sum_{j=1}^n \omega_j(u) du_j ~G.
 \ee
the latter admits holomorphic fundamental matrix solutions in $\mathbb{D}(u^0)$. Let $G(u)$ be one of them and define 
  \be
  \label{18agosto2020-1}
  \widehat{B}_k:=G(u)^{-1} B_k G(u).
  \ee
  By direct computation,  using \eqref{5aprile2021-1} and its integrability \eqref{10agosto2020-3},  it is verified that \eqref{10agosto2020-4}-\eqref{10agosto2020-5} (with $\gamma_j=\omega_j$) are equivalent to the normalized Schlesinger equations for the matrices  $\widehat{B}_k$,
$$
\partial_i\widehat{B}_k=\frac{[\widehat{B}_i,\widehat{B}_k]}{u_i-u_k}, \quad i\neq k;\quad\quad
   \partial_i\widehat{B}_i= -\sum_{k\neq i} \frac{[\widehat{B}_i,\widehat{B}_k]}{u_i-u_k}.
 $$
   The above equations imply that 
   $$
   \forall~i=1,...,n,\quad\quad \partial_i \widehat{B}_\infty =0,
   \quad \quad 
   \hbox{  where $\widehat{B}_\infty:= -\sum_{k=1}^n \widehat{B}_k$}
$$
  Being  $\widehat{B}_\infty$ constant, it can be put in constant Jordan form by a constant invertible matrix $P$, say $-\mathcal{J}=P^{-1} \widehat{B}_\infty P$. Since also $G(u)P$ solves \eqref{5aprile2021-1},  we can choose from the beginning $G(u)$  such that 
 \be
 \label{8novembre2020-3}
  G^{-1}(u) \left(\sum_{k=1}^n B_k(u)\right) G(u)=\mathcal{J} \quad\hbox{ constant Jordan form}. 
  \ee
 Now, observe that $\sum_{k=1}^nE_k=I$, so that 
 $$\sum_{k=1}^n B_k=-\sum_{k=1}^n E_k (A+I)= -A-I.
 $$ 
  Thus, $G(u)$ also puts $A$ in constant Jordan form, so that\footnote{Up to the freedom $G\mapsto G G_*$ where $G_*$ commutes with the Jordan form.}
  $$ G(u)=G^{(0)}(u),\quad\quad \hbox{ where $G^{(0)}$ is in \eqref{25nov2018-7}}.
  $$ 
    In particular, $G^{(0)}$ satisfies  \eqref{10agosto2020-2}.

\vskip 0.2 cm 
 The second part of the statement of Proposition \ref{6aprile2021-3} (Prop. 19.2 of \cite{CDG}) is now easily proved. Indeed, if $A(u)=G^{(0)}(u) J (G^{(0)})^{-1}$ holomorphically on $\mathbb{D}(u^c)$, where $J$ is Jordan,  then  $G^{(0)}$ satisfies  \eqref{10agosto2020-2} on $\mathbb{D}(u^0)$ (and $J$ is constant).  Since $G^{(0)}(u)$ is holomorphic on $\mathbb{D}(u^c)$, the $\omega_j$ must be as well, so that the vanishing conditions  \eqref{26nov2018-7} must hold. Conversely, if $A$ is holomorphic on on $\mathbb{D}(u^0)$  and satisfies the vanishing conditions  \eqref{26nov2018-7}  (or, more weakly, if $dA=\sum_j[\omega_j ,A]du_j$ on $\mathbb{D}\backslash \Delta$, which automatically implies \eqref{26nov2018-7} -- see Remark \ref{6aprile2021-2}), then $dG=\sum_j\omega_j du_j ~G$ is integrable with holomorphic coefficients on $\mathbb{D}(u^0)$,  and admits a fundamental matrix solution that can be chosen  so that $ (G^{(0)})^{-1}AG^{(0)}(u)=J$ (the proof is as done before on $\mathbb{D}(u^0)$).

\section{Appendix C. The normal form \eqref{17agosto2020-1}}
We prove the expression \eqref{17agosto2020-1} of Section \ref{16agosto2020-1}, where it was sufficient to only consider the generic case of all $\lp_1,...,\lp_n ~\not\in\mathbb{Z}$ and  no integer eigenvalues of $A$. 
    A fundamental matrix solution in Levelt form at $\lambda=\lambda_1$ for system \eqref{16agosto2020-3}   has structure 
       \be
       \label{27novembre2020-1} 
     \mathring{\Psi}(\lambda)=G^{(\boldsymbol{p}_1)}\Bigl(
      I+\sum_{l=1}^\infty \mathfrak{G}_l (\lambda-\lambda_1)^l
      \Bigr)(\lambda-\lambda_1)^{T^{(\boldsymbol{p}_1)}}(\lambda-\lambda_1)^R,
      \ee
      with 
      $$  R=R_1+R_2+\dots R_{\kappa},\quad \quad  \kappa:=\max\{T^{(\boldsymbol{p}_1)}_{ii}-T^{(\boldsymbol{p}_1)}_{jj} \hbox{ integer}\}.
      $$
      where $R$ is a nilpotent matrix with $R_{ij}\neq 0$  only if  $T^{(\boldsymbol{p}_1)}_{ii}-T^{(\boldsymbol{p}_1)}_{jj}$ is a  positive integer.   We prove that  $R=0$ in our case.  The  formulae  for $ ( \mathfrak{G}_l)_{ij}$ and $(R_l)_{ij}$ are obtained recursively by substituting the series into the differential system,  and are as follows. 
      \begin{itemize}
      
       \item If $T^{(\boldsymbol{p}_1)}_{ii}-T^{(\boldsymbol{p}_1)}_{jj}= l$ (positive integer),
       $ ( \mathfrak{G}_l)_{ij} $ is arbitrary, and 
       $$
       (R_l)_{ij}=\left(
      \sum_{p=1}^{l-1} (\mathfrak{D}_{l-p}\mathfrak{G}_l-\mathfrak{G}_l{R}_{l-p})+\mathfrak{D}_l\right)_{ij},$$
      \item If $T^{(\boldsymbol{p}_1)}_{ii}-T^{(\boldsymbol{p}_1)}_{jj}\neq l$ (positive integer)
      $$ 
      ( \mathfrak{G}_l)_{ij} =\frac{1}{T^{(\boldsymbol{p}_1)}_{jj}-T^{(\boldsymbol{p}_1)}_{ii}+l}\left(
      \sum_{p=1}^{l-1} (\mathfrak{D}_{l-p}\mathfrak{G}_l-\mathfrak{G}_l{R}_{l-p})+\mathfrak{D}_l
      \right)_{ij}      $$

      \end{itemize}
     The claim that $R=0$ follows from two facts. First, if we evaluate at $u=u^c$  the isomonodromic  fundamental matrix solution \eqref{19settembre2020-3-new}, we receive   a fundamental matrix solution of \eqref{03} at $u=u^c$, 
    \be
    \label{22novembre2020-8}
     \Psi^{(\boldsymbol{p}_1)}(\lambda,u^c)=G^{(\boldsymbol{p}_1)}\cdot U^{(\boldsymbol{p}_1)}(\lambda,u^c)\cdot  (\lambda-\lambda_1)^{T^{(\boldsymbol{p}_1)}},
     \ee
      which has   $R=0$, because in the generic case here considered all  $R^{(j)}=0$ in \eqref{19settembre2020-3-new}.  
 The expression \eqref{22novembre2020-8}  belongs to the class of solutions \eqref{27novembre2020-1}.  
      
          The second fact is  that other solutions in the class \eqref{27novembre2020-1} may have different matrix-exponents (see \cite{Guzz-notes} and \cite{CDG}; see also \cite{Dub2,CDG1} for the case of Frobenius manifolds), but  if  $R$ corresponds to one solution, all the other solutions in the class can only have exponent   
          \be
          \label{2settembre2020-5}
         \widetilde{R}= \mathcal{D}^{-1} R \mathcal{D},
          \ee
                   where $\mathcal{D}$ is an invertible matrix explained below. Now, since $R=0$ in \eqref{22novembre2020-8}, then by \eqref{2settembre2020-5} all the other $\widetilde{R}=0$. This  proves that \eqref{17agosto2020-1} is the correct form. 
    
   Finally, we explain \eqref{2settembre2020-5}. System \eqref{03} at $u=u^c$ is holomorphically equivalent to   "Birkhoff-normal forms"
                   $$
                   \frac{d\Psi}{d\lambda}=\left(\frac{T^{(\boldsymbol{p}_1)}}{\lambda-\lambda_1}+\sum_{l=1}^\kappa R_l(\lambda-\lambda_1)^l\right) \Psi 
                   \quad
                    \hbox{ and }
                    \quad
              \frac{d\widetilde{\Psi}}{d\lambda}=\left(\frac{T^{(\boldsymbol{p}_1)}}{\lambda-\lambda_1}+\sum_{l=1}^\kappa \widetilde{R}_l(\lambda-\lambda_1)^l\right) \widetilde{\Psi},
                                      $$
                  which are related to each other by a 
                    gauge transformations $\Psi=\mathcal{D}(\lambda) \widetilde{\Psi}$, with $\mathcal{D}(\lambda)=\mathcal{D}_0(I+\mathcal{D}_0(\lambda-\lambda_1)+\dots + \mathcal{D}_\kappa(\lambda-\lambda_1)^\kappa)$, where $\hbox{det}(\mathcal{D}_0)\neq 0$ and $[\mathcal{D}_0,T^{(\boldsymbol{p}_1)}]=0$. Then, $\mathcal{D}:=\mathcal{D}_0(I+\mathcal{D}_0+\dots + \mathcal{D}_\kappa)$.

     \bre {\rm In our case, the equations $R_l=0$, $l=1,2,...,\kappa$ are conditions on the entries of $A(u^c)$. The above discussion shows that,   in the isomonodromic case, such conditions turn out to be automatically satisfied with the only vanishing assumption  $(A(u^c))_{ab}=0$ for $u_a^c=u_b^c$. 
     These conditions are equivalent to the conditions  (4.24)-(4.25) of Proposition 4.2 in \cite{CDG}, and probably more conveninent. We will not enter into the tedious verification of the equivalence. 
    }
    \ere
\begin{small}

\end{small}
\end{document}